\documentclass[12pt,letterpaper,leqno]{amsart}
\usepackage{microtype}
\usepackage{amssymb}
\usepackage{amsfonts}
\usepackage{geometry}
\usepackage{graphicx}

\newtheorem{theorem}{Theorem}
\theoremstyle{plain}

\newtheorem{corollary}{Corollary}

\newtheorem{definition}{Definition}

\newtheorem{lemma}{Lemma}
\newtheorem{notation}{Notation}
\newtheorem{problem}{Problem}

\numberwithin{equation}{section}
\numberwithin{theorem}{section}  
\numberwithin{proposition}{section}  
\numberwithin{lemma}{section}  
\numberwithin{corollary}{section}  
\input{tcilatex}
\begin{document}
\title[1D Focusing NLS from 1D N-Body Focusing Quantum System]{Focusing
Quantum Many-body Dynamics: \\
The Rigorous Derivation of the 1D Focusing Cubic Nonlinear Schr\"{o}dinger
Equation }
\author{Xuwen Chen}
\address{Department of Mathematics, Brown University, 151 Thayer Street,
Providence, RI 02912}
\email{chenxuwen@math.brown.edu}
\urladdr{http://www.math.brown.edu/\symbol{126}chenxuwen/}
\author{Justin Holmer}
\address{Department of Mathematics, Brown University, 151 Thayer Street,
Providence, RI 02912}
\email{holmer@math.brown.edu}
\urladdr{http://www.math.brown.edu/\symbol{126}holmer/}

\begin{abstract}
We consider the dynamics of $N$ bosons in one dimension. We assume that the
pair interaction is attractive and given by $N^{\beta -1}V(N^{\beta }\cdot )$
where $\int V\leqslant 0.$ We develop new techniques in treating the $N-$%
body Hamiltonian so that we overcome the difficulties generated by the
attractive interaction and establish new energy estimates. We also prove the
optimal 1D collapsing estimate which reduces the regularity requirement in
the uniqueness argument by half a derivative. We derive rigorously the one
dimensional focusing cubic NLS with a quadratic trap as the $N\rightarrow
\infty $ limit of the $N$-body dynamic and hence justify the mean-field
limit and prove the propagation of chaos for the focusing quantum many-body
system.
\end{abstract}

\maketitle
\tableofcontents

\section{Introduction}

In 1925, Einstein predicted that, at low temperatures, non-interacting
bosons in a gas could all reside in the same quantum state. This peculiar
gaseous state in trapped interacting atomic clouds, a Bose-Einstein
condensate (BEC), was produced in the laboratory for the first time in 1995
using the laser-cooling methods \cite{Anderson, Davis}. E. A. Cornell, W.
Ketterle, and C. E. Wieman were awarded the 2001 Nobel Prize in physics for
observing BEC. Many similar successful experiments \cite{Cornish, Ketterle,
Stamper} were performed later. These condensates exhibit quantum phenomena
on a large scale, and investigating them has become one of the most active
areas of contemporary research.

Let $t\in \mathbb{R}$ be the time variable and $\mathbf{x}_{N}=\left(
x_{1},x_{2},...,x_{N}\right) \in \mathbb{R}^{nN}$ be the position vector of $%
N$ particles in $\mathbb{R}^{n}$. Then BEC naively means that the $N$-body
wave function $\psi _{N}(t,\mathbf{x}_{N})$ satisfies 
\begin{equation*}
\psi _{N}(t,\mathbf{x}_{N})\sim \dprod\limits_{j=1}^{N}\phi (t,x_{j})
\end{equation*}%
up to a phase factor solely depending on $t$, for some one particle state $%
\phi .$ In other words, every particle is in the same quantum state.
Equivalently, there is the Penrose-Onsager formulation of BEC: if we define $%
\gamma _{N}^{(k)}$ to be the $k$-particle marginal densities associated with 
$\psi _{N}$ by 
\begin{equation}
\gamma _{N\,}^{(k)}(t,\mathbf{x}_{k};\mathbf{x}_{k}^{\prime })=\int \psi
_{N}(t,\mathbf{x}_{k},\mathbf{x}_{N-k})\overline{\psi _{N}}(t,\mathbf{x}%
_{k}^{\prime },\mathbf{x}_{N-k})d\mathbf{x}_{N-k},\quad \mathbf{x}_{k},%
\mathbf{x}_{k}^{\prime }\in \mathbb{R}^{nk}  \label{E:marginal}
\end{equation}%
then, equivalently, BEC means 
\begin{equation}
\gamma _{N}^{(k)}(t,\mathbf{x}_{k};\mathbf{x}_{k}^{\prime })\sim
\dprod\limits_{j=1}^{k}\phi (t,x_{j})\bar{\phi}(t,x_{j}^{\prime }).
\label{formula:BEC state}
\end{equation}%
It is widely believed that the one particle state $\phi $ in (\ref%
{formula:BEC state}), also called the condensate wave function since it
describes the whole condensate, satisfies the cubic nonlinear Schr\"{o}%
dinger equation (NLS)%
\begin{equation*}
i\partial _{t}\phi =L\phi +\mu \left\vert \phi \right\vert ^{2}\phi ,
\end{equation*}%
where $L$ is the Laplacian $-\triangle $ or the Hermite operator $-\triangle
+\omega ^{2}\left\vert x\right\vert ^{2}$. Such a belief is one of the
motivations for studying the cubic NLS. Here, the nonlinear term $\mu
\left\vert \phi \right\vert ^{2}\phi $ represents a mean-field approximation
of the pair interactions between the particles: a repelling interaction
gives a positive $\mu $ while an attractive interaction yields a $\mu <0$.
Gross and Pitaevskii proposed such a description of the many-body effect.
Naturally, the validity of the cubic NLS needs to be established rigorously
from the many body system which it is supposed to characterize because it is
a phenomenological mean-field type equation.

In a series of works \cite{Lieb2, AGT, E-E-S-Y1, E-S-Y1,E-S-Y2,E-S-Y5,
E-S-Y3,TChenAndNP,ChenAnisotropic,TChenAndNPSpace-Time, Chen3DDerivation,
SchleinNew, C-H3Dto2D, GM1}, it has been proven rigorously that, for a
repelling interaction potential with suitable assumptions, relation %
\eqref{formula:BEC state} holds, moreover, the one-particle state $\phi $
satisfies the defocusing cubic NLS ($\mu >0$).

It is then natural to wonder, whether BEC happens (whether relation %
\eqref{formula:BEC state} holds) when the interaction potential is
attractive, and whether the condensate wave function $\phi $ satisfies a
focusing cubic NLS ($\mu <0$) if relation \eqref{formula:BEC
state} does hold. In contemporary experiments, both positive \cite%
{Khaykovich,Streker} and negative \cite{Cornish} results exist. To present
the mathematical interpretations of the experiments, we investigate the
procedure of laboratory experiments of BEC subject to attractive
interactions according to \cite{Cornish,Khaykovich,Streker}.

\begin{itemize}
\item[Step A.] Confine a large number of bosons, whose interactions are
originally \emph{repelling}, inside a trap. Reduce the temperature of the
system so that the many-body system reaches its ground state. It is expected
that this ground state is a BEC state / factorized state. This step
corresponds to the following mathematical problem.

\begin{problem}
\label{Problem:Lieb}Show that if $\psi _{N,0}$ is the ground state of the
N-body Hamiltonian $H_{N,0}$ defined by%
\begin{equation}
H_{N,0}=\sum_{j=1}^{N}\left( -\frac{1}{2}\triangle _{x_{j}}+\frac{\omega
_{0}^{2}}{2}\left\vert x_{j}\right\vert ^{2}\right) +\frac{1}{N}%
\sum_{1\leqslant i<j\leqslant N}N^{n\beta }V_{0}\left( N^{\beta }\left(
x_{i}-x_{j}\right) \right)  \label{Hamiltonian:Step A}
\end{equation}%
where $V_{0}\geqslant 0$, then the marginal densities $\left\{ \gamma
_{N,0}^{(k)}\right\} $ associated with $\psi _{N,0}$, defined in (\ref%
{E:marginal}), satisfy relation \eqref{formula:BEC state}$.$
\end{problem}

Here, the factor $1/N$ is to make sure that the interactions are
proportional to the number of particles, the pair interaction $N^{n\beta
}V_{0}(N^{\beta }\cdot )$ is an approximation to the Dirac $\delta $
function so that it matches the Gross-Pitaevskii description of BEC that the
many-body effect should be modeled by a strong on-site self-interaction, and
the quadratic potential $\omega _{0}^{2}\left\vert x\right\vert ^{2}$
represents the trapping since \cite{Cornish,Khaykovich,Streker} and many
other experiments of BEC use the harmonic trap and measure the strength of
the trap with $\omega _{0}$. This step is exactly the same as the
preparation of experiments with repelling interactions and satisfactory
answers to Problem \ref{Problem:Lieb} have been given in \cite{Lieb2}.

\item[Step B.] Strengthen the trap (increase $\omega _{0}$) to make the
interaction attractive and observe the evolution of the many-body system.
This technique which continuously controls the sign and the size of the
interaction in a certain range is called the Feshbach resonance.\footnote{%
See \cite[Fig.1]{Cornish}, \cite[Fig.2]{Khaykovich}, or \cite[Fig.1]{Streker}
for graphs of the relation between $\omega $ and $V$.} The system is then
time dependent. In order to observe BEC, the factorized structure obtained
in Step A must be preserved in time. Assuming this to be the case, we then
reset the time so that $t=0$ represents the point at which this Feshbach
resonance phase is complete. The subsequent evolution should then be
governed by a focusing time-dependent $N$-body Schr\"odinger equation with
an attractive pair interaction $V$ subject to an asymptotically factorized
initial datum.  Moreover, the confining strength is different from Step A,
and we denote it by $\omega$. A mathematically precise statement is the
following:

\begin{problem}
\label{Problem:Ours}Let $\psi _{N}\left( t,\mathbf{x}_{N}\right) $ be the
solution to the $N-body$ Schr\"{o}dinger equation%
\begin{equation}
i\partial _{t}\psi _{N}=\sum_{j=1}^{N}\left( -\frac{1}{2}\triangle _{x_{j}}+%
\frac{\omega ^{2}}{2}\left\vert x_{j}\right\vert ^{2}\right) \psi _{N}+\frac{%
1}{N}\sum_{1\leqslant i<j\leqslant N}N^{n\beta }V\left( N^{\beta }\left(
x_{i}-x_{j}\right) \right) \psi _{N}  \label{Hamiltonian:Step B}
\end{equation}%
where $V\leqslant 0,$ with $\psi _{N,0}$ from Step A as initial datum. Prove
that the marginal densities $\left\{ \gamma _{N}^{(k)}(t)\right\} $
associated with $\psi _{N}\left( t,\mathbf{x}_{N}\right) $ satisfies
relation \eqref{formula:BEC state}.\footnote{%
Since $\omega \neq \omega _{0}$, $V\neq V_{0}$, one could not expect that $%
\psi _{N,0},$ the ground state of (\ref{Hamiltonian:Step A}), is close to
the ground state of (\ref{Hamiltonian:Step B}).}
\end{problem}
\end{itemize}

In the experiment by Cornell and Wieman et.al \cite{Cornish}, once the
interaction is tuned attractive, the condensate suddenly shrinks to below
the resolution limit, then after $\sim 5ms$, the many-body system blows up.
That is, there is no BEC once the interaction becomes attractive. Moreover,
there is no condensate wave function due to the absence of the condensate.
Whence, the current\ NLS theory, which is about the condensate wave function
when there is a condensate, cannot explain this $5ms$ of time or the blow
up. This is currently an open problem in the study of quantum many systems.

In \cite{Khaykovich,Streker}, the particles are confined in a strongly
anisotropic cigar-shape trap. That is, the confinement is very strong in two
spatial directions to simulate a 1D system. In this case, the experiment is
a success in the sense that one obtains a persistent BEC after the
interaction is switched to attractive. Moreover, a soliton is observed in 
\cite{Khaykovich} and a soliton train is observed in \cite{Streker}. The
solitons in \cite{Khaykovich,Streker} have different motion patterns.

In this paper, we consider the 1D model in \cite{Khaykovich,Streker}: we
take $n=1$ in (\ref{Hamiltonian:Step B}). We derive rigorously the 1D cubic
focusing NLS from a 1D quantum many-body system. We establish the following
theorem.

\begin{theorem}[Main Theorem]
\label{THM:MainTHMFiniteKinetic}Assume that the pair interaction $V\ $is an
even Schwartz class function, which has a nonpositive integration, that is, $%
\int_{\mathbb{R}}V(x)dx\leqslant 0$, but may not be negative everywhere. Let 
$\psi _{N}\left( t,\mathbf{x}_{N}\right) $ be the $N-body$ Hamiltonian
evolution $e^{itH_{N}}\psi _{N}(0)$, where 
\begin{equation}
H_{N}=\sum_{j=1}^{N}\left( -\frac{1}{2}\partial _{x_{j}}^{2}+\frac{\omega
^{2}}{2}x_{j}^{2}\right) +\frac{1}{N}\sum_{1\leqslant i<j\leqslant
N}N^{\beta }V\left( N^{\beta }\left( x_{i}-x_{j}\right) \right)
\label{Hamiltonian:1D N-body}
\end{equation}%
for some $\omega \in \mathbb{R}$ which could be zero and for some $\beta \in
\left( 0,1\right) ,$ and let $\left\{ \gamma _{N}^{(k)}\right\} $ be the
family of marginal densities associated with $\psi _{N}$. Suppose that the
initial datum $\psi _{N}(0)$ verifies the following conditions:

(a) the initial datum is normalized, that is 
\begin{equation*}
\left\Vert \psi _{N}(0)\right\Vert _{L^{2}}=1,
\end{equation*}

(b) the initial datum is asymptotically factorized, in the sense that,%
\begin{equation}
\lim_{N\rightarrow \infty }\limfunc{Tr}\left\vert \gamma
_{N}^{(1)}(0,x_{1};x_{1}^{\prime })-\phi _{0}(x_{1})\overline{\phi _{0}}%
(x_{1}^{\prime })\right\vert =0,  \label{eqn:asym factorized}
\end{equation}%
for some one particle wave function $\phi _{0}$ s.t. $\left\Vert \left(
1-\partial _{x}^{2}+\omega ^{2}x^{2}\right) ^{\frac{1}{2}}\phi
_{0}\right\Vert _{L^{2}\left( \mathbb{R}\right) }<\infty $.

(c) the initial datum has finite kinetic energy and variance each particle%
\footnote{%
Finite variance can be dropped when $\omega $ is zero.} 
\begin{equation}
\sup_{j,N}\left\langle \psi _{N}(0),\left( -\partial _{x_{j}}^{2}+\omega
^{2}x_{j}^{2}\right) \psi _{N}(0)\right\rangle <\infty .
\label{Condition:FiniteKineticOnManyBodyInitialData}
\end{equation}%
Then $\forall t\geqslant 0$, $\forall k\geqslant 1$, we have the convergence
in the trace norm or the propagation of chaos that%
\begin{equation*}
\lim_{N\rightarrow \infty }\limfunc{Tr}\left\vert \gamma _{N}^{(k)}(t,%
\mathbf{x}_{k};\mathbf{x}_{k}^{\prime })-\dprod_{j=1}^{k}\phi (t,x_{j})%
\overline{\phi }(t,x_{j}^{\prime })\right\vert =0,
\end{equation*}%
where $\phi (t,x)$ is the solution to the 1D focusing cubic NLS%
\begin{eqnarray}
i\partial _{t}\phi &=&\left( -\frac{1}{2}\partial _{x_{j}}^{2}+\frac{\omega
^{2}}{2}x_{j}^{2}\right) \phi -b_{0}\left\vert \phi \right\vert ^{2}\phi 
\text{ in }\mathbb{R}^{1+1}  \label{equation:TargetCubicNLS} \\
\phi (0,x) &=&\phi _{0}(x)  \notag
\end{eqnarray}%
and the coupling constant $b_{0}=\left\vert \int_{\mathbb{R}%
}V(x)dx\right\vert .$
\end{theorem}

Theorem \ref{THM:MainTHMFiniteKinetic} is equivalent to the following
theorem.

\begin{theorem}[Main Theorem]
\label{THM:Main Theorem}Assume that the pair interaction $V\ $is an even
Schwartz class function, which has a nonpositive integration, that is, $%
\int_{\mathbb{R}}V(x)dx\leqslant 0$, but may not be negative everywhere. Let 
$\psi _{N}\left( t,\mathbf{x}_{N}\right) $ be the $N-body$ Hamiltonian
evolution $e^{itH_{N}}\psi _{N}(0)$ with $H_{N}$ given by (\ref%
{Hamiltonian:1D N-body}) for some $\omega \in \mathbb{R}$ which could be
zero and for some $\beta \in \left( 0,1\right) ,$ and let $\left\{ \gamma
_{N}^{(k)}\right\} $ be the family of marginal densities associated with $%
\psi _{N}$. Suppose that the initial datum $\psi _{N}(0)$ is normalized and
asymptotically factorized in the sense of (a) and (b) in Theorem \ref%
{THM:MainTHMFiniteKinetic} and verifies the following energy condition:

(c') there is a $C>0$ independent of $N$ or $k$ such that 
\begin{equation}
\left\langle \psi _{N}(0),H_{N}^{k}\psi _{N}(0)\right\rangle <C^{k}N^{k},%
\text{ }\forall k\geqslant 1.\footnote{%
Here, the energies $\left\langle \psi _{N}(0),H_{N}^{k}\psi
_{N}(0)\right\rangle $ are allowed to be negative. Estimate (\ref%
{estimate:why negative energy is fine}), in which we use (\ref%
{Condition:EnergyBoundOnManyBodyInitialData}), does not depend on the signs
of $\left\langle \psi _{N}(0),H_{N}^{k}\psi _{N}(0)\right\rangle .$ This is
not surprising because we are working in one dimension.}
\label{Condition:EnergyBoundOnManyBodyInitialData}
\end{equation}%
Then $\forall t\geqslant 0$, $\forall k\geqslant 1$, we have the convergence
in the trace norm or the propagation of chaos that%
\begin{equation*}
\lim_{N\rightarrow \infty }\limfunc{Tr}\left\vert \gamma _{N}^{(k)}(t,%
\mathbf{x}_{k};\mathbf{x}_{k}^{\prime })-\dprod_{j=1}^{k}\phi (t,x_{j})%
\overline{\phi }(t,x_{j}^{\prime })\right\vert =0,
\end{equation*}%
where $\phi (t,x)$ is the solution to the 1D focusing cubic NLS (\ref%
{equation:TargetCubicNLS}).
\end{theorem}

The equivalence of Theorems \ref{THM:MainTHMFiniteKinetic} and \ref{THM:Main
Theorem} for asymptotically factorized initial data has been used in all
defocusing works. In the main part of this paper, we prove Theorem \ref%
{THM:Main Theorem} in full detail. For completeness, we discuss briefly how
to deduce Theorem \ref{THM:MainTHMFiniteKinetic} from Theorem \ref{THM:Main
Theorem} in Appendix \ref{Appendix:equivalence}.

To our knowledge, Theorems \ref{THM:MainTHMFiniteKinetic} and \ref{THM:Main
Theorem} offer the first rigorous derivation of the focusing cubic NLS (\ref%
{equation:TargetCubicNLS}) from the $N$-body dynamic (\ref{Hamiltonian:Step
B}).\footnote{%
If one replaces the Hermite operator with $-\triangle $ or $(1-\triangle )^{%
\frac{1}{2}}$ in (\ref{Hamiltonian:Step B}) and let $\beta =0$, then there
are works by Erd\"{o}s and Yau \cite{E-Y1}, Michelangeli and Schlein \cite%
{MichelangeliSchlein}, and Ammari and Nier \cite{Ammari2,Ammari1} which
derive the corresponding focusing Hartree equations. For works on defocusing
Hartree dynamic ($\beta =0$)$,$ see \cite%
{Frolich,KnowlesAndPickl,RodnianskiAndSchlein,GMM1,GMM2,Chen2ndOrder, LChen}.%
} The main tool used in establishing Theorem \ref{THM:Main Theorem} is the
analysis of the focusing Bogoliubov--Born--Green--Kirkwood--Yvon hierarchy
(BBGKY) hierarchy of $\left\{ \gamma _{N}^{(k)}\right\} _{k=1}^{N}$ as $%
N\rightarrow \infty .$ With our definition, the sequence of the marginal
densities $\left\{ \gamma _{N}^{(k)}\right\} _{k=1}^{N}$ associated with $%
\psi _{N}$ solves the 1D BBGKY hierarchy with a quadratic trap%
\begin{eqnarray}
i\partial _{t}\gamma _{N}^{(k)} &=&\left[ -\frac{1}{2}\triangle _{\mathbf{x}%
_{k}}+\omega ^{2}\frac{\left\vert \mathbf{x}_{k}\right\vert ^{2}}{2},\gamma
_{N}^{(k)}\right] +\frac{1}{N}\sum_{1\leqslant i<j\leqslant k}\left[
N^{\beta }V\left( N^{\beta }\left( x_{i}-x_{j}\right) \right) ,\gamma
_{N}^{(k)}\right]   \label{hierarchy:TheBBGKYHierarchyWithTraps} \\
&&+\frac{N-k}{N}\sum_{j=1}^{k}\limfunc{Tr}\nolimits_{k+1}\left[ N^{\beta
}V\left( N^{\beta }\left( x_{j}-x_{k+1}\right) \right) ,\gamma _{N}^{(k+1)}%
\right] .  \notag
\end{eqnarray}

In the classical setting, deriving mean-field type equations by studying the
limit of the BBGKY hierarchy was proposed by Kac and demonstrated by
Landford's work on the Boltzmann equation. In the quantum setting, the usage
of the BBGKY hierarchy was suggested by Spohn \cite{Spohn} and has been
proven to be successful by Elgart, Erd\"{o}s, Schlein, and Yau in their
fundamental papers \cite{E-E-S-Y1, E-S-Y1,E-S-Y2,E-S-Y5, E-S-Y3}\footnote{%
Around the same time, there was the 1D defocusing work \cite{AGT}.} which
rigorously derives the 3D cubic defocusing NLS from a 3D quantum many-body
dynamic without trapping. The Elgart-Erd\"{o}s-Schlein-Yau program\footnote{%
See \cite{SchleinNew,GM1,Pickl} for different approaches.} consists of two
principal parts: in one part, they consider the sequence of the marginal
densities $\left\{ \gamma _{N}^{(k)}\right\} $ and prove that its
appropriate limit as $N\rightarrow \infty $ solves the 3D defocusing
Gross-Pitaevskii (GP) hierarchy 
\begin{equation}
i\partial _{t}\gamma ^{(k)}=\left[ -\frac{1}{2}\triangle _{\mathbf{x}%
_{k}},\gamma ^{(k)}\right] +b_{0}\sum_{j=1}^{k}\limfunc{Tr}%
\nolimits_{_{k+1}}[\delta (x_{j}-x_{k+1}),\gamma ^{(k+1)}],\text{ }%
b_{0}\geqslant 0.  \label{equation:Gross-Pitaevskii hiearchy without a trap}
\end{equation}%
In another part, they show that hierarchy 
\eqref{equation:Gross-Pitaevskii
hiearchy without a trap} has a unique solution which is therefore a
completely factorized state. However, the uniqueness theory for hierarchy 
\eqref{equation:Gross-Pitaevskii hiearchy
without a trap} is surprisingly delicate due to the fact that it is a system
of infinitely many coupled equations over an unbounded number of variables.
In \cite{KlainermanAndMachedon}, by assuming a space-time bound on the limit
of $\left\{ \gamma _{N}^{(k)}\right\} $, Klainerman and Machedon gave
another uniqueness theorem regarding (\ref{equation:Gross-Pitaevskii
hiearchy without a trap}) through a collapsing estimate originating from the
multilinear Strichartz estimates and a board game argument inspired by the
Feynman graph argument in \cite{E-S-Y2}.

Later, the method in Klainerman and Machedon \cite{KlainermanAndMachedon}
was taken up by Kirkpatrick, Schlein, and Staffilani \cite{Kirpatrick}, who
derived the 2D cubic defocusing NLS from the 2D quantum many-body dynamic;
by Chen and Pavlovi\'{c} \cite{TChenAndNpGP1, TChenAndNP}, who considered
the 1D and 2D 3-body repelling interaction problem and the general existence
theory of hierarchy $%
\eqref{equation:Gross-Pitaevskii hiearchy without a
trap}$; by X.C. \cite{ChenAnisotropic, Chen3DDerivation}, who investigated
the defocusing problem with trapping in 2D and 3D; and by X.C. and J.H. \cite%
{C-H3Dto2D}, who proved the effectiveness of the 3D to 2D reduction problem.
In \cite{TCNPNT, TCNPNT1}, Chen, Pavlovi\'{c} and Tzirakis worked out the
virial and Morawetz identities for hierarchy 
\eqref{equation:Gross-Pitaevskii
hiearchy without a trap} and showed the blow up for hierarchy 
\eqref{equation:Gross-Pitaevskii
hiearchy without a trap} in 2D and 3D in the case of negative energy initial
data and negative $b_{0}$. In \cite{Sohinger}, Gressman, Sohinger, and
Staffilani have obtained a uniqueness theorem of solution to hierarchy 
\eqref{equation:Gross-Pitaevskii
hiearchy without a trap} in 3D subject to periodic boundary condition.

Recently, in \cite{TChenAndNPSpace-Time}, for the 3D defocusing problem
without traps, Chen and Pavlovi\'{c} showed that, for $\beta \in (0,1/4)$,
the limit of the BBGKY sequence satisfies the space-time bound assumed by
Klainerman and Machedon \cite{KlainermanAndMachedon} as $N\rightarrow \infty 
$. In \cite{Chen3DDerivation}, X.C. extended and simplified their method to
study the 3D trapping problem for $\beta \in (0,2/7].$ X.C. and J.H. \cite%
{C-H2/3}\ then extended the $\beta \in (0,2/7]$ result by X.C. to $\beta \in
(0,2/3)$ using $X_{b}$ spaces and Littlewood-Paley theory.

We use the Klainerman-Machedon framework for the uniqueness argument in this
paper. While the known uniqueness theorems \cite{AGT, TChenAndNpGP1,
TChenAndNP} regarding the 1D GP hierarchy need $H^{\frac{1}{2}+\varepsilon }$
smoothness, i.e. more than the continuity in 1D, our Theorem \ref%
{THM:OptimalUniqueness of 1D GP} requires merely $H^{\varepsilon }$
regularity to establish uniqueness$.$ To achieve this reduction, we prove
the optimal 1D collapsing estimate which has been open for a while.

\begin{theorem}
\footnote{%
For more estimates of this type, see \cite%
{KlainermanAndMachedon,Kirpatrick,GM,ChenDie,ChenAnisotropic,Beckner,Sohinger}%
}\label{THM:Optimal1DCollpasing}Let $U^{(k)}\left( \tau \right)
=\dprod\limits_{j=1}^{k}e^{i\tau \partial _{y_{j}}^{2}}e^{-i\tau \partial
_{y_{j}^{\prime }}^{2}}$ and $R_{\varepsilon
}^{(k)}=\prod_{j=1}^{k}\left\langle \partial _{y_{j}}\right\rangle
^{\varepsilon }\left\langle \partial _{y_{j}^{\prime }}\right\rangle
^{\varepsilon }$. Define the collision operator $B_{j,k+1}$ by 
\begin{equation*}
B_{j,k+1}u^{(k+1)}=\limfunc{Tr}\nolimits_{k+1}\left[ \delta \left(
y_{j}-y_{k+1}\right) ,u^{(k+1)}\right] .
\end{equation*}%
Given any finite time $T$ and any $\varepsilon >0$, there is a constant $%
C_{T}>0$ independent of $j,k$ and $\phi ^{(k+1)}$, such that%
\begin{equation*}
\left\Vert R_{\varepsilon }^{(k)}B_{j,k+1}U^{(k+1)}(\tau )\phi
^{(k+1)}\right\Vert _{L_{T}^{2}L_{\mathbf{y,y}^{\prime }}^{2}}\leqslant
C_{T}\left\Vert R_{\varepsilon }^{(k+1)}\phi ^{(k+1)}\right\Vert _{L_{%
\mathbf{y,y}^{\prime }}^{2}}.
\end{equation*}%
This estimate is optimal in the sense that it fails whenever $T=\infty $ or $%
\varepsilon =0.$
\end{theorem}

It is surprising that the 1D scale-invariant global-in-time collapsing
estimate $\left( T=\infty \text{ and }\varepsilon =0\right) $ fails while
the scale-invariant global-in-time estimates are true in 2D \cite%
{ChenAnisotropic, Beckner}\ and 3D \cite{KlainermanAndMachedon}. The failure
of Theorem \ref{THM:Optimal1DCollpasing} when $T=\infty $ for any $%
\epsilon\geq 0$ indicates that we do not have enough decay in time in 1D.
Since collapsing estimates like Theorem \ref{THM:Optimal1DCollpasing}
determine many features of the corresponding GP hierarchies, we wonder if
this is related to the fact that there is no $L^{2}$ small data scattering
theory for the ordinary 1D focusing GP hierarchy 
\begin{equation*}
i\partial _{t}\gamma ^{(k)}=\left[ -\frac{1}{2}\triangle _{\mathbf{x}%
_{k}},\gamma ^{(k)}\right] -b_{0}\sum_{j=1}^{k}\left[ \delta
(x_{j}-x_{k+1}),\gamma ^{(k+1)}\right] .
\end{equation*}%
Specifically, we can write down a tensor product of 1D NLS solitons
arbitrarily small in any unweighted $H^{s}$ norm, $s\geq 0$, which shows the
lack of small data scattering for GP. If we conjecture that in a general
setting, scale-invariant global-in-time collapsing estimates from \cite%
{KlainermanAndMachedon,ChenAnisotropic,Beckner}\ could be part of a proof of
small data scattering, then the above mentioned lack of scattering in 1D
implies the nonexistence of global-in-time collapsing estimates in 1D. This
heuristically implies the optimality of Theorem \ref{THM:Optimal1DCollpasing}%
. On the other hand, the known global-in-time collapsing estimates in 2D and
3D could eventually be used to prove small data scattering for 2D and 3D GP.
All of these remarks pertain to unweighted Sobolev spaces at or above the
critical (scale invariant) level; in the setting of weighted Sobolev spaces,
small-data scattering for 1D cubic NLS is known \cite{Hayashi-Naumkin}.

Theorem \ref{THM:Optimal1DCollpasing} also reduces the regularity
requirement by $1/2$ for the current local existence theory \cite%
{TChenAndNpGP1}\ of the GP hierarchy (\ref{equation:Gross-Pitaevskii
hiearchy without a trap}) subject to general initial data in 1D. In fact,
plugging Theorem \ref{THM:Optimal1DCollpasing} into \cite{TChenAndNpGP1}
yields the following corollary.

\begin{corollary}
For every initial data in $H^{\varepsilon }$ which is not necessarily
factorized, there is a time $T>0$ such that there exists a unique solution
in $H^{\varepsilon }$ for $t\in \lbrack 0,T]$ to the GP hierarchy (\ref%
{equation:Gross-Pitaevskii hiearchy without a trap}) in 1D regardless of the
sign of $b_{0}.$
\end{corollary}

\subsection{Organization of the Paper}

We first review the lens transform and its relevant properties in \S \ref%
{Section:Lens Transform}. It aids in the proof of the main theorem in the
sense that it links the analysis of $-\partial _{x}^{2}+\omega ^{2}x^{2}$ to
the analysis of $-\partial _{y}^{2}$ which is easier to deal with using the
Fourier transform. With the lens transform, we then outline the proof of our
main theorem, Theorem \ref{THM:Main Theorem}, in \S \ref{Section:Outline}.
The components of the proof are in \S \ref{Section:EnergyEstimate}, \ref%
{Section:Compactness and convergence}, and \ref{Section:Uniqueness}.

In \S \ref{Section:EnergyEstimate}, we prove the needed energy estimate for
the focusing $N$-body Schr\"odinger evolution. The key obstacle here,
compared to earlier versions of such estimates in the defocusing works \cite%
{AGT, E-E-S-Y1, E-S-Y1,E-S-Y2,E-S-Y5,
E-S-Y3,TChenAndNP,ChenAnisotropic,TChenAndNPSpace-Time, Chen3DDerivation,
C-H3Dto2D}, is to accommodate the negativity of the potential. We first
observe a new decomposition of the Hamiltonian $H_N$ given by 
\begin{equation*}
N^{-1} H_N + \|V\|_{L^1}^2+ 1= \frac{1}{2N(N-1)} \sum_{1\leq i,j \leq N}
H_{+ij}
\end{equation*}
where 
\begin{equation*}
H_{+ij} = S_i^2 +S_j^2 + \frac{N-1}{N} N^\beta V(N^\beta(x_i-x_j)) +
2\|V\|_{L^1}^2
\end{equation*}
and 
\begin{equation*}
S_j = ( 1- \tfrac12 \partial_{x_j}^2 +\tfrac12 \omega^2 x_j^2)^{1/2}
\end{equation*}
In the expansion of $(N^{-1}H_N + \|V\|_{L^1}^2+ 1)^k$, the terms that occur
most frequently are of the form 
\begin{equation*}
H_{+i_1j_1} \cdots H_{+i_kj_k}
\end{equation*}
with all $i_1, j_1, \ldots, i_k,j_k$ distinct. Since these operators
pairwise commute, we can exploit the positivity of each $H_{+ij}$. In
particular, we have 
\begin{equation*}
H_{+ij} \geq \tfrac12(S_i^2+S_j^2)
\end{equation*}
We justify the above heuristic by induction.

In \S \ref{Section:Compactness and convergence}, we use the energy estimates
derived in \S \ref{Section:EnergyEstimate} and duality to prove weak*
compactness and convergence of the corresponding BBGKY hierarchy. This
follows the similar procedure in the defocusing works.

Finally, in \S \ref{Section:Uniqueness}, we prove Theorem \ref%
{THM:Optimal1DCollpasing}, the optimal 1D collapsing estimate. As discussed
previously, we need to include a time-localization. On the Fourier side, the
time localization mollifies the resulting surface measure and makes it
integrable. Without the time localization, the surface measure remains
unmollified and is not integrable, and the estimate fails. The optimality
statement essentially follows.

\subsection{Acknowledgements}

J.H. was supported in part by NSF grant DMS-0901582 and a Sloan Research
Fellowship (BR-4919), and X.C. received travel support in part from the same
Sloan Fellowship to visit U. Maryland during May 12th to 17th. We would like
to thank A. Soffer for encouraging us to work on this problem during X.C's
visit to U. Maryland.

\section{Lens Transform\label{Section:Lens Transform}}

In this section, we review the lens transform and its relevant properties.
Everything here comes from \cite{Chen3DDerivation}. (See also \cite%
{Carles,ChenAnisotropic}.) We include it solely for completeness. The lens
transform aids in the proof of the main theorem in the sense that it links
the analysis of $-\partial _{x}^{2}+\omega ^{2}x^{2}$ to the analysis of $%
-\partial _{y}^{2}$ which is a better understood operator. We remark that
the lens transform is exactly the identity when $\omega =0$ i.e. this
section is trivial when $\omega =0$.

We denote $(t,x)$ the space-time on the Hermite side and $(\tau ,y)$ the
space-time on the Laplacian side. We define the lens transform in
Definitions \ref{def:functionG-Lens} and \ref{def:KernelG-Lens}. We then
explain how the lens transform acts on the BBGKY hierarchy and the GP
hierarchy via Lemmas \ref{Lemma:BBGKY hierarchy under GLens} and \ref%
{Lemma:GP hierarchy under GLens}. Finally, we relate the trace norms and the
energies of the two sides of the lens transform through Lemmas \ref%
{Lemma:TraceNormPreservation} and \ref{Lemma:EnergyBetweenLensTransform}.

\begin{definition}[\protect\cite{Chen3DDerivation}]
\label{def:functionG-Lens}Let $\mathbf{x}_{N},\mathbf{y}_{N}\in \mathbb{R}%
^{N}$. We define the lens transform for $L^{2}$ functions $M_{N}:L^{2}(d%
\mathbf{y}_{N})\rightarrow $ $L^{2}(d\mathbf{x}_{N})$ and its inverse by%
\begin{eqnarray*}
\left( M_{N}u_{N}\right) (t,\mathbf{x}_{N}) &=&\frac{e^{-i\omega \tan \omega
t\frac{\left\vert \mathbf{x}_{N}\right\vert ^{2}}{2}}}{(\cos \omega t)^{%
\frac{N}{2}}}u_{N}(\frac{\tan \omega t}{\omega },\frac{\mathbf{x}_{N}}{\cos
\omega t}) \\
\left( M_{N}^{-1}\psi _{N}\right) (\tau ,\mathbf{y}_{N}) &=&\frac{e^{i\frac{%
\omega ^{2}\tau }{1+\omega ^{2}\tau ^{2}}\frac{\left\vert \mathbf{y}%
_{N}\right\vert ^{2}}{2}}}{\left( 1+\omega ^{2}\tau ^{2}\right) ^{\frac{N}{4}%
}}\psi _{N}(\frac{\arctan \left( \omega \tau \right) }{\omega },\frac{%
\mathbf{y}_{N}}{\sqrt{1+\omega ^{2}\tau ^{2}}}).
\end{eqnarray*}%
$M_{N}$ is unitary by definition and the variables are related by%
\begin{equation*}
\tau =\frac{\tan \omega t}{\omega },\text{ }\mathbf{y}_{N}=\frac{\mathbf{x}%
_{N}}{\cos \omega t}.
\end{equation*}
\end{definition}

\begin{definition}[\protect\cite{Chen3DDerivation}]
\label{def:KernelG-Lens}Let $\mathbf{x}_{k},\mathbf{x}_{k}^{\prime },\mathbf{%
y}_{k},\mathbf{y}_{k}^{\prime }\in \mathbb{R}^{k}$. We define the lens
transform for Hilbert-Schmidt kernels $T_{k}:$ $L^{2}(d\mathbf{y}_{k}d%
\mathbf{y}_{k}^{\prime })\rightarrow $ $L^{2}(d\mathbf{x}_{k}d\mathbf{x}%
_{k}^{\prime })$ and its inverse by%
\begin{eqnarray*}
&&\left( T_{k}u^{(k)}\right) (t,\mathbf{x}_{k};\mathbf{x}_{k}^{\prime }) \\
&=&\frac{e^{-i\omega \tan \omega t\frac{\left( \left\vert \mathbf{x}%
_{k}\right\vert ^{2}-\left\vert \mathbf{x}_{k}^{\prime }\right\vert
^{2}\right) }{2}}}{(\cos \omega t)^{k}}u^{(k)}(\frac{\tan \omega t}{\omega },%
\frac{\mathbf{x}_{k}}{\cos \omega t};\frac{\mathbf{x}_{k}^{\prime }}{\cos
\omega t})
\end{eqnarray*}%
\begin{eqnarray*}
&&\left( T_{k}^{-1}\gamma ^{(k)}\right) (\tau ,\mathbf{y}_{k};\mathbf{y}%
_{k}^{\prime }) \\
&=&\frac{e^{i\frac{\omega ^{2}\tau }{1+\omega ^{2}\tau ^{2}}\frac{\left(
\left\vert \mathbf{y}_{k}\right\vert ^{2}-\left\vert \mathbf{y}_{k}^{\prime
}\right\vert ^{2}\right) }{2}}}{\left( 1+\omega ^{2}\tau ^{2}\right) ^{\frac{%
k}{2}}}\gamma ^{(k)}(\frac{\arctan \left( \omega \tau \right) }{\omega },%
\frac{\mathbf{y}_{k}}{\sqrt{1+\omega ^{2}\tau ^{2}}};\frac{\mathbf{y}%
_{k}^{\prime }}{\sqrt{1+\omega ^{2}\tau ^{2}}}).
\end{eqnarray*}%
$T_{k}$ is unitary by definition as well and the variables are again related
by%
\begin{equation*}
\tau =\frac{\tan \omega t}{\omega },\text{ }\mathbf{y}_{k}=\frac{\mathbf{x}%
_{k}}{\cos \omega t}\text{ and }\mathbf{y}_{k}^{\prime }=\frac{\mathbf{x}%
_{k}^{\prime }}{\cos \omega t}.
\end{equation*}%
In particular, if $u_{N}(\tau ,\mathbf{y}_{N})=M_{N}^{-1}\left( \psi
_{N}\right) ,$ then $\left\{ u_{N}^{(k)}=T_{k}^{-1}\gamma _{N}^{(k)}\right\} 
$ is exactly the family of marginal densities associated with $u_{N}.$
\end{definition}

\begin{lemma}[\protect\cite{Chen3DDerivation}]
\label{Lemma:BBGKY hierarchy under GLens}Write $V_{N}\left( x\right)
=N^{\beta }V\left( N^{\beta }x\right) .$ $\left\{ \gamma _{N}^{(k)}\right\} $
solves the 1D BBGKY hierarchy with a quadratic trap (\ref%
{hierarchy:TheBBGKYHierarchyWithTraps}) in $\left[ -T_{0},T_{0}\right] $ if
and only if $\left\{ u_{N}^{(k)}=T_{k}^{-1}\gamma _{N}^{(k)}\right\} $
solves the hierarchy%
\begin{eqnarray}
i\partial _{\tau }u_{N}^{(k)} &=&\left[ -\frac{1}{2}\triangle _{\mathbf{y}%
_{k}},\gamma _{N}^{(k)}\right] +\frac{1}{\left( 1+\omega ^{2}\tau
^{2}\right) }\frac{1}{N}\sum_{1\leqslant i<j\leqslant k}\left[ V_{N}(\frac{%
y_{i}-y_{j}}{\left( 1+\omega ^{2}\tau ^{2}\right) ^{\frac{1}{2}}}%
),u_{N}^{(k)}\right]  \label{hierarchy:LensBBGKY} \\
&&+\frac{N-k}{N}\frac{1}{\left( 1+\omega ^{2}\tau ^{2}\right) }\sum_{j=1}^{k}%
\limfunc{Tr}\nolimits_{k+1}\left[ V_{N}\left( \frac{y_{j}-y_{k+1}}{\left(
1+\omega ^{2}\tau ^{2}\right) ^{\frac{1}{2}}}\right) ,u_{N}^{(k+1)}\right] 
\notag
\end{eqnarray}%
in $\left[ -\frac{\tan \omega T_{0}}{\omega },\frac{\tan \omega T_{0}}{%
\omega }\right] .$
\end{lemma}

\begin{lemma}[\protect\cite{Chen3DDerivation}]
\label{Lemma:GP hierarchy under GLens}$\left\{ \gamma ^{(k)}\right\} $
solves the 1D focusing GP hierarchy with a quadratic trap 
\begin{equation}
i\partial _{t}\gamma ^{(k)}=\left[ -\frac{1}{2}\triangle _{\mathbf{x}%
_{k}}+\omega ^{2}\frac{\left\vert \mathbf{x}_{k}\right\vert ^{2}}{2},\gamma
^{(k)}\right] -b_{0}\sum_{j=1}^{k}\left[ \delta (x_{j}-x_{k+1}),\gamma
^{(k+1)}\right] ,  \label{hierarchy:TheGPHierarchyWithTraps}
\end{equation}%
in $\left[ -T_{0},T_{0}\right] $ if and only if $\left\{
u^{(k)}=T_{k}^{-1}\gamma ^{(k)}\right\} $ solves the focusing hierarchy%
\begin{equation}
i\partial _{\tau }u^{(k)}=\left[ -\frac{1}{2}\triangle _{\mathbf{y}%
_{k}},u^{(k)}\right] -\frac{1}{\left( 1+\omega ^{2}\tau ^{2}\right) ^{\frac{1%
}{2}}}b_{0}\sum_{j=1}^{k}\left[ \delta (y_{j}-y_{k+1}),u^{(k+1)}\right] ,
\label{hierarchy:LensGP}
\end{equation}%
in $\left[ -\frac{\tan \omega T_{0}}{\omega },\frac{\tan \omega T_{0}}{%
\omega }\right] .$
\end{lemma}

\begin{lemma}[\protect\cite{Chen3DDerivation}]
\label{Lemma:TraceNormPreservation}If $K(\mathbf{y}_{k},\mathbf{y}%
_{k}^{\prime })$ is the kernel of a self-adjoint trace class operator on $%
L^{2}\left( \mathbb{R}^{k}\right) $, then the eigenvectors of the kernel $%
\left( T_{k}K\right) (\mathbf{x}_{k},\mathbf{x}_{k}^{\prime })$ are exactly
the lens transform of the eigenvectors of the kernel $K(\mathbf{y}_{k},%
\mathbf{y}_{k}^{\prime })$ with the same eigenvalues. In particular, we have%
\begin{equation*}
\limfunc{Tr}\left\vert T_{k}K\right\vert =\limfunc{Tr}\left\vert
K\right\vert .
\end{equation*}
\end{lemma}

\begin{lemma}[\protect\cite{Chen3DDerivation}]
\label{Lemma:EnergyBetweenLensTransform}There is a $C>0$ such that 
\begin{equation*}
\left\langle u_{N}(\tau ),\dprod\limits_{j=1}^{k}\left( 1-\partial
_{y_{j}}^{2}\right) u_{N}(\tau )\right\rangle \leqslant C^{k}\left\langle
\psi _{N}(t),\dprod\limits_{j=1}^{k}\left( 1-\frac{1}{2}\partial
_{x_{j}}^{2}+\frac{1}{2}\omega ^{2}x_{j}^{2}\right) \psi _{N}(t)\right\rangle
\end{equation*}%
for all $\psi _{N}(t,\mathbf{x}_{N})$, where $u_{N}(\tau ,\mathbf{y}%
_{N})=M_{N}^{-1}\left( \psi _{N}\right) $. In particular, if $%
u_{N}^{(k)}=T_{k}^{-1}\gamma _{N}^{(k)}$, we have%
\begin{eqnarray*}
&&\limfunc{Tr}\left( \dprod\limits_{j=1}^{k}\left( 1-\partial
_{y_{j}}^{2}\right) ^{\frac{1}{2}}\right) u_{N}^{(k)}(\tau )\left(
\dprod\limits_{j=1}^{k}\left( 1-\partial _{y_{j}}^{2}\right) ^{\frac{1}{2}%
}\right) \\
&\leqslant &C^{k}\limfunc{Tr}\left( \dprod\limits_{j=1}^{k}\left( 1-\frac{1}{%
2}\partial _{x_{j}}^{2}+\frac{1}{2}\omega ^{2}x_{j}^{2}\right) ^{\frac{1}{2}%
}\right) \gamma _{N}^{(k)}(t)\left( \dprod\limits_{j=1}^{k}\left( 1-\frac{1}{%
2}\partial _{x_{j}}^{2}+\frac{1}{2}\omega ^{2}x_{j}^{2}\right) ^{\frac{1}{2}%
}\right) .
\end{eqnarray*}
\end{lemma}

\begin{notation}
From here on out, to make formulas shorter, we write%
\begin{equation*}
L_{j}=\left( 1-\partial _{y_{j}}^{2}\right) ^{\frac{1}{2}},\text{ }%
S_{j}=\left( 1-\frac{1}{2}\partial _{x_{j}}^{2}+\frac{1}{2}\omega
^{2}x_{j}^{2}\right) ^{\frac{1}{2}},
\end{equation*}%
\begin{equation*}
L^{(k)}=\dprod\limits_{j=1}^{k}L_{j},\text{ }S^{(k)}=\dprod%
\limits_{j=1}^{k}S_{j},
\end{equation*}%
\begin{equation*}
g(\tau )=\left( 1+\omega ^{2}\tau ^{2}\right) ^{-\frac{1}{2}},V_{N,\tau
}(y)=N^{\beta }g(\tau )V(N^{\beta }g(\tau )y).
\end{equation*}%
The only properties we need are $0<g(\tau )\leqslant 1,$ and $\int $ $%
V_{N,\tau }(y)dy=b_{0}.$
\end{notation}

\section{Proof of the Main Theorem (Theorem \protect\ref{THM:Main Theorem}) 
\label{Section:Outline}}

We start by introducing an appropriate topology on the density matrices as
was previously done in \cite{E-E-S-Y1, E-Y1, E-S-Y1,E-S-Y2,E-S-Y5,
E-S-Y3,Kirpatrick,TChenAndNP,ChenAnisotropic,Chen3DDerivation,C-H3Dto2D,C-H2/3}%
. Denote the spaces of compact operators and trace class operators on $%
L^{2}\left( \mathbb{R}^{k}\right) $ as $\mathcal{K}_{k}$ and $\mathcal{L}%
_{k}^{1}$, respectively. Then $\left( \mathcal{K}_{k}\right) ^{\prime }=%
\mathcal{L}_{k}^{1}$. By the fact that $\mathcal{K}_{k}$ is separable, we
select a dense countable subset $\{J_{i}^{(k)}\}_{i\geqslant 1}\subset 
\mathcal{K}_{k}$ in the unit ball of $\mathcal{K}_{k}$ (so $\Vert
J_{i}^{(k)}\Vert _{\func{op}}\leqslant 1$ where $\left\Vert \cdot
\right\Vert _{\func{op}}$ is the operator norm). For $\gamma ^{(k)},\tilde{%
\gamma}^{(k)}\in \mathcal{L}_{k}^{1}$, we then define a metric $d_{k}$ on $%
\mathcal{L}_{k}^{1}$ by 
\begin{equation*}
d_{k}(\gamma ^{(k)},\tilde{\gamma}^{(k)})=\sum_{i=1}^{\infty
}2^{-i}\left\vert \limfunc{Tr}J_{i}^{(k)}\left( \gamma ^{(k)}-\tilde{\gamma}%
^{(k)}\right) \right\vert .
\end{equation*}%
A uniformly bounded sequence $\gamma _{N}^{(k)}\in \mathcal{L}_{k}^{1}$
converges to $\gamma ^{(k)}\in \mathcal{L}_{k}^{1}$ with respect to the
weak* topology if and only if 
\begin{equation*}
\lim_{N\rightarrow \infty }d_{k}(\gamma _{N}^{(k)},\gamma ^{(k)})=0.
\end{equation*}%
For fixed $T>0$, let $C\left( \left[ 0,T\right] ,\mathcal{L}_{k}^{1}\right) $
be the space of functions of $t\in \left[ 0,T\right] $ with values in $%
\mathcal{L}_{k}^{1}$ which are continuous with respect to the metric $d_{k}.$
On $C\left( \left[ 0,T\right] ,\mathcal{L}_{k}^{1}\right) ,$ we define the
metric 
\begin{equation*}
\hat{d}_{k}(\gamma ^{(k)}\left( \cdot \right) ,\tilde{\gamma}^{(k)}\left(
\cdot \right) )=\sup_{t\in \left[ 0,T\right] }d_{k}(\gamma ^{(k)}\left(
t\right) ,\tilde{\gamma}^{(k)}\left( t\right) ),
\end{equation*}%
and denote by $\tau _{prod}$ the topology on the space $\oplus _{k\geqslant
1}C\left( \left[ 0,T\right] ,\mathcal{L}_{k}^{1}\right) $ given by the
product of topologies generated by the metrics $\hat{d}_{k}$ on $C\left( %
\left[ 0,T\right] ,\mathcal{L}_{k}^{1}\right) .$

With the above topology on the space of marginal densities, we can begin the
proof of Theorem \ref{THM:Main Theorem}. We divide the proof into four steps.

\begin{itemize}
\item[Step I] (Energy estimate) Before we apply the lens transform to our
problem, we first establish, through an elaborate calculation in Theorem \ref%
{THM:EnergyEstimate}, that one can absorb the negativity of the interaction
in (\ref{Hamiltonian:1D N-body}). Henceforth we transform the energy
condition (\ref{Condition:EnergyBoundOnManyBodyInitialData}) into a $H^{1}$
type bound. Due to the fact that the quantity $\left\langle \psi
_{N}(0),H_{N}^{k}\psi _{N}(0)\right\rangle $ in (\ref%
{Condition:EnergyBoundOnManyBodyInitialData}) is conserved by the evolution,
we deduce the \emph{a priori} bound on the marginal densities%
\begin{equation*}
\sup_{t}\limfunc{Tr}S^{(k)}\gamma _{N}^{(k)}(t)S^{(k)}\leqslant C^{k}.
\end{equation*}%
In Corollary \ref{Corollary:EnergyBoundForMarginals}, we then combine the
above bound and Lemma \ref{Lemma:EnergyBetweenLensTransform} to obtain the $%
H^{1}$ bound%
\begin{equation}
\sup_{\tau \in \left[ -\frac{\tan \omega T_{0}}{\omega },\frac{\tan \omega
T_{0}}{\omega }\right] }\limfunc{Tr}L^{(k)}u_{N}^{(k)}(\tau
)L^{(k)}\leqslant C^{k}\text{, if }T_{0}<\frac{\pi }{2\omega },
\label{EnergyBound:LensBBGKY}
\end{equation}%
where $u_{N}^{(k)}=T_{k}^{-1}\gamma _{N}^{(k)}$.

\item[Step II] (Compactness and Convergence) Fix $T_{0}<\frac{\pi }{2\omega }
$ and employ (\ref{EnergyBound:LensBBGKY}), we prove, in Theorem \ref%
{Theorem:CompactnessOfBBGKY}, that the sequence $\Gamma _{N}(\tau )=\left\{
u_{N}^{(k)}\right\} _{k=1}^{N}$ which satisfies the 1D BBGKY hierarchy (\ref%
{hierarchy:LensBBGKY}) is compact with respect to the product topology $\tau
_{prod}$. Moreover, we prove, in Theorem \ref{THM:convergence to GP}, that
if $\Gamma (\tau )=\left\{ u^{(k)}\right\} _{k=1}^{\infty }$ is a limit
point of $\Gamma _{N}(\tau )$ with respect to the product topology $\tau
_{prod}$, then $\Gamma (\tau )$ is a solution to the focusing GP hierarchy (%
\ref{hierarchy:LensGP}) subject to initial data $u^{(k)}\left( 0\right)
=\left\vert \phi _{0}\right\rangle \left\langle \phi _{0}\right\vert
^{\otimes k}$ and the coupling constant is given by $b_{0}=$ $\left\vert
\int V\left( x\right) dx\right\vert .$ This is a well-known argument used in 
\cite{E-E-S-Y1, E-Y1, E-S-Y1,E-S-Y2,E-S-Y5,
E-S-Y3,Kirpatrick,TChenAndNP,ChenAnisotropic,Chen3DDerivation,C-H3Dto2D}, we
include the proof in \S \ref{Section:Compactness and convergence} for
completeness since it is the first time such an argument is used in the
focusing setting.

\item[Step III] (Uniqueness) When $u^{(k)}\left( 0\right) =\left\vert \phi
_{0}\right\rangle \left\langle \phi _{0}\right\vert ^{\otimes k}$, we know
that there is a special solution to the focusing GP hierarchy (\ref%
{hierarchy:LensGP}), namely%
\begin{equation}
u^{(k)}(\tau ,\mathbf{y}_{k},\mathbf{y}_{k}^{\prime
})=\dprod\limits_{j=1}^{k}\tilde{\phi}(\tau ,y_{j})\overline{\tilde{\phi}%
(\tau ,y_{j}^{\prime })}  \label{limit:GlensTarget}
\end{equation}%
where $\tilde{\phi}$ solves%
\begin{eqnarray}
i\partial _{\tau }\tilde{\phi} &=&-\partial _{y}^{2}\tilde{\phi}-g(\tau
)b_{0}\left\vert \tilde{\phi}\right\vert ^{2}\tilde{\phi}
\label{equation:LensTargetCubicNLS} \\
\tilde{\phi}(0,y) &=&\phi _{0}.  \notag
\end{eqnarray}%
A suitable uniqueness theorem regarding (\ref{hierarchy:LensGP}) will then
identify all limit points of $\Gamma _{N}(\tau )$ obtained in Step II with (%
\ref{limit:GlensTarget}) for us. The Klainerman-Machedon scheme, introduced
in \cite{KlainermanAndMachedon} and used in \cite%
{Kirpatrick,TChenAndNP,ChenAnisotropic,Chen3DDerivation,C-H3Dto2D,C-H2/3},
transforms Theorem \ref{THM:Optimal1DCollpasing} into the following
uniqueness theorem.
\end{itemize}

\begin{theorem}
\label{THM:OptimalUniqueness of 1D GP}Let $R_{\varepsilon }^{(k)}$ and $%
B_{j,k+1}$ be defined in Theorem \ref{THM:Optimal1DCollpasing}. Suppose that 
$\left\{ u^{(k)}\right\} _{k=1}^{\infty }$ solves the 1D focusing GP
hierarchy (\ref{hierarchy:LensGP}) subject to zero initial data and the
space-time bound 
\begin{equation}
\int_{0}^{T}\left\Vert R_{\varepsilon }^{(k)}B_{j,k+1}u^{(k+1)}(\tau ,%
\mathbf{\cdot };\mathbf{\cdot })\right\Vert _{L_{\mathbf{y,y}^{\prime
}}^{2}}d\tau \leqslant C^{k}  \label{Condition:Space-Time Bound}
\end{equation}%
for some $\varepsilon ,C>0$ and all $1\leqslant j\leqslant k.$ Then $\forall
k,\tau \in \lbrack 0,T]$, 
\begin{equation*}
\left\Vert R_{\varepsilon }^{(k)}u^{(k)}(\tau ,\mathbf{\cdot };\mathbf{\cdot 
})\right\Vert _{L_{\mathbf{y,y}^{\prime }}^{2}}=0.
\end{equation*}
\end{theorem}

\begin{proof}
Once we prove Theorem \ref{THM:Optimal1DCollpasing}, Theorem \ref%
{THM:OptimalUniqueness of 1D GP} follows from the proof of \cite[Theorem 6]%
{Chen3DDerivation} line by line.
\end{proof}

To apply Theorem \ref{THM:OptimalUniqueness of 1D GP}, we need to check (\ref%
{Condition:Space-Time Bound}). As the spatial dimension is one, the
following trace lemma and (\ref{EnergyBound:LensBBGKY}) takes care of (\ref%
{Condition:Space-Time Bound}) for us.\footnote{%
Verifying (\ref{Condition:Space-Time Bound}) in 3D is highly nontrivial and
is merely partially solved so far. See \cite%
{TChenAndNPSpace-Time,Chen3DDerivation,C-H2/3}}

\begin{lemma}[{\protect\cite[Theorem 4.3]{TChenAndNP}}]
For $\alpha >\frac{1}{2}$, we have%
\begin{equation*}
\left\Vert R_{\alpha }^{(k)}B_{j,k+1}u^{(k+1)}\right\Vert _{L_{\mathbf{y,y}%
^{\prime }}^{2}}\leqslant C\left\Vert R_{\alpha
}^{(k+1)}u^{(k+1)}\right\Vert _{L_{\mathbf{y,y}^{\prime }}^{2}}.
\end{equation*}
\end{lemma}

\begin{itemize}
\item[Step IV] (Conclusion) By Step III, the compact sequence $\left\{
\Gamma _{N}(\tau )\right\} $ has only one limit point, thus it converges,
that is, as trace class operators kernels,%
\begin{equation*}
u_{N}^{(k)}(\tau )\rightarrow \dprod\limits_{j=1}^{k}\tilde{\phi}(\tau
,y_{j})\overline{\tilde{\phi}(\tau ,y_{j}^{\prime })}\text{ weak* as }%
N\rightarrow \infty ,\text{ }\forall \tau \in \left[ 0,\frac{\tan \omega
T_{0}}{\omega }\right] .
\end{equation*}%
Notice that the above weak* limit is an orthogonal projection, the argument
in the bottom of \cite[p. 296]{E-S-Y3} which uses the Gr\"{u}mm's
convergence theorem \cite[Theorem 2.19]{Simon}\footnote{%
One can also use the argument in \cite[Appendix A]{Chen3DDerivation}.} then
implies the strong convergence in trace norm%
\begin{equation*}
\lim_{N\rightarrow \infty }\limfunc{Tr}\left\vert u_{N}^{(k)}(\tau ,\mathbf{y%
}_{k},\mathbf{y}_{k}^{\prime })-\dprod\limits_{j=1}^{k}\tilde{\phi}(\tau
,y_{j})\overline{\tilde{\phi}(\tau ,y_{j}^{\prime })}\right\vert =0,\text{ }%
\forall \tau \in \left[ 0,\frac{\tan \omega T_{0}}{\omega }\right] .
\end{equation*}%
Recall $\gamma _{N}^{(k)}=T_{k}u_{N}^{(k)}$ and $\phi =M_{1}\tilde{\phi}$,
we utilize Lemma $\ref{Lemma:TraceNormPreservation}\ $and infer that%
\begin{equation*}
\lim_{N\rightarrow \infty }\limfunc{Tr}\left\vert \gamma _{N}^{(k)}(t,%
\mathbf{x}_{k},\mathbf{x}_{k}^{\prime })-\dprod\limits_{j=1}^{k}\phi
(t,x_{j})\overline{\phi (t,x_{j}^{\prime })}\right\vert =0,\text{ }\forall
t\in \left[ 0,T_{0}\right] ,
\end{equation*}%
where $\phi $ solves (\ref{equation:TargetCubicNLS}). So far, we have proved
Theorem \ref{THM:Main Theorem} for every $T_{0}<\frac{\pi }{2\omega }$, a
bootstrapping argument then establishes Theorem \ref{THM:Main Theorem} for
all time. Thence we conclude the proof of Theorem \ref{THM:Main Theorem}.
\end{itemize}

\section{Energy Estimates for the Focusing $N$-body Hamiltonian \label%
{Section:EnergyEstimate}}

\begin{theorem}
\label{THM:EnergyEstimate}Let $H_{N}$ be defined as in (\ref{Hamiltonian:1D
N-body}). For every $k,$ there exists $N_{0}(k)$ such that, we have 
\begin{equation*}
\langle \psi ,\left( H_{N}+N\left\Vert V\right\Vert _{L^{1}}^{2}+N\right)
^{k}\psi \rangle \geq 2^{-k}N^{k}\Vert S^{(k)}\psi \Vert _{L^{2}}^{2},
\end{equation*}%
for all $N\geq N_{0}(k)$ and $\psi \in L_{s}^{2}(\mathbb{R}^{N})$ with $%
\left\Vert \psi \right\Vert _{L^{2}}=1.$
\end{theorem}

We prove Theorem \ref{THM:EnergyEstimate} in \S \ref%
{Sec:ProofOfEnergyEstimate}. At the moment, we present the following
corollary of Theorem \ref{THM:EnergyEstimate}.

\begin{corollary}
\label{Corollary:EnergyBoundForMarginals}Let $\psi _{N}(t,\mathbf{x}%
_{N})=e^{itH_{N}}\psi _{N}(0)$ for some $\beta \in \left( 0,1\right) $
subject to initial $\psi _{N}(0)$ which satisfies energy condition (\ref%
{Condition:EnergyBoundOnManyBodyInitialData}). If $u_{N}(\tau ,\mathbf{y}%
_{N})=M_{N}^{-1}\psi _{N}$, where $M_{N}^{-1}$ is the inverse lens transform
for functions in Definition \ref{def:functionG-Lens}, then there is a $%
C\geqslant 0$, for all $k\geqslant 0$, there exists $N_{0}\left( k\right) $
such that 
\begin{equation*}
\left\langle u_{N}(\tau ),\dprod\limits_{j=1}^{k}L_{j}^{2}u_{N}(\tau
)\right\rangle \leqslant C^{k},
\end{equation*}%
for all $N\geqslant N_{0}$ and all $\tau \in \left[ -\frac{\tan \omega T_{0}%
}{\omega },\frac{\tan \omega T_{0}}{\omega }\right] $ provided that $T_{0}<%
\frac{\pi }{2\omega }.$ Thus, for $u_{N}^{(k)}=T_{k}^{-1}\gamma _{N}^{(k)}$,
the inverse lens transform of $\gamma _{N}^{(k)}$, 
\begin{equation*}
\sup_{\tau \in \left[ -\frac{\tan \omega T_{0}}{\omega },\frac{\tan \omega
T_{0}}{\omega }\right] }\limfunc{Tr}L^{(k)}u_{N}^{(k)}\left( \tau \right)
L^{(k)}\leqslant C^{k}
\end{equation*}%
where the inverse lens transform for kernels is given by Definition \ref%
{def:KernelG-Lens}.
\end{corollary}

\begin{proof}
By Lemma \ref{Lemma:EnergyBetweenLensTransform}, we have%
\begin{equation*}
\left\langle u_{N}(\tau ),\dprod\limits_{j=1}^{k}L_{j}^{2}u_{N}(\tau
)\right\rangle \leqslant C^{k}\left\langle \psi
_{N}(t),\dprod\limits_{j=1}^{k}S_{j}^{2}\psi _{N}(t)\right\rangle .
\end{equation*}%
With Theorem \ref{THM:EnergyEstimate} and the conservation law, we get%
\begin{eqnarray*}
\left\langle u_{N}(\tau ),\dprod\limits_{j=1}^{k}L_{j}^{2}u_{N}(\tau
)\right\rangle &\leqslant &\frac{C^{k}}{N^{k}}\langle \psi _{N}(t),\left(
H_{N}+N\left\Vert V\right\Vert _{L^{1}}^{2}+N\right) ^{k}\psi _{N}(t)\rangle
\\
&=&\frac{C^{k}}{N^{k}}\langle \psi _{N}(0),\left( H_{N}+N\left\Vert
V\right\Vert _{L^{1}}^{2}+N\right) ^{k}\psi _{N}(0)\rangle
\end{eqnarray*}%
The binomial theorem and (\ref{Condition:EnergyBoundOnManyBodyInitialData})
give 
\begin{eqnarray}
&&\langle \psi _{N}(0),\left( H_{N}+N\left\Vert V\right\Vert
_{L^{1}}^{2}+N\right) ^{k}\psi _{N}(0)\rangle
\label{estimate:why negative energy is fine} \\
&=&\sum_{j=0}^{k}%
\begin{pmatrix}
k \\ 
j%
\end{pmatrix}%
\left( N\left\Vert V\right\Vert _{L^{1}}^{2}+N\right) ^{j}\langle \psi
_{N}(0),H_{N}^{k-j}\psi _{N}(0)\rangle  \notag \\
&\leqslant &\sum_{j=0}^{k}%
\begin{pmatrix}
k \\ 
j%
\end{pmatrix}%
\left( N\left\Vert V\right\Vert _{L^{1}}^{2}+N\right) ^{j}C^{k-j}N^{k-j} 
\notag \\
&=&\left( CN+N\left\Vert V\right\Vert _{L^{1}}^{2}+N\right) ^{k}  \notag \\
&\leqslant &C^{k}N^{k}.  \notag
\end{eqnarray}%
Thus%
\begin{equation*}
\left\langle u_{N}(\tau ),\dprod\limits_{j=1}^{k}L_{j}^{2}u_{N}(\tau
)\right\rangle \leqslant \frac{C^{k}}{N^{k}}C^{k}N^{k}\leqslant C^{k},
\end{equation*}%
as claimed.
\end{proof}

\subsection{Proof of Theorem \protect\ref{THM:EnergyEstimate}\label%
{Sec:ProofOfEnergyEstimate}}

For convenience, we let $\alpha =\Vert V\Vert _{L^{1}}^{2}$ and rewrite the
desired estimate as%
\begin{equation}
\langle \psi ,\left( N^{-1}H_{N}+1+\alpha \right) ^{k}\psi \rangle \geqslant
2^{-k}\Vert S^{(k)}\psi \Vert _{L^{2}}^{2}.
\label{estimate:target in energy estimate}
\end{equation}%
Note that estimate (\ref{estimate:target in energy estimate}) is trivial for 
$k=0$. To establish estimate (\ref{estimate:target in energy estimate}) for
general $k$, we first prove the $k=1$ case which is already nontrivial in \S %
\ref{Sec:EnergyEstimate:k=1}, we then prove estimate (\ref{estimate:target
in energy estimate}) for $k+2$ assuming that it holds for $k$ in \S \ref%
{Sec:EnergyEstimate:k=k+2}, thence a two-step induction based on the $k=0$
and $k=1$ cases proves (\ref{estimate:target in energy estimate}) for all $k$%
.

The only technical tool we need is the 1D estimate: for $f(x)$ 
\begin{equation}
\Vert f\Vert _{L_{x}^{\infty }}\leq \Vert f^{\prime }\Vert _{L_{x}^{1}}
\label{E:1-120}
\end{equation}%
which is a direct consequence of the fundamental theorem of calculus. We
also utilize the ordinary Sobolev estimates: 
\begin{eqnarray}
\Vert f\Vert _{L_{x}^{\infty }} &\leqslant &C\Vert S_{x}f\Vert _{L_{x}^{2}}%
\text{ for }f(x),  \label{Estimate:1DSoblev} \\
\Vert f\Vert _{L_{xy}^{\infty }} &\leqslant &C\Vert S_{x}S_{y}f\Vert
_{L_{xy}^{2}}\text{ for }f(x,y).
\end{eqnarray}%
when the sizes of the controlling constants do not matter.\footnote{%
The only place in which we apply (\ref{E:1-120}) is the proof of Lemma \ref%
{L:1}. We use (\ref{E:1-120}) to determine $\alpha =\Vert V\Vert
_{L^{1}}^{2} $. With the elementary inequality: $\left\vert ab\right\vert
\leqslant \varepsilon a^{2}+\varepsilon ^{-1}b^{2},$ one can use (\ref%
{Estimate:1DSoblev}) instead and get another $\alpha $, namely, $\alpha
=C\left\Vert V\right\Vert _{L^{1}}^{2}$ for some $C$ depending on the
controlling constant in (\ref{Estimate:1DSoblev}). We are not using (\ref%
{Estimate:1DSoblev}) because we would like to give an exact $\alpha $ and
keep track of one less constant.} We will use the shorthand $%
L_{c}^{1}L_{x_{2}}^{\infty }$ for $L_{x_{1}x_{3}x_{4}\cdots
x_{N}}^{1}L_{x_{2}}^{\infty }$. Here, $c$ stands for \textquotedblleft
complementary coordinates\textquotedblright .

\subsubsection{The $k=1$ Case\label{Sec:EnergyEstimate:k=1}}

Recall $V_{N}(x)=N^{\beta }V(N^{\beta }x)$ and 
\begin{equation*}
S_{j}=\left( 1-\frac{1}{2}\partial _{x_{j}}^{2}+\frac{1}{2}\omega
^{2}x_{j}^{2}\right) ^{\frac{1}{2}}.
\end{equation*}
We write 
\begin{equation*}
H_{N}+N=\sum_{j=1}^{N}S_{j}^{2}+\frac{1}{2N}\sum_{\substack{ i,j=1,\ldots ,N 
\\ i\neq j}}V_{N}(x_{i}-x_{j}).
\end{equation*}%
We next introduce a convenient decomposition of $H_{N}+N$. Let%
\begin{equation*}
H_{ij}=S_{i}^{2}+S_{j}^{2}+\frac{N-1}{N}V_{N}(x_{i}-x_{j})
\end{equation*}%
Note that $H_{ij}=H_{ji}$ because $V$ is even, and 
\begin{equation*}
H_{N}+N=\frac{1}{2(N-1)}\sum_{\substack{ i,j=1,\ldots ,N  \\ i\neq j}}H_{ij}.
\end{equation*}%
It follows that 
\begin{equation}
N^{-1}H_{N}+1+\alpha =\frac{1}{2N(N-1)}\sum_{\substack{ i,j=1,\ldots ,N  \\ %
i\neq j}}(H_{ij}+2\alpha )  \label{E:1-101}
\end{equation}

\begin{lemma}
\label{L:1}Recall $\alpha =\Vert V\Vert _{L^{1}}^{2}$, we have 
\begin{equation*}
(H_{12}+2\alpha )\geqslant \tfrac{1}{2}(S_{1}^{2}+S_{2}^{2})
\end{equation*}
\end{lemma}

\begin{proof}
Apply the well-known change of variable $y_{1}=x_{1}-x_{2}$, $%
y_{2}=x_{1}+x_{2}$ which is also compatible with the Hermite operator, then 
\begin{eqnarray*}
H_{12} &=&2-\partial _{y_{1}}^{2}-\partial _{y_{2}}^{2}+\omega
^{2}\left\vert y_{1}\right\vert ^{2}+\omega ^{2}\left\vert y_{2}\right\vert
^{2}+(1-N^{-1})V_{N}(y_{1}) \\
&=&K_{y_{1}}+2-\partial _{y_{2}}^{2}+\omega ^{2}\left\vert y_{1}\right\vert
^{2}+\omega ^{2}\left\vert y_{2}\right\vert ^{2}
\end{eqnarray*}%
where 
\begin{equation*}
K_{y}=-\partial _{y}^{2}+(1-N^{-1})V_{N}(y).
\end{equation*}%
We claim that 
\begin{equation}
(K+2\alpha )\geq -\frac{1}{2}\partial _{y}^{2}  \label{E:1-103}
\end{equation}%
Indeed, 
\begin{align*}
\langle K\phi ,\phi \rangle & \geqslant \Vert \phi ^{\prime }\Vert
_{L^{2}}^{2}-\Vert V_{N}\Vert _{L^{1}}\Vert |\phi |^{2}\Vert _{L_{y}^{\infty
}} \\
& \geqslant \Vert \phi ^{\prime }\Vert _{L^{2}}^{2}-\Vert V\Vert
_{L^{1}}\Vert \partial _{y}(|\phi |^{2})\Vert _{L_{y}^{1}} \\
& \geqslant \Vert \phi ^{\prime }\Vert _{L^{2}}^{2}-2\Vert V\Vert
_{L^{1}}\Vert \phi ^{\prime }\Vert _{L^{2}}\Vert \phi \Vert _{L^{2}} \\
& \geqslant \Vert \phi ^{\prime }\Vert _{L^{2}}^{2}-(\frac{1}{2}\Vert \phi
^{\prime }\Vert _{L^{2}}^{2}+2\Vert V\Vert _{L^{1}}^{2}\Vert \phi \Vert
_{L^{2}}^{2}) \\
& =\frac{1}{2}\Vert \phi ^{\prime }\Vert _{L^{2}}^{2}-2\Vert V\Vert
_{L^{1}}^{2}
\end{align*}%
from which \eqref{E:1-103} follows.

We clearly have%
\begin{equation*}
2-\partial _{y_{2}}^{2}+\omega ^{2}\left\vert y_{1}\right\vert ^{2}+\omega
^{2}\left\vert y_{2}\right\vert ^{2}\geqslant 1-\frac{1}{2}\partial
_{y_{2}}^{2}+\frac{1}{2}\omega ^{2}\left\vert y_{1}\right\vert ^{2}+\frac{1}{%
2}\omega ^{2}\left\vert y_{2}\right\vert ^{2}.
\end{equation*}%
By this and \eqref{E:1-103}, we have 
\begin{eqnarray*}
H_{12}+2\alpha &=&(K_{y_{1}}+2\alpha +2-\partial _{y_{2}}^{2}+\omega
^{2}\left\vert y_{1}\right\vert ^{2}+\omega ^{2}\left\vert y_{2}\right\vert
^{2}) \\
&\geqslant &-\frac{1}{2}\partial _{y_{1}}^{2}+1-\frac{1}{2}\partial
_{y_{2}}^{2}+\frac{1}{2}\omega ^{2}\left\vert y_{1}\right\vert ^{2}+\frac{1}{%
2}\omega ^{2}\left\vert y_{2}\right\vert ^{2} \\
&=&1-\frac{1}{4}\partial _{x_{1}}^{2}-\frac{1}{4}\partial _{x_{2}}^{2}+\frac{%
1}{4}\omega ^{2}\left\vert x_{1}\right\vert ^{2}+\frac{1}{4}\omega
^{2}\left\vert x_{2}\right\vert ^{2} \\
&=&\tfrac{1}{2}(S_{1}^{2}+S_{2}^{2}).
\end{eqnarray*}
\end{proof}

In light of \eqref{E:1-101}, symmetry, and Lemma \ref{L:1}, we readily see
that 
\begin{eqnarray}
2\langle \psi ,(N^{-1}H_{N}+1+\alpha )\psi \rangle &=&\langle \psi
,(H_{12}+2\alpha )\psi \rangle  \label{E:1-106} \\
&\geqslant &\tfrac{1}{2}\langle \psi ,(S_{1}^{2}+S_{2}^{2})\psi \rangle 
\notag \\
&=&\Vert S_{1}\psi \Vert _{L^{2}}^{2}.  \notag
\end{eqnarray}%
Thus we have proved (\ref{estimate:target in energy estimate}) for $k=1$.

\subsubsection{The $k+2$ Case\label{Sec:EnergyEstimate:k=k+2}}

For convenience, let us introduce some notation. For any function $f$, let 
\begin{equation*}
f_{Nij}=N^{\beta }f(N^{\beta }(x_{i}-x_{j}))
\end{equation*}%
Also, let 
\begin{equation}
H_{+ij}=H_{ij}+2\alpha =S_{i}^{2}+S_{j}^{2}+(1-N^{-1})V_{Nij}+2\alpha
\label{E:1-125}
\end{equation}%
Then \eqref{E:1-101} can be written more compactly as 
\begin{equation}
N^{-1}H_{N}+1+\alpha =\frac{1}{2N(N-1)}\sum_{\substack{ i,j=1,\ldots ,N  \\ %
i\neq j}}H_{+ij}=\frac{1}{N(N-1)}\sum_{1\leq i<j\leq N}H_{+ij}
\label{E:1-102}
\end{equation}

Before delving into the proof of the $k+2$ case, we give an idea for why (%
\ref{estimate:target in energy estimate}) is true for all $k.$ Note that we
have 
\begin{equation*}
2^{k}(N^{-1}H_{N}+1+\alpha )^{k}=\frac{1}{N^{k}(N-1)^{k}}\sum_{\substack{ %
1\leq i_{1},j_{1},\ldots ,i_{k},j_{k}\leq N  \\ i_{1}\neq i_{2},\ldots
,i_{k}\neq j_{k}}}H_{+i_{1}j_{1}}\cdots H_{+i_{k}j_{k}}
\end{equation*}%
The dominant term in this expression occurs when all indices $%
i_{1},j_{1},\ldots ,i_{k},j_{k}$ are distinct, since it occurs with
frequency $\sim N^{2k}$. The other terms occur with lower frequency -- for
example, the terms in which exactly two of the indices are equal and all
others are distinct occur with frequency $\sim N^{2k-1}$. By symmetry, the
terms in which all indices are distinct can be rearranged so that formally,
we have 
\begin{equation}
2^{k}\langle \psi ,(N^{-1}H_{N}+1+\alpha )^{k}\psi \rangle \approx \langle
H_{+12}\cdots H_{+(2k-1)(2k)}\psi ,\psi \rangle  \label{E:1-105}
\end{equation}%
Moreover, by symmetry 
\begin{equation}
2^{-k}\langle \psi ,\dprod\limits_{i=1}^{k}(S_{2i-1}^{2}+S_{2i}^{2})\psi
\rangle =\Vert S^{(k)}\psi \Vert _{L^{2}}^{2}  \label{E:1-107}
\end{equation}%
Since $H_{+ij}\geq 2^{-1}\left( S_{i}^{2}+S_{j}^{2}\right) $ for each $i,j$
by Lemma \ref{L:1}, Lemma \ref{L:commuting-positive-ops} implies 
\begin{equation*}
H_{+12}\cdots H_{+(2k-1)(2k)}\geq
2^{-k}\dprod\limits_{i=1}^{k}(S_{2i-1}^{2}+S_{2i}^{2})
\end{equation*}%
This, together with \eqref{E:1-105} and \eqref{E:1-107} suggest (not
rigorously) that a statement like \eqref{estimate:target in energy
estimate} should hold.

We now establish \eqref{estimate:target in energy estimate} for $k+2$
rigorously, assuming it holds for $k$. To be precise, we will prove that, if %
\eqref{estimate:target in energy estimate} holds for $k$, then 
\begin{eqnarray}
&&2^{k+2}\langle \psi ,(N^{-1}H_{N}+\alpha +1)^{k+2}\psi \rangle
\label{E:1-104} \\
&\geqslant &\left( 1-C_{k+2}N^{\beta -1}\right) \left( \left\Vert
S^{(k+2)}\psi \right\Vert _{L^{2}}^{2}+N^{-1}\left\Vert S_{1}S^{(k+1)}\psi
\right\Vert _{L^{2}}^{2}\right)  \notag
\end{eqnarray}%
We remind the reader that we already have the $k=0$ case which is trivial
and the $k=1$ case proved in \S \ref{Sec:EnergyEstimate:k=1}, thus 
\eqref{estimate:target in energy
estimate} is proved for all $k$ once we prove that \eqref{E:1-104} holds as
long as \eqref{estimate:target in energy
estimate} is true for $k$.

Using the induction hypothesis, we arrive at%
\begin{eqnarray}
&&2^{k+2}\langle \psi ,(N^{-1}H_{N}+\alpha +1)^{k+2}\psi \rangle
\label{E:1-122} \\
&=&4(2^{k}\langle (N^{-1}H_{N}+\alpha +1)\psi ,(N^{-1}H_{N}+\alpha
+1)^{k}(N^{-1}H_{N}+\alpha +1)\psi \rangle )  \notag \\
&\geqslant &4\langle S^{(k)}(N^{-1}H_{N}+\alpha +1)\psi
,S^{(k)}(N^{-1}H_{N}+\alpha +1)\psi \rangle .  \notag
\end{eqnarray}%
We start with the following decomposition of the rightmost sum in %
\eqref{E:1-102}: 
\begin{equation}
(N^{-1}H_{N}+\alpha +1)=\frac{1}{N(N-1)}\sum_{\substack{ 1\leq i<j\leq N  \\ %
i\leq k}}H_{+ij}+\frac{1}{N(N-1)}\sum_{\substack{ 1\leq i<j\leq N  \\ i>k}}%
H_{+ij}.  \label{E:1-124}
\end{equation}%
Note that in the first term $i\leq k$ and $j$ can be either $\leq k$ or $>k$%
. We have ordered the indices $i_{1}<j_{1}$ and $i_{2}<j_{2}$ for
convenience. In the unordered setting, the above decomposition would be
characterized as follows: the first sum consists of terms in which \emph{at
least} one index is $\leq k$, and the second term consists of terms in which
both indices are $>k$. The decomposition \eqref{E:1-124} is similar to the
one used in the \cite[Proposition 1]{E-E-S-Y1}, although the authors of \cite%
{E-E-S-Y1} do not use the $H_{ij}$ decomposition of the Hamiltonian. There
are $\sim N$ terms in the first sum and $\sim N^{2}$ terms in the second
sum. Note that in the $k=0$ case, the decomposition \eqref{E:1-124} contains
only the second term since the first term is an empty sum.

Plug the decomposition \eqref{E:1-124} into the end of \eqref{E:1-122} to
obtain 
\begin{equation*}
2^{k+2}\langle \psi ,(N^{-1}H_{N}+\alpha +1)^{k+2}\psi \rangle \geq
A_{1}+A_{2}+A_{3}
\end{equation*}%
where 
\begin{equation*}
A_{1}=\frac{4}{N^{2}(N-1)^{2}}\sum_{\substack{ 1\leq i_{1}<j_{1}\leq N  \\ %
1\leq i_{2}<j_{2}\leq N  \\ \text{such that }i_{1}>k,\;i_{2}>k}}\langle
S^{(k)}H_{+i_{1}j_{1}}\psi ,\;S^{(k)}H_{+i_{2}j_{2}}\psi \rangle ,
\end{equation*}%
\begin{equation*}
A_{2}=\frac{4}{N^{2}(N-1)^{2}}\sum_{\substack{ 1\leq i_{1}<j_{1}\leq N  \\ %
1\leq i_{2}<j_{2}\leq N  \\ \text{such that }i_{1}\leq k,\;i_{2}>k}}2\func{Re%
}\langle S^{(k)}H_{+i_{1}j_{1}}\psi ,\;S^{(k)}H_{+i_{2}j_{2}}\psi \rangle ,
\end{equation*}%
\begin{eqnarray*}
A_{3} &=&\frac{4}{N^{2}(N-1)^{2}}\sum_{\substack{ 1\leq i_{1}<j_{1}\leq N 
\\ 1\leq i_{2}<j_{2}\leq N  \\ \text{such that }i_{1}\leq k,\;i_{2}\leq k}}%
\langle S^{(k)}H_{+i_{1}j_{1}}\psi ,\;S^{(k)}H_{+i_{2}j_{2}}\psi \rangle \\
&=&\frac{4}{N^{2}(N-1)^{2}}\langle S^{(k)}\sum_{\substack{ 1\leq i<j\leq N 
\\ i\leqslant k}}H_{+ij}\psi ,\;\sum_{\substack{ 1\leq i<j\leq N  \\ %
i\leqslant k}}S^{(k)}H_{+ij}\psi \rangle \geqslant 0.
\end{eqnarray*}%
Since $A_{3}\geq 0$, we drop this term to obtain 
\begin{equation}
2^{k+2}\langle \psi ,(N^{-1}H_{N}+\alpha +1)^{k+2}\psi \rangle \geq
A_{1}+A_{2}.  \label{E:1-145}
\end{equation}%
Note that $A_{1}$ contains $\sim N^{4}$ terms and the cross term $A_{2}$
contains $\sim N^{3}$ terms.\footnote{%
In the case $k=0$, $A_{2}=0$ and $A_{3}=0$ since they are both empty sums.}
In other words, $A_{1}$ is the dominant term and $A_{2}$ is the error term.
In below, we deal with $A_{1}$ and $A_{2}$ one by one.

In $A_{1}$, we can commute both terms $H_{+i_{1}j_{1}}$ and $H_{+i_{2}j_{2}}$
with $S^{(k)}$. Then 
\begin{equation}
A_{1}=\frac{4}{N^{2}(N-1)^{2}}\sum_{\substack{ 1\leq i_{1}<j_{1}\leq N  \\ %
1\leq i_{2}<j_{2}\leq N  \\ \text{such that }i_{1}>k,\;i_{2}>k}}\langle
S^{(k)}\psi ,\;H_{+i_{1}j_{1}}H_{+i_{2}j_{2}}S^{(k)}\psi \rangle
\label{E:1-126}
\end{equation}%
We decompose 
\begin{equation}
A_{1}=A_{11}+A_{12}+A_{13}  \label{E:1-152}
\end{equation}%
where

\begin{itemize}
\item $A_{11}$ consists of those terms for which all indices $i_1$, $j_1$, $%
i_2$, $j_2$ are different. There are $\frac14 a_{N,k} N^4$ such terms, where 
\begin{equation*}
a_{N,k} \overset{\mathrm{def}}{=} N^{-4}(N-k)(N-k-1)(N-k-2)(N-k-3)
\end{equation*}

\item $A_{12}$ consists of those terms for which exactly one pair of indices 
$i_1$, $j_1$, $i_2$, $j_2$ are the same. There are $b_N N^3$ such terms,
where 
\begin{equation*}
b_{N,k} \overset{\mathrm{def}}{=} N^{-3}(N-k)(N-k-1)(N-k-2)
\end{equation*}

\item $A_{13}$ consists of those terms for which exactly two pairs of
indices $i_1$, $j_1$, $i_2$, $j_2$ are the same. There are $\frac12 c_N N^2$
such terms, where 
\begin{equation*}
c_{N,k} \overset{\mathrm{def}}{=} N^{-2}(N-k)(N-k-1)
\end{equation*}
\end{itemize}

Note that 
\begin{equation*}
1-C_{k}N^{-1}\leq a_{N,k}\leq 1+C_{k}N^{-1}
\end{equation*}%
for some\footnote{%
We allow that $C_{k}$ changes from one line to the next.} $C_{k}$, and
similarly for $b_{N,k}$, $c_{N,k}$. Likewise, the coefficient in %
\eqref{E:1-126} satisfies 
\begin{equation*}
4N^{-4}(1-CN^{-1})\leq \frac{4}{N^{2}(N-1)^{2}}\leq 4N^{-4}(1+CN^{-1})
\end{equation*}%
The $O(N^{-1})$ corrections are easily absorbed into the error term in %
\eqref{E:1-104} and we drop them in the calculations that follow, for
expositional convenience.

By symmetry, we have 
\begin{equation*}
A_{11}=\langle H_{+(k+1)(k+2)}S^{(k)}\psi ,H_{+(k+3)(k+4)}S^{(k)}\psi \rangle
\end{equation*}%
\begin{equation*}
A_{12}=4N^{-1}\langle H_{+(k+1)(k+2)}S^{(k)}\psi ,H_{+(k+2)(k+3)}S^{(k)}\psi
\rangle
\end{equation*}%
\begin{equation*}
A_{13}=2N^{-2}\langle H_{+(k+1)(k+2)}S^{(k)}\psi ,H_{+(k+1)(k+2)}S^{(k)}\psi
\rangle
\end{equation*}%
Since 
\begin{equation}
A_{13}\geqslant 0  \label{E:1-151}
\end{equation}%
we can discard it. By Lemmas \ref{L:1} and \ref{L:commuting-positive-ops},
we have 
\begin{equation*}
A_{11}\geqslant \tfrac{1}{4}\langle \left( S_{k+1}^{2}+S_{k+2}^{2}\right)
S^{(k)}\psi ,\left( S_{k+3}^{2}+S_{k+4}^{2}\right) S^{(k)}\psi \rangle
\end{equation*}%
By integration by parts and symmetry, we obtain 
\begin{equation}
A_{11}\geqslant \Vert S^{(k+2)}\psi \Vert _{L^{2}}^{2}.  \label{E:1-146}
\end{equation}

Plugging in the definition \eqref{E:1-125} of $H_{+ij}$ and expanding 
\begin{equation}
A_{12}=A_{121}+A_{122}+A_{123}  \label{E:1-147}
\end{equation}%
where 
\begin{equation*}
A_{121}=4N^{-1}\langle \left( S_{k+1}^{2}+S_{k+2}^{2}\right) \left(
S_{k+2}^{2}+S_{k+3}^{2}\right) S^{(k)}\psi ,S^{(k)}\psi \rangle
\end{equation*}%
\begin{equation*}
A_{122}=8\func{Re}N^{-1}\langle \left( S_{k+1}^{2}+S_{k+2}^{2}\right)
(V_{N(k+2)(k+3)}+2\alpha )S^{(k)}\psi ,S^{(k)}\psi \rangle
\end{equation*}%
\begin{equation*}
A_{123}=4N^{-1}\langle (V_{N(k+1)(k+2)}+2\alpha )(V_{N(k+2)(k+3)}+2\alpha
)S^{(k)}\psi ,S^{(k)}\psi \rangle
\end{equation*}

For $A_{121}$, we only need to keep one term: 
\begin{equation}
A_{121}\geqslant 4N^{-1}\Vert S_{k+2}^{2}S^{(k)}\psi \Vert
_{L^{2}}^{2}=4N^{-1}\Vert S_{1}S^{(k+1)}\psi \Vert _{L^{2}}^{2}.
\label{E:1-148}
\end{equation}

For $A_{122}$, integration by parts gives 
\begin{eqnarray*}
A_{122} &=&8N^{-1}\langle (V_{N(k+2)(k+3)}+2\alpha )S^{(k+1)}\psi
,S^{(k+1)}\psi \rangle \\
&&+8N^{-1}\langle (V_{N(k+2)(k+3)}+2\alpha )S_{k+2}S^{(k)}\psi
,S_{k+2}S^{(k)}\psi \rangle \\
&&+8\func{Re}N^{\beta -1}\langle \left( V^{\prime }\right)
_{N(k+2)(k+3)}S^{(k)}\psi ,S_{k+2}S^{(k)}\psi \rangle
\end{eqnarray*}%
where we used the fact that $\partial _{x_{j}}$ is the only thing inside $%
S_{j}$ which needs the Leibniz's rule. Estimating, 
\begin{eqnarray*}
\left\vert A_{122}\right\vert &\lesssim &N^{-1}\left\Vert
V_{N(k+2)(k+3)}\right\Vert _{L_{x_{k+3}}^{1}}\left\Vert S^{(k+1)}\psi
\right\Vert _{L_{c}^{2}L_{x_{k+3}}^{\infty }}^{2}+N^{-1}\alpha \Vert
S^{(k+1)}\psi \Vert _{L^{2}}^{2} \\
&&+N^{-1}\left\Vert V_{N(k+2)(k+3)}\right\Vert _{L_{x_{k+3}}^{1}}\left\Vert
S_{k+2}S^{(k)}\psi \right\Vert _{L_{c}^{2}L_{x_{k+3}}^{\infty
}}^{2}+N^{-1}\alpha \Vert S_{k+2}S^{(k)}\psi \Vert _{L^{2}}^{2} \\
&&+N^{\beta -1}\left\Vert \left( V^{\prime }\right)
_{N(k+2)(k+3)}\right\Vert _{L_{x_{k+3}}^{1}}\left\Vert S^{(k)}\psi
\right\Vert _{L_{c}^{2}L_{x_{k+3}}^{\infty }}\left\Vert S_{k+2}S^{(k)}\psi
\right\Vert _{L_{c}^{2}L_{x_{k+3}}^{\infty }} \\
&\lesssim &N^{\beta -1}\left( \left\Vert S^{(k+1)}\psi \right\Vert
_{L_{c}^{2}L_{x_{k+3}}^{\infty }}^{2}+\Vert S^{(k+1)}\psi \Vert
_{L^{2}}^{2}+\left\Vert S^{(k)}\psi \right\Vert
_{L_{c}^{2}L_{x_{k+3}}^{\infty }}^{2}\right) .
\end{eqnarray*}%
By the 1D estimate (\ref{Estimate:1DSoblev}), Cauchy-Schwarz, and symmetry,
we have 
\begin{eqnarray}
\left\vert A_{122}\right\vert &\leqslant &CN^{\beta -1}\left( \Vert
S^{(k+2)}\psi \Vert _{L^{2}}^{2}+\Vert S^{(k+1)}\psi \Vert
_{L^{2}}^{2}\right)  \label{E:1-149} \\
&\leqslant &CN^{\beta -1}\Vert S^{(k+2)}\psi \Vert _{L^{2}}^{2}.  \notag
\end{eqnarray}

For $A_{123}$, 
\begin{eqnarray*}
\left\vert A_{123}\right\vert &\leqslant &CN^{-1}\left\Vert
V_{N(k+1)(k+2)}\right\Vert _{L_{x_{k+1}}^{1}}\left\Vert
V_{N(k+2)(k+3)}\right\Vert _{L_{x_{k+3}}^{1}}\left\Vert S^{(k)}\psi
\right\Vert _{L_{c}^{2}L_{x_{k+1}}^{\infty }L_{x_{k+3}}^{\infty }}^{2} \\
&&+CN^{-1}\left\Vert V_{N(k+1)(k+2)}\right\Vert _{L_{x_{k+1}}^{1}}\left\Vert
S^{(k)}\psi \right\Vert _{L_{c}^{2}L_{x_{k+1}}^{\infty
}}^{2}+CN^{-1}\left\Vert S^{(k)}\psi \right\Vert _{L^{2}}^{2}.
\end{eqnarray*}%
Using (\ref{Estimate:1DSoblev}) twice, we obtain 
\begin{eqnarray}
\left\vert A_{123}\right\vert &\leqslant &CN^{-1}(\Vert S^{(k+2)}\psi \Vert
_{L^{2}}^{2}+\Vert S^{(k+1)}\psi \Vert _{L^{2}}^{2}+\left\Vert S^{(k)}\psi
\right\Vert _{L^{2}}^{2})  \label{E:1-150} \\
&\leqslant &CN^{-1}\Vert S^{(k+2)}\psi \Vert _{L^{2}}^{2}.  \notag
\end{eqnarray}

By \eqref{E:1-147}, \eqref{E:1-148}, \eqref{E:1-149}, and \eqref{E:1-150}, 
\begin{equation}
A_{12}\geqslant 4N^{-1}\Vert S_{1}S^{(k+1)}\psi \Vert _{L^{2}}^{2}-CN^{\beta
-1}\Vert S^{(k+2)}\psi \Vert _{L^{2}}^{2}.  \label{E:1-153}
\end{equation}

Collecting \eqref{E:1-152}, \eqref{E:1-151}, \eqref{E:1-146}, and %
\eqref{E:1-153}, we have the estimate for $A_{1}$: 
\begin{equation}
A_{1}\geqslant \left( 1-CN^{\beta -1}\right) \left( \Vert S^{(k+2)}\psi
\Vert _{L^{2}}^{2}+4N^{-1}\Vert S_{1}S^{(k+1)}\psi \Vert _{L^{2}}^{2}\right)
\label{E:1-154}
\end{equation}%
The above estimate yields the positive contribution on the right-side of %
\eqref{E:1-104}.

Next we turn our attention to estimating $A_{2}$. We will prove that%
\begin{equation*}
A_{2}\geqslant -CN^{\beta -1}\left( \Vert S^{(k+2)}\psi \Vert
_{L^{2}}^{2}+N^{-1}\Vert S_{1}S^{(k+1)}\psi \Vert _{L^{2}}^{2}\right) .
\end{equation*}%
Recall that in the case $k=0$, $A_{2}=0$, so we can assume $k\geq 1$. We
decompose 
\begin{equation}
A_{2}=A_{21}+A_{22}+A_{23}  \label{E:1-155}
\end{equation}%
where

\begin{itemize}
\item $A_{21}$ contains those terms with $j_1 \leq k$. There are $\sim N^2$
such terms. (In the case $k=1$, there are no terms of this type, so $%
A_{21}=0 $)

\item $A_{22}$ contains those terms with $j_{1}>k$, and ($j_{1}=i_{2}$ OR $%
j_{1}=j_{2}$). There are $\sim N^{2}$ such terms.

\item $A_{23}$ contains those terms with $j_1>k$, $j_1 \neq i_2$ and $%
j_1\neq j_2$. There are $\sim N^3$ such terms.
\end{itemize}

By symmetry of $\psi $ and $H_{+ij}=H_{+ji}$, 
\begin{equation*}
A_{21}=N^{-2}\langle S^{(k)}H_{+12}\psi ,S^{(k)}H_{+(k+1)(k+2)}\psi \rangle
\end{equation*}%
\begin{equation*}
A_{22}=N^{-2}\langle S^{(k)}H_{+1(k+1)}\psi ,S^{(k)}H_{+(k+1)(k+2)}\psi
\rangle
\end{equation*}%
\begin{equation*}
A_{23}=N^{-1}\langle S^{(k)}H_{+1(k+1)}\psi ,S^{(k)}H_{+(k+2)(k+3)}\psi
\rangle
\end{equation*}

First, we address $A_{21}.$ We decompose 
\begin{equation}
A_{21}=A_{211}+A_{212}+A_{213},  \label{E:1-156}
\end{equation}%
where 
\begin{equation*}
A_{211}=N^{-2}\langle H_{+12}S^{(k)}\psi ,H_{+(k+1)(k+2)}S^{(k)}\psi \rangle
\end{equation*}%
\begin{equation*}
A_{212}=N^{-2}\langle \lbrack S_{1},H_{+12}]S_{2}\cdots S_{k}\psi
,H_{+(k+1)(k+2)}S^{(k)}\psi \rangle
\end{equation*}%
\begin{eqnarray*}
A_{213} &=&N^{-2}\langle S_{1}[S_{2},H_{+12}]S_{3}\cdots S_{k}\psi
,H_{+(k+1)(k+2)}S^{(k)}\psi \rangle \\
&=&N^{-2}\langle \lbrack S_{2},H_{+12}]S_{3}\cdots S_{k}\psi
,H_{+(k+1)(k+2)}S_{1}S^{(k)}\psi \rangle .
\end{eqnarray*}

By Lemmas \ref{L:1} and \ref{L:commuting-positive-ops}, 
\begin{equation}
A_{211}\geqslant 0.  \label{E:1-157}
\end{equation}%
Since $[S_{1},H_{+12}]=N^{\beta }(V^{\prime })_{N12}$, integrating by parts
half the Hermite terms in $H_{+(k+1)(k+2)}$ and using symmetry, 
\begin{eqnarray*}
A_{212} &=&2N^{\beta -2}\langle (V^{\prime })_{N12}S_{2}\cdots
S_{k}S_{k+1}\psi ,S^{(k+1)}\psi \rangle \\
&&+2\alpha N^{\beta -2}\langle (V^{\prime })_{N12}S_{2}\cdots S_{k}\psi
,S^{(k)}\psi \rangle \\
&&+N^{\beta -2}\langle (V)_{N\left( k+1\right) \left( k+2\right) }(V^{\prime
})_{N12}S_{2}\cdots S_{k}\psi ,S^{(k)}\psi \rangle
\end{eqnarray*}%
Estimating 
\begin{eqnarray*}
\left\vert A_{212}\right\vert &\leqslant &CN^{\frac{3\beta }{2}-2}\left\Vert
V^{\prime }\right\Vert _{L_{x_{1}}^{2}}\left\Vert S_{2}\cdots
S_{k}S_{k+1}\psi \right\Vert _{L_{c}^{2}L_{x_{1}}^{\infty }}\left\Vert
S^{(k+1)}\psi \right\Vert _{L^{2}} \\
&&+CN^{\frac{3\beta }{2}-2}\left\Vert V^{\prime }\right\Vert
_{L_{x_{1}}^{2}}\left\Vert S_{2}\cdots S_{k}\psi \right\Vert
_{L_{c}^{2}L_{x_{1}}^{\infty }}\left\Vert S^{(k)}\psi \right\Vert _{L^{2}} \\
&&+CN^{\frac{3\beta }{2}-2}\left\Vert V\right\Vert
_{L_{x_{k+1}}^{1}}\left\Vert V^{\prime }\right\Vert
_{L_{x_{1}}^{2}}\left\Vert S_{2}\cdots S_{k}\psi \right\Vert
_{L_{c}^{2}L_{x_{1}}^{\infty }L_{x_{k+1}}^{\infty }}\left\Vert S^{(k)}\psi
\right\Vert _{L_{c}^{2}L_{x_{k+1}}^{\infty }}
\end{eqnarray*}%
Using (\ref{Estimate:1DSoblev}) and symmetry, 
\begin{eqnarray}
\left\vert A_{212}\right\vert &\leqslant &CN^{\frac{3\beta }{2}-2}\left(
\left\Vert S^{(k+1)}\psi \right\Vert _{L^{2}}^{2}+\left\Vert S^{(k)}\psi
\right\Vert _{L^{2}}^{2}+\left\Vert S^{(k+1)}\psi \right\Vert
_{L^{2}}^{2}\right)  \label{E:1-158} \\
&\leqslant &CN^{\frac{3\beta }{2}-2}\left\Vert S^{(k+1)}\psi \right\Vert
_{L^{2}}^{2}  \notag
\end{eqnarray}

For $A_{213}$, we use $[S_{2},H_{+12}]=-N^{\beta }(V^{\prime })_{N12}$ to
get 
\begin{equation*}
A_{213}=-N^{\beta -2}\langle (V^{\prime })_{N12}S_{3}\cdots S_{k}\psi
,H_{+(k+1)(k+2)}S_{1}S^{(k)}\psi \rangle
\end{equation*}%
Split up the terms of $H_{+(k+1)(k+2)}$ via integration by parts and use
symmetry to obtain 
\begin{eqnarray*}
A_{213} &=&-2N^{\beta -2}\langle (V^{\prime })_{N12}S_{3}\cdots
S_{k}S_{k+1}\psi ,S_{1}S^{(k+1)}\psi \rangle \\
&&-2\alpha N^{\beta -2}\langle (V^{\prime })_{N12}S_{3}\cdots S_{k}\psi
,S_{1}S^{(k)}\psi \rangle \\
&&-N^{\beta -2}\langle V_{N(k+1)(k+2)}(V^{\prime })_{N12}S_{3}\cdots
S_{k}\psi ,S_{1}S^{(k)}\psi \rangle
\end{eqnarray*}%
We now implement the same estimates used to treat $A_{212}$ but carry a
factor $N^{-1/2}$ with $S_{1}S^{(k)}\psi $ and $S_{1}S^{(k+1)}\psi $ 
\begin{eqnarray*}
\left\vert A_{213}\right\vert &\leqslant &CN^{\frac{3\beta }{2}-\frac{3}{2}%
}\left\Vert V^{\prime }\right\Vert _{L_{x_{1}}^{2}}\left\Vert S_{3}\cdots
S_{k}S_{k+1}\psi \right\Vert _{L_{c}^{2}L_{x_{1}}^{\infty
}}N^{-1/2}\left\Vert S_{1}S^{(k+1)}\psi \right\Vert _{L^{2}} \\
&&+CN^{\frac{3\beta }{2}-\frac{3}{2}}\left\Vert V^{\prime }\right\Vert
_{L_{x_{1}}^{2}}\left\Vert S_{3}\cdots S_{k}\psi \right\Vert
_{L_{c}^{2}L_{x_{1}}^{\infty }}N^{-1/2}\left\Vert S_{1}S^{(k)}\psi
\right\Vert _{L^{2}} \\
&&+CN^{\frac{3\beta }{2}-\frac{3}{2}}\left\Vert V\right\Vert
_{L_{x_{k+1}}^{2}}\left\Vert V^{\prime }\right\Vert
_{L_{x_{1}}^{2}}\left\Vert S_{3}\cdots S_{k}\psi \right\Vert
_{L_{c}^{2}L_{x_{1}}^{\infty }L_{x_{k+1}}^{\infty }}N^{-1/2}\left\Vert
S_{1}S^{(k)}\psi \right\Vert _{L_{c}^{2}L_{x_{k+1}}^{\infty }}
\end{eqnarray*}%
Arguing as above using (\ref{Estimate:1DSoblev}) and symmetry 
\begin{eqnarray}
\left\vert A_{213}\right\vert &\leqslant &CN^{\frac{3\beta }{2}-\frac{3}{2}%
}\left\Vert S^{(k)}\psi \right\Vert _{L^{2}}N^{-1/2}\left\Vert
S_{1}S^{(k+1)}\psi \right\Vert _{L^{2}}  \label{E:1-159} \\
&&+CN^{\frac{3\beta }{2}-\frac{3}{2}}\left\Vert S^{(k-1)}\psi \right\Vert
_{L^{2}}N^{-1/2}\left\Vert S_{1}S^{(k)}\psi \right\Vert _{L^{2}}  \notag \\
&&+CN^{\frac{3\beta }{2}-\frac{3}{2}}\left\Vert S^{(k)}\psi \right\Vert
_{L^{2}}N^{-1/2}\left\Vert S_{1}S^{(k+1)}\psi \right\Vert _{L^{2}}  \notag \\
&\leqslant &CN^{\frac{3\beta }{2}-\frac{3}{2}}\left( \left\Vert S^{(k)}\psi
\right\Vert _{L^{2}}^{2}+N^{-1}\left\Vert S_{1}S^{(k+1)}\psi \right\Vert
_{L^{2}}^{2}\right)  \notag
\end{eqnarray}

By \eqref{E:1-156}, \eqref{E:1-157}, \eqref{E:1-158}, and \eqref{E:1-159},
we obtain 
\begin{equation}
A_{21}\geqslant -CN^{\frac{3\beta }{2}-\frac{3}{2}}\left( \left\Vert
S^{(k+1)}\psi \right\Vert _{L^{2}}^{2}+N^{-1}\left\Vert S_{1}S^{(k+1)}\psi
\right\Vert _{L^{2}}^{2}\right)  \label{Estimate:A_21Final}
\end{equation}

Next, we address $A_{22}$. Recall 
\begin{equation*}
A_{22}=N^{-2}\langle S^{(k)}H_{+1(k+1)}\psi ,H_{+(k+1)(k+2)}S^{(k)}\psi
\rangle
\end{equation*}%
Decompose 
\begin{equation}
A_{22}=A_{221}+A_{222}  \label{E:1-142}
\end{equation}%
where 
\begin{equation*}
A_{221}=N^{-2}\langle H_{+1(k+1)}S^{(k)}\psi ,H_{+(k+1)(k+2)}S^{(k)}\psi
\rangle
\end{equation*}%
\begin{equation*}
A_{222}=N^{-2}\langle \lbrack S_{1},H_{+1(k+1)}]S_{2}...S_{k}\psi
,H_{+(k+1)(k+2)}S^{(k)}\psi \rangle
\end{equation*}

For $A_{221}$, plug in the definition \eqref{E:1-125} of $H_{+ij}$ to obtain
the decomposition 
\begin{equation*}
A_{221}=A_{2211}+A_{2212}+A_{2213}+A_{2214}
\end{equation*}%
where 
\begin{equation*}
A_{2211}=N^{-2}\langle \left( S_{1}^{2}+S_{k+1}^{2}\right) S^{(k)}\psi
,\left( S_{k+1}^{2}+S_{k+2}^{2}\right) S^{(k)}\psi \rangle
\end{equation*}%
\begin{equation*}
A_{2212}=N^{-2}\langle \left( S_{1}^{2}+S_{k+1}^{2}\right) S^{(k)}\psi
,(V_{N(k+1)(k+2)}+2\alpha )S^{(k)}\psi \rangle
\end{equation*}%
\begin{equation*}
A_{2213}=N^{-2}\langle (V_{N1(k+1)}+2\alpha )S^{(k)}\psi ,\left(
S_{k+1}^{2}+S_{k+2}^{2}\right) S^{(k)}\psi \rangle
\end{equation*}%
\begin{equation*}
A_{2214}=N^{-2}\langle (V_{N1(k+1)}+2\alpha )S^{(k)}\psi
,(V_{N(k+1)(k+2)}+2\alpha )S^{(k)}\psi \rangle
\end{equation*}%
Note that $A_{2211}\geq 0$, so we can discard this term. Integrating by
parts, 
\begin{eqnarray*}
A_{2212} &=&N^{-2}\langle S_{1}S^{(k)}\psi ,(V_{N(k+1)(k+2)}+2\alpha
)S_{1}S^{(k)}\psi \rangle \\
&&+N^{-2}\langle S^{(k+1)}\psi ,(V_{N(k+1)(k+2)}+2\alpha )S^{(k+1)}\psi
\rangle \\
&&+N^{\beta -2}\langle S^{(k+1)}\psi ,\left( V^{\prime }\right)
_{N(k+1)(k+2)}S^{(k)}\psi \rangle
\end{eqnarray*}%
Putting every instance of $V$ or $V^{\prime }$ in $L^{\infty }$, we obtain
the estimate 
\begin{eqnarray*}
\left\vert A_{2212}\right\vert &\leqslant &CN^{\beta -2}\left\Vert
S_{1}S^{(k)}\psi \right\Vert _{L^{2}}^{2}+CN^{\beta -2}\left\Vert
S^{(k+1)}\psi \right\Vert _{L^{2}}^{2} \\
&&+CN^{2\beta -2}\left\Vert S^{(k+1)}\psi \right\Vert _{L^{2}}\left\Vert
S^{(k)}\psi \right\Vert _{L^{2}}
\end{eqnarray*}%
Using that $\max \left( N^{2\beta -2},N^{\beta -2}\right) \leq N^{\beta -1}$%
, 
\begin{equation*}
\left\vert A_{2212}\right\vert \leqslant CN^{\beta -1}(N^{-1}\left\Vert
S_{1}S^{(k)}\psi \right\Vert _{L^{2}}^{2}+\left\Vert S^{(k+1)}\psi
\right\Vert _{L^{2}}^{2})
\end{equation*}%
By integration by parts, 
\begin{eqnarray*}
A_{2213} &=&N^{-2}\langle (V_{N1(k+1)}+2\alpha )S^{(k+1)}\psi ,S^{(k+1)}\psi
\rangle \\
&&-N^{\beta -2}\langle \left( V^{\prime }\right) _{N1(k+1)}S^{(k)}\psi
,S^{(k+1)}\psi \rangle \\
&&+N^{-2}\langle (V_{N1(k+1)}+2\alpha )S_{k+2}S^{(k)}\psi
,S_{k+2}S^{(k)}\psi \rangle
\end{eqnarray*}%
Putting every instance of $V$ or $V^{\prime }$ in $L^{\infty }$, 
\begin{eqnarray*}
\left\vert A_{2213}\right\vert &\leqslant &CN^{\beta -2}\left\Vert
S^{(k+1)}\psi \right\Vert _{L^{2}}^{2}+CN^{2\beta -2}\left\Vert
S^{(k+1)}\psi \right\Vert _{L^{2}}\left\Vert S^{(k)}\psi \right\Vert _{L^{2}}
\\
&&+CN^{\beta -2}\left\Vert S^{(k+1)}\psi \right\Vert _{L^{2}}^{2} \\
&\leqslant &CN^{2\beta -2}\left\Vert S^{(k+1)}\psi \right\Vert _{L^{2}}^{2}
\end{eqnarray*}%
For $A_{2214}$, we put both $V$ terms in $L^{\infty }$ to obtain 
\begin{equation*}
\left\vert A_{2214}\right\vert \leqslant CN^{2\beta -2}\Vert S^{(k)}\psi
\Vert _{L^{2}}^{2}
\end{equation*}%
This completes the bound for $A_{221}$. Specifically, 
\begin{equation}
A_{221}\geqslant -CN^{\beta -1}(\left\Vert S^{(k+1)}\psi \right\Vert
_{L^{2}}^{2}+N^{-1}\left\Vert S_{1}S^{(k)}\psi \right\Vert _{L^{2}}^{2})
\label{E:1-140}
\end{equation}

For $A_{222}$, substitute $[S_{1},H_{+1(k+1)}]=N^{\beta }(V^{\prime
})_{N1(k+1)}$ and plug in the definition \eqref{E:1-125} of $H_{+(k+1)(k+2)}$
to obtain%
\begin{equation*}
A_{222}=A_{2221}+A_{2222}+A_{2223}
\end{equation*}%
where 
\begin{equation*}
A_{2221}=N^{\beta -2}\langle (V^{\prime })_{N1(k+1)}S_{2}...S_{k}\psi
,S_{k+1}^{2}S^{(k)}\psi \rangle
\end{equation*}%
\begin{equation*}
A_{2222}=N^{\beta -2}\langle (V^{\prime })_{N1(k+1)}S_{2}...S_{k}\psi
,S_{k+2}^{2}S^{(k)}\psi \rangle
\end{equation*}%
\begin{equation*}
A_{2223}=N^{\beta -2}\langle (V^{\prime })_{N1(k+1)}S_{2}...S_{k}\psi
,(V_{N(k+1)(k+2)}+2\alpha )S^{(k)}\psi \rangle
\end{equation*}%
For $A_{2221}$, we apply H\"{o}lder in $x_{1}$ as follows: 
\begin{equation*}
\left\vert A_{2221}\right\vert \leqslant N^{\beta -2}\Vert (V^{\prime
})_{N1(k+1)}\Vert _{L_{x_{1}}^{2}}\Vert S_{2}...S_{k}\psi \Vert
_{L_{c}^{2}L_{x_{1}}^{\infty }}\Vert S_{k+1}^{2}S^{(k)}\psi \Vert _{L^{2}}
\end{equation*}%
By (\ref{Estimate:1DSoblev}) and symmetry, 
\begin{eqnarray*}
\left\vert A_{2221}\right\vert &\leqslant &CN^{\frac{3}{2}\beta -\frac{3}{2}%
}\Vert V^{\prime }\Vert _{L^{2}}\Vert S^{(k)}\psi \Vert
_{L^{2}}N^{-1/2}\Vert S_{1}S^{(k+1)}\psi \Vert _{L^{2}} \\
&\leqslant &CN^{\frac{3}{2}\beta -\frac{3}{2}}\left( \Vert S^{(k)}\psi \Vert
_{L^{2}}^{2}+N^{-1}\Vert S_{1}S^{(k+1)}\psi \Vert _{L^{2}}^{2}\right) .
\end{eqnarray*}%
Argue the same for $A_{2222}$, we get%
\begin{eqnarray*}
\left\vert A_{2222}\right\vert &\leqslant &N^{\beta -2}\Vert (V^{\prime
})_{N1(k+1)}\Vert _{L_{x_{1}}^{2}}\Vert S_{2}...S_{k}\psi \Vert
_{L_{c}^{2}L_{x_{1}}^{\infty }}\Vert S_{k+2}^{2}S^{(k)}\psi \Vert _{L^{2}} \\
&\leqslant &CN^{\frac{3}{2}\beta -\frac{3}{2}}\Vert V^{\prime }\Vert
_{L^{2}}\Vert S^{(k)}\psi \Vert _{L^{2}}N^{-1/2}\Vert S_{1}S^{(k+1)}\psi
\Vert _{L^{2}} \\
&\leqslant &CN^{\frac{3}{2}\beta -\frac{3}{2}}\left( \Vert S^{(k)}\psi \Vert
_{L^{2}}^{2}+N^{-1}\Vert S_{1}S^{(k+1)}\psi \Vert _{L^{2}}^{2}\right) .
\end{eqnarray*}%
For $A_{2223}$, we use H\"{o}lder in $x_{k+1}$ to obtain%
\begin{equation*}
N^{\beta -2}\langle (V^{\prime })_{N1(k+1)}S_{2}...S_{k}\psi
,(V_{N(k+1)(k+2)}+2\alpha )S^{(k)}\psi \rangle
\end{equation*}%
\begin{eqnarray*}
\left\vert A_{2223}\right\vert &\leqslant &CN^{\beta -2}\left\Vert
(V^{\prime })_{N1(k+1)}\right\Vert _{L_{x_{k+1}}^{1}}\left\Vert
V_{N(k+1)(k+2)}+2\alpha \right\Vert _{L_{x_{k+1}}^{\infty }} \\
&&\times \left\Vert S_{2}...S_{k}\psi \right\Vert
_{L_{c}^{2}L_{x_{k+1}}^{\infty }}\left\Vert S^{(k)}\psi \right\Vert
_{L_{c}^{2}L_{x_{k+1}}^{\infty }} \\
&\leqslant &CN^{2\beta -2}\left\Vert S^{(k)}\psi \right\Vert
_{L^{2}}\left\Vert S^{(k+1)}\psi \right\Vert _{L^{2}} \\
&\leqslant &CN^{2\beta -2}\left\Vert S^{(k+1)}\psi \right\Vert _{L^{2}}^{2}.
\end{eqnarray*}%
This completes the estimate for $A_{222}$. Specifically, collecting the
estimates for $A_{2221}\sim A_{2223}$, we obtain 
\begin{equation}
A_{222}\geqslant -CN^{\beta -1}(\left\Vert S^{(k+1)}\psi \right\Vert
_{L^{2}}^{2}+N^{-1}\Vert S_{1}S^{(k+1)}\psi \Vert _{L^{2}}^{2})
\label{E:1-141}
\end{equation}%
By \eqref{E:1-142}, \eqref{E:1-140} and \eqref{E:1-141}, we complete the
estimate for $A_{22}$ as 
\begin{equation}
A_{22}\geq -CN^{\beta -1}(\left\Vert S^{(k+1)}\psi \right\Vert
_{L^{2}}^{2}+N^{-1}\Vert S_{1}S^{(k+1)}\psi \Vert _{L^{2}}^{2})
\label{Estimate:A_22Final}
\end{equation}

Finally, for $A_{23}$, we have%
\begin{equation*}
A_{23}=N^{-1}\langle S^{(k)}H_{+1(k+1)}\psi ,S^{(k)}H_{+(k+2)(k+3)}\psi
\rangle
\end{equation*}%
\begin{equation}
A_{23}=A_{231}+A_{232}  \label{E:1-160}
\end{equation}%
where 
\begin{equation*}
A_{231}=N^{-1}\langle H_{+1(k+1)}S^{(k)}\psi ,H_{+(k+2)(k+3)}S^{(k)}\psi
\rangle
\end{equation*}%
\begin{equation*}
A_{232}=N^{-1}\langle \lbrack S_{1},H_{+1(k+1)}]S_{2}...S_{k}\psi
,H_{+(k+2)(k+3)}S^{(k)}\psi \rangle
\end{equation*}%
By Lemmas \ref{L:1} and \ref{L:commuting-positive-ops}, 
\begin{equation}
A_{231}\geqslant 0,  \label{E:1-161}
\end{equation}%
so we discard it. For $A_{232}$, we plug in $[S_{1},H_{+1(k+1)}]=N^{\beta
}(V^{\prime })_{N1(k+1)}$, the definition \eqref{E:1-125} of $%
H_{+(k+2)(k+3)} $, integrate by parts and use symmetry to obtain 
\begin{eqnarray*}
A_{232} &=&2N^{\beta -1}\langle (V^{\prime
})_{N1(k+1)}S_{2}...S_{k}S_{k+2}\psi ,S_{k+2}S^{(k)}\psi \rangle \\
&&+2\alpha N^{\beta -1}\langle (V^{\prime })_{N1(k+1)}S_{2}...S_{k}\psi
,S^{(k)}\psi \rangle \\
&&+N^{\beta -1}\langle (V^{\prime })_{N1(k+1)}S_{2}...S_{k}\psi
,V_{N(k+2)(k+3)}S^{(k)}\psi \rangle
\end{eqnarray*}%
For the first two terms, we apply H\"{o}lder in $x_{k+1}$, and for the third
term, we apply H\"{o}lder in both $x_{k+1}$ and $x_{k+2}$ to obtain 
\begin{eqnarray*}
\left\vert A_{232}\right\vert &\leqslant &CN^{\beta -1}\left\Vert V^{\prime
}\right\Vert _{L_{x_{k+1}}^{1}}\left\Vert S_{2}...S_{k}S_{k+2}\psi
\right\Vert _{L_{c}^{2}L_{x_{k+1}}^{\infty }}\left\Vert S_{k+2}S^{(k)}\psi
\right\Vert _{L_{c}^{2}L_{x_{k+1}}^{\infty }} \\
&&+CN^{\beta -1}\left\Vert V^{\prime }\right\Vert
_{L_{x_{k+1}}^{1}}\left\Vert S_{2}...S_{k}\psi \right\Vert
_{L_{c}^{2}L_{x_{k+1}}^{\infty }}\left\Vert S^{(k)}\psi \right\Vert
_{L_{c}^{2}L_{x_{k+1}}^{\infty }} \\
&&+CN^{\beta -1}\left\Vert V^{\prime }\right\Vert
_{L_{x_{k+1}}^{1}}\left\Vert V\right\Vert _{L_{x_{k+2}}^{1}}\left\Vert
S_{2}...S_{k}\psi \right\Vert _{L_{c}^{2}L_{x_{k+1}}^{\infty
}L_{x_{k+2}}^{\infty }}\left\Vert S^{(k)}\psi \right\Vert
_{L_{c}^{2}L_{x_{k+1}}^{\infty }L_{x_{k+2}}^{\infty }}.
\end{eqnarray*}%
Again, use (\ref{Estimate:1DSoblev}), 
\begin{eqnarray}
\left\vert A_{232}\right\vert &\leqslant &CN^{\beta -1}\left\Vert
S^{(k+1)}\psi \right\Vert _{L^{2}}\left\Vert S^{(k+2)}\psi \right\Vert
_{L^{2}}  \label{E:1-162} \\
&&+CN^{\beta -1}\left\Vert S^{(k)}\psi \right\Vert _{L^{2}}\left\Vert
S^{(k+1)}\psi \right\Vert _{L^{2}}  \notag \\
&&+CN^{\beta -1}\left\Vert S^{(k+1)}\psi \right\Vert _{L^{2}}\left\Vert
S^{(k+2)}\psi \right\Vert _{L^{2}}  \notag \\
&\leqslant &CN^{\beta -1}\left\Vert S^{(k+2)}\psi \right\Vert _{L^{2}}^{2} 
\notag
\end{eqnarray}%
Collecting \eqref{E:1-160}, \eqref{E:1-161}, and \eqref{E:1-162}, we obtain 
\begin{equation}
A_{23}\geq -CN^{\beta -1}\left\Vert S^{(k+2)}\psi \right\Vert _{L^{2}}^{2}.
\label{Estimate:A_23Final}
\end{equation}%
By \eqref{E:1-155}, \eqref{Estimate:A_21Final}, \eqref{Estimate:A_22Final},
and \eqref{Estimate:A_23Final}, we obtain 
\begin{equation}
A_{2}\geq -CN^{\beta -1}\left( \Vert S^{(k+2)}\psi \Vert
_{L^{2}}^{2}+N^{-1}\Vert S_{1}S^{(k+1)}\psi \Vert _{L^{2}}^{2}\right)
\label{Estimate:A_2Final}
\end{equation}%
Finally, combining \eqref{E:1-145}, \eqref{E:1-154}, and %
\eqref{Estimate:A_2Final}, we complete the proof of \eqref{E:1-104}
(assuming (\ref{estimate:target in energy estimate}) for $k$). Whence, we
have proved (\ref{estimate:target in energy estimate}) for all $k$ and
established Theorem \ref{THM:EnergyEstimate}.

\section{Proof of Compactness and Convergence\label{Section:Compactness and
convergence}}

\begin{theorem}[Compactness]
\label{Theorem:CompactnessOfBBGKY}For $T\in \left[ -\frac{\tan \omega T_{0}}{%
\omega },\frac{\tan \omega T_{0}}{\omega }\right] $, the sequence 
\begin{equation*}
\Gamma _{N}(\tau )=\left\{ u_{N}^{(k)}\right\} _{k=1}^{N}\in
\bigoplus_{k\geqslant 1}C\left( \left[ 0,T\right] ,\mathcal{L}%
_{k}^{1}\right) ,
\end{equation*}%
which satisfies the 1D BBGKY hierarchy (\ref{hierarchy:LensBBGKY}) subject
to energy condition (\ref{EnergyBound:LensBBGKY}) is compact with respect to
the product topology $\tau _{prod}$. For any limit point $\Gamma (t)=\left\{
u^{(k)}\right\} _{k=1}^{N},$ $\gamma ^{(k)}$ is a symmetric nonnegative
trace class operator with trace bounded by $1,$ and it verifies the energy
bound%
\begin{equation}
\sup_{\tau \in \left[ 0,T\right] }\limfunc{Tr}L^{(k)}u^{(k)}\left( \tau
\right) L^{(k)}\leqslant C^{k}.  \label{EnergyBound:ForLensGP}
\end{equation}
\end{theorem}

\begin{theorem}[Convergence]
\label{THM:convergence to GP}Let $\Gamma (\tau )=\left\{ u^{(k)}\right\}
_{k=1}^{\infty }$ be a limit point of $\Gamma _{N}(\tau )=\left\{
u_{N}^{(k)}\right\} _{k=1}^{N},$ the sequence in Theorem \ref%
{Theorem:CompactnessOfBBGKY}, with respect to the product topology $\tau
_{prod}$, then $\Gamma (\tau )$ is a solution to the focusing GP hierarchy (%
\ref{hierarchy:LensGP}) subject to initial data $u^{(k)}\left( 0\right)
=\left\vert \phi _{0}\right\rangle \left\langle \phi _{0}\right\vert
^{\otimes k}$ with coupling constant $b_{0}=$ $\int V\left( x\right) dx$,
which, written in integral form, is 
\begin{equation}
u^{(k)}\left( \tau \right) =U^{(k)}(\tau )u^{(k)}\left( 0\right)
+ib_{0}\sum_{j=1}^{k}\int_{0}^{\tau }U^{(k)}(\tau -s)\limfunc{Tr}%
\nolimits_{k+1}\left[ g(s)\delta \left( y_{j}-y_{k+1}\right)
,u^{(k+1)}\left( s\right) \right] ds.
\label{Lens GP hierarchy:integral form}
\end{equation}
\end{theorem}

\begin{proof}[Proof of Compactness]
By the standard diagonalization argument, it suffices to show the
compactness of $u_{N}^{(k)}$ for fixed $k$ with respect to the metric $\hat{d%
}_{k}$. By the Arzel\`{a}-Ascoli theorem, this is equivalent to the
equicontinuity of $u_{N}^{(k)}$, and by \cite[Lemma 6.2]{E-S-Y3}, this is
equivalent to the statement that for every observable $J^{(k)}$ from a dense
subset of $\mathcal{K}_{k}$ and for every $\varepsilon >0$, there exists $%
\delta (J^{(k)},\varepsilon )$ such that for all $\tau _{1},\tau _{2}\in %
\left[ 0,T\right] $ with $\left\vert \tau _{1}-\tau _{2}\right\vert
\leqslant \delta $, we have 
\begin{equation*}
\sup_{N}\left\vert \func{Tr}J^{(k)}u_{N}^{(k)}(\tau _{1})-\func{Tr}%
J^{(k)}u_{N}^{(k)}(\tau _{2})\right\vert \leqslant \varepsilon \,.
\end{equation*}%
We select the observables $J^{(k)}\in \mathcal{K}_{k}$ which satisfy 
\begin{equation*}
\left\Vert L_{i}L_{j}J^{(k)}L_{i}^{-1}L_{j}^{-1}\right\Vert _{\func{op}%
}+\left\Vert L_{i}^{-1}L_{j}^{-1}J^{(k)}L_{i}L_{j}\right\Vert _{\func{op}%
}<\infty ,
\end{equation*}%
where $L_{j}=\left( 1-\partial _{j}^{2}\right) ^{\frac{1}{2}}$. Assume $%
0\leqslant \tau _{1}\leqslant \tau _{2}\leqslant T$, hierarchy (\ref%
{hierarchy:LensBBGKY}) yields%
\begin{eqnarray*}
&&\left\vert \func{Tr}J^{(k)}u_{N}^{(k)}(\tau _{1})-\func{Tr}%
J^{(k)}u_{N}^{(k)}(\tau _{2})\right\vert \\
&\leqslant &\sum_{j=1}^{k}\int_{\tau _{1}}^{\tau _{2}}\left\vert \func{Tr}%
J^{(k)}\left[ -\partial _{j}^{2},u_{N}^{(k)}\left( s\right) \right]
\right\vert ds \\
&&+\frac{1}{N}\sum_{1\leqslant i<j\leqslant k}\int_{\tau _{1}}^{\tau
_{2}}\left\vert \func{Tr}J^{(k)}\left[ g(s)V_{N,s}\left( y_{i}-y_{j}\right)
,u_{N}^{(k)}\left( s\right) \right] \right\vert ds \\
&&+\frac{N-k}{N}\sum_{j=1}^{k}\int_{\tau _{1}}^{\tau _{2}}\left\vert \func{Tr%
}J^{(k)}\left[ g(s)V_{N,s}\left( y_{i}-y_{j}\right) ,u_{N}^{(k+1)}\left(
s\right) \right] \right\vert ds \\
&\leqslant &\sum_{j=1}^{k}\int_{\tau _{1}}^{\tau _{2}}Ids+\frac{1}{N}%
\sum_{1\leqslant i<j\leqslant k}\int_{\tau _{1}}^{\tau _{2}}IIds+\frac{N-k}{N%
}\sum_{j=1}^{k}\int_{\tau _{1}}^{\tau _{2}}IIIds.
\end{eqnarray*}%
For $I,$ we have, by (\ref{EnergyBound:LensBBGKY}), that,%
\begin{eqnarray*}
I &=&\left\vert \func{Tr}J^{(k)}L_{j}^{2}u_{N}^{(k)}\left( s\right) -\func{Tr%
}J^{(k)}u_{N}^{(k)}\left( s\right) L_{j}^{2}\right\vert \\
&=&\left\vert \func{Tr}L_{j}^{-1}J^{(k)}L_{j}L_{j}u_{N}^{(k)}\left( s\right)
L_{j}-\func{Tr}L_{j}J^{(k)}L_{j}^{-1}L_{j}u_{N}^{(k)}\left( s\right)
L_{j}\right\vert \\
&\leqslant &\left( \left\Vert L_{j}^{-1}J^{(k)}L_{j}\right\Vert _{\func{op}%
}+\left\Vert L_{j}J^{(k)}L_{j}^{-1}\right\Vert _{\func{op}}\right) \func{Tr}%
L_{j}u_{N}^{(k)}\left( s\right) L_{j} \\
&\leqslant &C_{J}.
\end{eqnarray*}

Lemma \ref{Lemma:SoblevLemma} and (\ref{EnergyBound:LensBBGKY}) will handle $%
II$ and $III.$ Write%
\begin{equation*}
W_{ij}=\left( L_{i}^{-1}L_{j}^{-1}V_{N,s}\left( y_{i}-y_{j}\right)
L_{i}^{-1}L_{j}^{-1}\right)
\end{equation*}%
which, by Lemma \ref{Lemma:SoblevLemma}, is a bounded operator with the bound%
\begin{equation*}
\left\Vert W_{ij}\right\Vert _{op}\leqslant C\left\Vert V\right\Vert
_{L^{1}},
\end{equation*}%
uniformly in $s.$ So then%
\begin{eqnarray*}
II &=&\left\vert \func{Tr}J^{(k)}g(s)V_{N,s}\left( y_{i}-y_{j}\right)
u_{N}^{(k)}\left( s\right) -\func{Tr}J^{(k)}u_{N}^{(k)}\left( s\right)
g(s)V_{N,s}\left( y_{i}-y_{j}\right) \right\vert \\
&=&\left\vert g(s)\right\vert |\func{Tr}%
L_{i}^{-1}L_{j}^{-1}J^{(k)}L_{i}L_{j}W_{ij}L_{i}L_{j}u_{N}^{(k)}\left(
s\right) L_{i}L_{j} \\
&&-\func{Tr}L_{i}L_{j}J^{(k)}L_{i}^{-1}L_{j}^{-1}L_{i}L_{j}u_{N}^{(k)}\left(
s\right) L_{i}L_{j}W_{ij}| \\
&\leqslant &C\left( \left\Vert
L_{i}^{-1}L_{j}^{-1}J^{(k)}L_{i}L_{j}\right\Vert _{op}+\left\Vert
L_{i}L_{j}J^{(k)}L_{i}^{-1}L_{j}^{-1}\right\Vert _{op}\right) \left\Vert
V\right\Vert _{L^{1}}\func{Tr}L_{i}L_{j}u_{N}^{(k)}\left( s\right) L_{i}L_{j}
\\
&\leqslant &C_{J}
\end{eqnarray*}%
and%
\begin{eqnarray*}
III &=&\left\vert \func{Tr}J^{(k)}g(s)V_{N,s}\left( y_{j}-y_{k+1}\right)
u_{N}^{(k+1)}\left( s\right) -\func{Tr}J^{(k)}u_{N}^{(k+1)}\left( s\right)
g(s)V_{N,s}\left( y_{j}-y_{k+1}\right) \right\vert \\
&=&\left\vert g(s)\right\vert |\func{Tr}%
L_{j}^{-1}J^{(k)}L_{j}W_{j(k+1)}L_{j}L_{k+1}u_{N}^{(k+1)}\left( s\right)
L_{j}L_{k+1} \\
&&-\func{Tr}L_{j}J^{(k)}L_{j}^{-1}L_{j}L_{k+1}u_{N}^{(k+1)}\left( s\right)
L_{j}L_{k+1}W_{j(k+1)}| \\
&\leqslant &C\left( \left\Vert L_{j}^{-1}J^{(k)}L_{j}\right\Vert
_{op}+\left\Vert L_{j}J^{(k)}L_{j}^{-1}\right\Vert _{op}\right) \left\Vert
V\right\Vert _{L^{1}}\func{Tr}L_{j}L_{k+1}u_{N}^{(k)}\left( s\right)
L_{j}L_{k+1} \\
&\leqslant &C_{J}.
\end{eqnarray*}%
Putting together the estimates of $I$, $II$, and $III$, we have%
\begin{equation*}
\sup_{N}\left\vert \func{Tr}J^{(k)}u_{N}^{(k)}(\tau _{1})-\func{Tr}%
J^{(k)}u_{N}^{(k)}(\tau _{2})\right\vert \leqslant C_{J}^{(k)}\left\vert
\tau _{1}-\tau _{2}\right\vert
\end{equation*}%
which is enough to end the proof of Theorem \ref{Theorem:CompactnessOfBBGKY}.
\end{proof}

\begin{proof}[Proof of Convergence]
By Theorem \ref{Theorem:CompactnessOfBBGKY}, passing to subsequences if
necessary, we have 
\begin{equation}
\lim_{N\rightarrow \infty }\sup_{\tau \in \lbrack 0,T]}\limfunc{Tr}%
J^{(k)}\left( u_{N}^{(k)}-u^{(k)}\right) =0\text{, }\forall J^{(k)}\in 
\mathcal{K}_{k}.  \label{condition:fast convergence}
\end{equation}

We test (\ref{Lens GP hierarchy:integral form}) against the observables $%
J^{(k)}$ in Theorem \ref{Theorem:CompactnessOfBBGKY}. We prove that the
limit point verifies%
\begin{equation}
\limfunc{Tr}J^{(k)}u^{(k)}\left( 0\right) =\limfunc{Tr}J^{(k)}\left\vert
\phi _{0}\right\rangle \left\langle \phi _{0}\right\vert ^{\otimes k},
\label{limit:testing initial data}
\end{equation}%
and%
\begin{eqnarray}
\limfunc{Tr}J^{(k)}u^{(k)} &=&\limfunc{Tr}J^{(k)}U^{(k)}(\tau )u^{(k)}\left(
0\right)  \label{limit:testing the limit pt.} \\
&&+ib_{0}\sum_{j=1}^{k}\int_{0}^{\tau }\limfunc{Tr}J^{(k)}U^{(k)}(\tau -s)%
\left[ g(s)\delta \left( y_{j}-y_{k+1}\right) ,u^{(k+1)}\left( s\right) %
\right] ds.  \notag
\end{eqnarray}%
We use the BBGKY hierarchy (\ref{hierarchy:LensBBGKY}) for this purpose.
Written in the form we need here, it becomes%
\begin{eqnarray*}
\limfunc{Tr}J^{(k)}u_{N}^{(k)} &=&\limfunc{Tr}J^{(k)}U^{(k)}(\tau
)u_{N}^{(k)}\left( 0\right) \\
&&+\frac{i}{N}\sum_{1\leqslant i<j\leqslant k}\int_{0}^{\tau }\limfunc{Tr}%
J^{(k)}U^{(k)}(\tau -s)\left[ -g(s)V_{N,s}(y_{i}-y_{j}),u_{N}^{(k)}(s)\right]
ds \\
&&+i\frac{N-k}{N}\sum_{j=1}^{k}\int_{0}^{\tau }\limfunc{Tr}%
J^{(k)}U^{(k)}(\tau -s)\left[ -g(s)V_{N,s}(y_{j}-y_{k+1}),u_{N}^{(k+1)}(s)%
\right] ds \\
&=&A+\frac{i}{N}\sum_{1\leqslant i<j\leqslant k}B+i\left( 1-\frac{k}{N}%
\right) \sum_{j=1}^{k}D.
\end{eqnarray*}%
We put a minus sign in front of $V_{N,s}$ so that the above takes the same
form as (\ref{limit:testing the limit pt.}) because $b_{0}=-\int
V_{N,s}(x)dx.$

First of all, (\ref{condition:fast convergence}) yields%
\begin{eqnarray*}
\lim_{N\rightarrow \infty }\limfunc{Tr}J^{(k)}u_{N}^{(k)} &=&\limfunc{Tr}%
J^{(k)}u^{(k)} \\
\lim_{N\rightarrow \infty }\limfunc{Tr}J^{(k)}U^{(k)}(\tau
)u_{N}^{(k)}\left( 0\right) &=&\limfunc{Tr}J^{(k)}U^{(k)}(\tau
)u^{(k)}\left( 0\right) .
\end{eqnarray*}%
Since%
\begin{equation*}
u_{N}^{(1)}\left( 0\right) =\gamma _{N}^{(1)}\left( 0\right) \rightarrow
\left\vert \phi _{0}\right\rangle \left\langle \phi _{0}\right\vert \text{
strongly as trace operators,}
\end{equation*}%
we obtain through the argument on \cite[p.64]{Lieb2} that%
\begin{equation*}
u_{N}^{(k)}\left( 0\right) =\gamma _{N}^{(k)}\left( 0\right) \rightarrow
\left\vert \phi _{0}\right\rangle \left\langle \phi _{0}\right\vert
^{\otimes k}\text{ strongly as trace operators.}
\end{equation*}%
So far, we have checked relation (\ref{limit:testing initial data}) and the
left hand side and the first term on the right hand side of (\ref%
{limit:testing the limit pt.}) for the limit point. We will prove%
\begin{equation*}
\lim_{N\rightarrow \infty }\frac{B}{N}=\lim_{N\rightarrow \infty }\frac{k}{N}%
D=0,
\end{equation*}%
\begin{equation}
\lim_{N\rightarrow \infty }D=\int_{0}^{\tau }g(s)\limfunc{Tr}%
J^{(k)}U^{(k)}(\tau -s)\left[ \delta \left( y_{j}-y_{k+1}\right)
,u^{(k+1)}\left( s\right) \right] ds.
\label{limit:emergence of the delta function}
\end{equation}%
A computation similar to the estimate of $II$ and $III$ in the proof of
Theorem \ref{Theorem:CompactnessOfBBGKY} shows that $\left\vert B\right\vert 
$ and $\left\vert D\right\vert $ are uniformly bounded for every finite
time, thus%
\begin{equation*}
\lim_{N\rightarrow \infty }\frac{B}{N}=\lim_{N\rightarrow \infty }\frac{k}{N}%
D=0.
\end{equation*}%
To acquire limit (\ref{limit:emergence of the delta function}), we use Lemma %
\ref{Lemma:ComparingDeltaFunctions}. Take a probability measure $\rho \in
L^{1}\left( \mathbb{R}\right) $, define $\rho _{\alpha }\left( y\right) =%
\frac{1}{\alpha }\rho \left( \frac{y}{\alpha }\right) .$ Use the short
notation $J_{s-\tau }^{(k)}=J^{(k)}U^{(k)}\left( \tau -s\right) $, we have%
\begin{align*}
\hspace{0.3in}& \hspace{-0.3in}\left\vert \limfunc{Tr}J^{(k)}U^{(k)}\left(
\tau -s\right) \left( -V_{N,s}\left( y_{j}-y_{k+1}\right)
u_{N}^{(k+1)}\left( s\right) -b_{0}\delta \left( y_{j}-y_{k+1}\right)
u^{(k+1)}\left( s\right) \right) \right\vert \\
& \leqslant \left\vert \limfunc{Tr}J_{s-\tau }^{(k)}\left( -V_{N,s}\left(
y_{j}-y_{k+1}\right) -b_{0}\delta \left( y_{j}-y_{k+1}\right) \right)
u_{N}^{(k+1)}\left( s\right) \right\vert \\
& \quad +b_{0}\left\vert \limfunc{Tr}J_{s-\tau }^{(k)}\left( \delta \left(
y_{j}-y_{k+1}\right) -\rho _{\alpha }\left( y_{j}-y_{k+1}\right) \right)
u_{N}^{(k+1)}\left( s\right) \right\vert \\
& \quad +b_{0}\left\vert \limfunc{Tr}J_{s-\tau }^{(k)}\rho _{\alpha }\left(
y_{j}-y_{k+1}\right) \left( u_{N}^{(k+1)}\left( s\right) -u^{(k+1)}\left(
s\right) \right) \right\vert \\
& \quad +b_{0}\left\vert \limfunc{Tr}J_{s-\tau }^{(k)}\left( \rho _{\alpha
}\left( y_{j}-y_{k+1}\right) -\delta \left( y_{j}-y_{k+1}\right) \right)
u^{(k+1)}\left( s\right) \right\vert \\
& =E+F+G+H.
\end{align*}%
A direct application of Lemma \ref{Lemma:ComparingDeltaFunctions} and the
energy condition (\ref{EnergyBound:LensBBGKY}) hands us%
\begin{eqnarray*}
E &\leqslant &\frac{C}{N^{\kappa \beta }\left( g(s)\right) ^{\kappa }}\left(
\left\Vert L_{j}^{-1}J^{(k)}L_{j}\right\Vert _{op}+\left\Vert
L_{j}J^{(k)}L_{j}^{-1}\right\Vert _{op}\right) \limfunc{Tr}%
L_{j}L_{k+1}u_{N}^{(k+1)}L_{j}L_{k+1} \\
&\leqslant &\frac{C_{J}}{N^{\kappa \beta }}\rightarrow 0\text{ as }%
N\rightarrow \infty \text{, uniformly for }s\in \lbrack 0,T]\text{ with }%
T<\infty .
\end{eqnarray*}%
Similarly, using Lemma \ref{Lemma:ComparingDeltaFunctions} and (\ref%
{EnergyBound:LensBBGKY}) and (\ref{EnergyBound:ForLensGP}), we have%
\begin{eqnarray*}
F &\leqslant &C_{\kappa }\alpha ^{\kappa }b_{0}\left( \left\Vert
L_{j}^{-1}J^{(k)}L_{j}\right\Vert _{op}+\left\Vert
L_{j}J^{(k)}L_{j}^{-1}\right\Vert _{op}\right) \limfunc{Tr}%
L_{j}L_{k+1}u_{N}^{(k+1)}L_{j}L_{k+1} \\
&\leqslant &C_{J}\alpha ^{\kappa }\rightarrow 0\text{ as }\alpha \rightarrow
0, \\
H &\leqslant &C_{\kappa }\alpha ^{\kappa }b_{0}\left( \left\Vert
L_{j}^{-1}J^{(k)}L_{j}\right\Vert _{op}+\left\Vert
L_{j}J^{(k)}L_{j}^{-1}\right\Vert _{op}\right) \limfunc{Tr}%
L_{j}L_{k+1}u^{(k+1)}L_{j}L_{k+1} \\
&\leqslant &C_{J}\alpha ^{\kappa }\rightarrow 0\text{ as }\alpha \rightarrow
0.
\end{eqnarray*}%
For G, 
\begin{eqnarray*}
G &\leqslant &b_{0}\left\vert \limfunc{Tr}J_{s-\tau }^{(k)}\rho _{\alpha
}\left( y_{j}-y_{k+1}\right) \frac{1}{1+\varepsilon L_{k+1}}\left(
u_{N}^{(k+1)}\left( s\right) -u^{(k+1)}\left( s\right) \right) \right\vert \\
&&+b_{0}\left\vert \limfunc{Tr}J_{s-\tau }^{(k)}\rho _{\alpha }\left(
y_{j}-y_{k+1}\right) \frac{\varepsilon L_{k+1}}{1+\varepsilon L_{k+1}}\left(
u_{N}^{(k+1)}\left( s\right) -u^{(k+1)}\left( s\right) \right) \right\vert .
\end{eqnarray*}%
The first term in the above estimate goes to zero as $N\rightarrow \infty $
for every $\varepsilon >0$, since we have assumed 
\eqref{condition:fast
convergence} and $J_{s-\tau }^{(k)}\rho _{\alpha }\left(
y_{j}-y_{k+1}\right) \left( 1+\varepsilon L_{k+1}\right) ^{-1}$ is a compact
operator. Due to the energy bounds on $u_{N}^{(k+1)}$ and $u^{(k+1)}$, the
second term tends to zero as $\varepsilon \rightarrow 0$, uniformly in $N$.

Combining the estimates for $E-H$, we have justified limit (\ref%
{limit:emergence of the delta function}) and thus limit (\ref{limit:testing
the limit pt.}). Hence, we have finished proving Theorem \ref%
{THM:convergence to GP}.
\end{proof}

\section{Proof of the Optimal 1D Collapsing Estimate (Theorem \protect\ref%
{THM:Optimal1DCollpasing})\label{Section:Uniqueness}}

We prove the optimality in \S \ref{Section:Optimality}. It suffices to prove
Theorem \ref{THM:OptimalUniqueness of 1D GP} for $k=1.$ We aim to prove
that, for each $\varepsilon >0$ and each bump function $\theta $, 
\begin{equation*}
\Vert \theta (\tau )R_{\varepsilon }^{(1)}U^{(1)}(-\tau )B_{1,2}U^{(2)}(\tau
)\phi ^{(2)}\Vert _{L_{\tau }^{2}L_{\mathbf{y,y}^{\prime }}^{2}}\leq
C_{\varepsilon ,\theta }\Vert R_{\varepsilon }^{(2)}\phi ^{(2)}\Vert _{L_{%
\mathbf{y,y}^{\prime }}^{2}}
\end{equation*}%
which is equivalent to 
\begin{equation*}
\Vert \theta (\tau )R_{\varepsilon }^{(1)}U^{(1)}(-\tau )B_{1,2}U^{(2)}(\tau
)R_{-\varepsilon }^{(2)}\phi ^{(2)}\Vert _{L_{\tau }^{2}L_{\mathbf{y,y}%
^{\prime }}^{2}}\leq C_{\varepsilon ,\theta }\Vert \phi ^{(2)}\Vert _{L_{%
\mathbf{y,y}^{\prime }}^{2}}
\end{equation*}%
The space-time Fourier transform of the operator on the left side is
(dropping the $x_{1}^{\prime }$ variable) 
\begin{equation*}
\iint_{\xi _{2},\xi _{2}^{\prime }}\frac{\langle \xi _{1}\rangle
^{\varepsilon }\hat{\theta}\left( \eta +\left( \xi _{1}-\xi _{2}-\xi
_{2}^{\prime }\right) ^{2}-\xi _{1}^{2}+\xi _{2}^{2}-\left( \xi _{2}^{\prime
}\right) ^{2}\right) }{\langle \xi _{1}-\xi _{2}-\xi _{2}^{\prime }\rangle
^{\varepsilon }\langle \xi _{2}\rangle ^{\varepsilon }\langle \xi
_{2}^{\prime }\rangle ^{\varepsilon }}\hat{\phi}(\xi _{1}-\xi _{2}-\xi
_{2}^{\prime },\xi _{2},\xi _{2}^{\prime })d\xi _{2}d\xi _{2}^{\prime }
\end{equation*}%
where $(\eta ,\xi )$ is the space-time Fourier variable. By the usual
Cauchy-Schwarz procedure, it suffices to prove the boundedness (independent
of $\eta $, $\xi _{1}$) of 
\begin{equation*}
I(\eta ,\xi _{1})\overset{\mathrm{def}}{=}\iint_{\xi _{2},\xi _{2}^{\prime }}%
\frac{\langle \xi _{1}\rangle ^{2\varepsilon }|\hat{\theta}\left( \eta -2\xi
_{1}\left( \xi _{2}+\xi _{2}^{\prime }\right) +\left( \xi _{2}+\xi
_{2}^{\prime }\right) ^{2}+\xi _{2}^{2}-\left( \xi _{2}^{\prime }\right)
^{2}\right) |}{\langle \xi _{1}-\xi _{2}-\xi _{2}^{\prime }\rangle
^{2\varepsilon }\langle \xi _{2}\rangle ^{2\varepsilon }\langle \xi
_{2}^{\prime }\rangle ^{2\varepsilon }}d\xi _{2}d\xi _{2}^{\prime }.
\end{equation*}

Changing variables $(\xi _{2},\xi _{2}^{\prime })\mapsto (u,v)$, where%
\begin{eqnarray}
u &=&\xi _{2}+\xi _{2}^{\prime }  \label{change of variable:u-v in 3*1D} \\
v &=&\xi _{2}-\xi _{2}^{\prime }  \notag
\end{eqnarray}
we obtain 
\begin{equation*}
I(\eta ,\xi _{1})=\iint_{u,v}|\hat{\theta}(\eta -2\xi _{1}u+u^{2}+uv)|\frac{%
\langle \xi _{1}\rangle ^{2\varepsilon }}{\langle \xi _{1}-u\rangle
^{2\varepsilon }\langle u+v\rangle ^{2\varepsilon }\langle u-v\rangle
^{2\varepsilon }}dudv
\end{equation*}%
Doing the $dv$ integral first and change $v\mapsto w$, where $w=\frac{\eta }{%
u}-2\xi _{1}+u+v$, we obtain 
\begin{eqnarray*}
I(\eta ,\xi _{1}) &=&\int_{u}\frac{\langle \xi _{1}\rangle ^{2\epsilon }}{%
\langle \xi _{1}-u\rangle ^{2\epsilon }}H(\eta ,\xi _{1},u)du \\
&=&\int_{|u|<1}+\int_{|u|>1} \\
&=&I_{1}(\eta ,\xi _{1})+I_{2}(\eta ,\xi _{1})
\end{eqnarray*}%
where 
\begin{equation}
H(\eta ,\xi _{1},u)=\int_{w}\frac{|\hat{\theta}(uw)|}{\langle w-\frac{\eta }{%
u}+2\xi _{1}\rangle ^{2\varepsilon }\langle w-2u-\frac{\eta }{u}+2\xi
_{1}\rangle ^{2\varepsilon }}dw  \label{E:2}
\end{equation}%
For convenience, we introduce 
\begin{equation}
\sigma (\eta ,\xi _{1},u)\overset{\mathrm{def}}{=}\frac{\eta }{u}-2\xi _{1}
\label{E:5}
\end{equation}

\subsubsection{Treating $I$}

We first address $I_{1}(\eta ,\xi _{1})$. For $\left\vert u\right\vert
\leqslant 1$, we have by \eqref{E:2} and \eqref{E:5} that 
\begin{eqnarray*}
H(\eta ,\xi _{1},u) &\leqslant &C\int_{w}\frac{|\hat{\theta}(uw)|}{\langle
w-\sigma \rangle ^{4\varepsilon }}\,dw \\
&\leqslant &C\int_{w}\frac{|\hat{\theta}(uw)|}{\left\vert w-\sigma
\right\vert ^{4\varepsilon }}\,dw.
\end{eqnarray*}%
Change variables, we get%
\begin{eqnarray*}
&=&C\int_{w}\frac{|\hat{\theta}(w)|}{\left\vert \frac{w}{u}-\sigma
\right\vert ^{4\varepsilon }}\frac{\,dw}{|u|} \\
&=&\frac{C}{|u|^{1-4\varepsilon }}\int_{w}\frac{|\hat{\theta}(w)|}{%
\left\vert w-u\sigma \right\vert ^{4\varepsilon }}dw.
\end{eqnarray*}%
Divide the integral into two pieces, 
\begin{eqnarray*}
&=&\frac{C}{|u|^{1-4\varepsilon }}\left( \int_{\left\vert w-u\sigma
\right\vert \leqslant 1}\frac{|\hat{\theta}(w)|}{\left\vert w-u\sigma
\right\vert ^{4\varepsilon }}dw+\int_{\left\vert w-u\sigma \right\vert
\geqslant 1}\frac{|\hat{\theta}(w)|}{\left\vert w-u\sigma \right\vert
^{4\varepsilon }}dw\right) \\
&\leqslant &\frac{C}{|u|^{1-4\varepsilon }}\left( \int_{\left\vert w-u\sigma
\right\vert \leqslant 1}\frac{1}{\left\vert w-u\sigma \right\vert
^{4\varepsilon }}dw+\int_{\left\vert w-u\sigma \right\vert \geqslant 1}|\hat{%
\theta}(w)|dw\right) .
\end{eqnarray*}%
Thus%
\begin{equation*}
H(\eta ,\xi _{1},u)\lesssim \frac{1}{|u|^{1-4\varepsilon }}.
\end{equation*}%
Therefore, plugging the above into $I_{1}$, we have 
\begin{equation*}
I_{1}(\eta ,\xi _{1})\lesssim \int_{|u|\leq 1}\frac{\langle \xi _{1}\rangle
^{2\varepsilon }}{\langle \xi _{1}-u\rangle ^{2\varepsilon
}|u|^{1-4\varepsilon }}\,du.
\end{equation*}%
Since $\left\vert u\right\vert \leqslant 1$, we have $\frac{\langle \xi
_{1}\rangle ^{2\varepsilon }}{\langle u-\xi _{1}\rangle ^{2\varepsilon }}%
\sim 1$ and therefore 
\begin{equation*}
I_{1}(\eta ,\xi _{1})\lesssim \int_{|u|\leqslant 1}\frac{du}{%
|u|^{1-4\varepsilon }}\lesssim 1.
\end{equation*}

\subsubsection{Treating $II$}

We turn our attention to $I_{2}(\eta ,\xi _{1})$. For $\left\vert
u\right\vert \geqslant 1$, by \eqref{E:2} and \eqref{E:5}, 
\begin{equation}
H(\eta ,\xi _{1},u)\lesssim \frac{1}{|u|\langle \sigma \rangle
^{2\varepsilon }\langle \sigma -2u\rangle ^{2\varepsilon }}  \label{E:6}
\end{equation}%
Indeed, in this case, the integral in \eqref{E:2} is effectively constrained
to the small interval $\left\vert w\right\vert \lesssim \left\vert
u\right\vert ^{-1}\leq 1$, and the extra factors $\langle \sigma \rangle
^{2\varepsilon }\langle \sigma -2u\rangle ^{2\varepsilon }$ in the
denominator in \eqref{E:6} come from the factors $\langle w-\sigma \rangle
^{2\varepsilon }\langle w-2u-\sigma \rangle ^{2\varepsilon }$ in the
denominator in \eqref{E:2}. Plugging \eqref{E:6} into $I_{2}(\eta ,\xi _{1})$%
, we get 
\begin{equation*}
I_{2}(\eta ,\xi _{1})\lesssim \int_{|u|\geq 1}\frac{\langle \xi _{1}\rangle
^{2\varepsilon }}{\langle \xi _{1}-u\rangle ^{2\varepsilon }|u|\langle
\sigma \rangle ^{2\varepsilon }\langle \sigma -2u\rangle ^{2\varepsilon }}%
\,du
\end{equation*}%
If $|\xi _{1}|\leq 1$, then $I_{2}(\eta ,\xi _{1})\lesssim \int_{|u|\geq 1}%
\frac{du}{|u|^{1+2\varepsilon }}$ (by neglecting the two terms $\langle
\sigma \rangle ^{2\varepsilon }\langle \sigma -2u\rangle ^{2\varepsilon }$
in the denominator), and this integral converges. If $|\xi _{1}|\geq 1$,
then $\langle \xi _{1}\rangle \sim |\xi _{1}|$ and hence 
\begin{equation*}
I_{2}(\eta ,\xi _{1})\lesssim \int_{|u|\geq 1}\frac{|\xi _{1}|^{2\varepsilon
}}{|\xi _{1}-u|^{2\varepsilon }|u||\sigma |^{2\varepsilon }|\sigma
-2u|^{2\varepsilon }}\,du
\end{equation*}%
Substituting \eqref{E:5}, 
\begin{eqnarray*}
I_{2}(\eta ,\xi _{1}) &\lesssim &\int_{|u|\geq 1}\frac{|\xi
_{1}|^{2\varepsilon }}{|\xi _{1}-u|^{2\varepsilon }|u||\frac{\eta }{u}-2\xi
_{1}|^{2\varepsilon }|u+\xi _{1}-\frac{\eta }{2u}|^{2\varepsilon }}\,du \\
&=&\int_{|u|\geq 1}\frac{|\xi _{1}|^{2\varepsilon }}{|\xi
_{1}-u|^{2\varepsilon }|u|^{1-4\varepsilon }|\eta -2\xi _{1}u|^{2\varepsilon
}|u^{2}+\xi _{1}u-\frac{\eta }{2}|^{2\varepsilon }}\,du
\end{eqnarray*}%
If $\frac{\eta }{2}+\frac{\xi _{1}^{2}}{4}\geq 0$, then let $a$, $b$ be the
roots of the quadratic $u^{2}+\xi _{1}u-\frac{\eta }{2}$ (which are real).
If $\frac{\eta }{2}+\frac{\xi _{1}^{2}}{4}<0$, then let $a=b=-\frac{\xi _{1}%
}{2}$. Then we obtain 
\begin{equation*}
I_{2}(\eta ,\xi _{1})\lesssim \int_{u}\frac{du}{|u-\xi _{1}|^{2\varepsilon
}\langle u\rangle ^{1-4\varepsilon }|u-\frac{\eta }{2\xi _{1}}%
|^{2\varepsilon }|u-a|^{2\varepsilon }|u-b|^{2\varepsilon }}
\end{equation*}%
The fact that this is bounded uniformly in $\xi _{1}$ and $\eta $ follows
from Lemma \ref{Lemma:EstimateOfIntegral}.

\begin{lemma}
\label{Lemma:EstimateOfIntegral}Suppose that $0<\varepsilon <\frac{1}{8}$.
Then 
\begin{equation*}
\int \frac{du}{|u-a|^{2\varepsilon }|u-b|^{2\varepsilon }|u-c|^{2\varepsilon
}|u-d|^{2\varepsilon }\langle u\rangle ^{1-4\varepsilon }}
\end{equation*}%
is bounded independently of $a$, $b$, $c$, $d$.
\end{lemma}

\begin{proof}
Call the given integral $G(a,b,c,d)$. Let 
\begin{equation}
F(e)\overset{\mathrm{def}}{=}\int_{u=-\infty }^{+\infty }\frac{du}{%
|u-e|^{8\varepsilon }\langle u\rangle ^{1-4\varepsilon }}  \label{E:8}
\end{equation}%
We claim that 
\begin{equation}
G(a,b,c,d)\leq F(a)+F(b)+F(c)+F(d)  \label{E:1}
\end{equation}%
To show \eqref{E:1}, we might as well assume that 
\begin{equation*}
-\infty <a\leq b\leq c\leq d<+\infty
\end{equation*}%
Divide the $u$-integration into the four intervals $-\infty <u\leq \frac{a+b%
}{2}$, $\frac{a+b}{2}\leq u\leq \frac{b+c}{2}$, $\frac{b+c}{2}\leq u\leq 
\frac{c+d}{2}$, $\frac{c+d}{2}\leq u<+\infty $. For $-\infty <u\leq \frac{a+b%
}{2}$, it is evident that the integral is bounded by $F(a)$. For $\frac{a+b}{%
2}\leq u\leq \frac{b+c}{2}$, the integral is bounded by $F(b)$, etc.

Hence it suffices to show that $F(e)$ is bounded independently of $e$. If $%
\left\vert e\right\vert \geq 1$, then we use that 
\begin{equation*}
F(e)\leq \int_{u}\frac{du}{|u-e|^{8\varepsilon }|u|^{1-4\varepsilon }}
\end{equation*}%
and then change variables $u\mapsto x$ where $u=ex$ to obtain 
\begin{equation*}
F(e)\leq \frac{1}{|e|^{4\varepsilon }}\int \frac{dx}{|x-1|^{8\varepsilon
}|x|^{1-4\varepsilon }}\lesssim 1
\end{equation*}%
If $\left\vert e\right\vert \leq 1$, then dividing the integration in $u$ in %
\eqref{E:8} into $\left\vert u\right\vert \leq 1$ and $\left\vert
u\right\vert \geq 1$ gives two integrals individually bounded independently
of $e$.
\end{proof}

\subsection{Proof of Optimality\label{Section:Optimality}}

We prove the failure of Theorem \ref{THM:Optimal1DCollpasing} for the $%
T=\infty $ and $\varepsilon \geqslant 0$ case and the $T<\infty $ and $%
\varepsilon =0$ case separately. We remark that both cases deduce to the
fact that $\int_{\left\vert u\right\vert \leqslant 1}\frac{1}{\left\vert
u\right\vert }du=\infty .$

\subsubsection{The $T=\infty $ and $\protect\varepsilon \geqslant 0$ Case}

We disprove the estimate:%
\begin{equation}
\Vert R_{\varepsilon }^{(1)}B_{1,2}U^{(2)}(\tau )R_{-\varepsilon }^{(2)}\phi
^{(2)}\Vert _{L_{\tau }^{2}L_{\mathbf{y,y}^{\prime }}^{2}}\leqslant C\Vert
\phi ^{(2)}\Vert _{L_{\mathbf{y,y}^{\prime }}^{2}}.
\label{Estimate:3*1d global in time}
\end{equation}%
By duality, it is equivalent to the estimate that 
\begin{equation}
\left\vert \int_{\mathbb{R}^{1+1}}J\left( \eta ,\xi _{1}\right) g(\eta ,\xi
_{1})d\eta d\xi _{1}\right\vert \leqslant C\Vert \phi ^{(2)}\Vert _{L_{%
\mathbf{y}}^{2}}\Vert g\Vert _{L_{\eta }^{2}L_{\xi _{1}}^{2}}
\label{Estimate:3*1d global in time (dual)}
\end{equation}%
for all $g\in L_{\eta }^{2}L_{\xi _{1}}^{2}$, where $J\left( \eta ,\xi
_{1}\right) $ is the space-time Fourier transform of $R_{\varepsilon
}^{(1)}B_{1,2}U^{(2)}(\tau )R_{-\varepsilon }^{(2)}\phi ^{(2)}$ (dropping
the $x_{1}^{\prime }$ variable) which is%
\begin{equation*}
J\left( \eta ,\xi _{1}\right) =\int \frac{\langle \xi _{1}\rangle
^{\varepsilon }\delta \left( \eta +\left( \xi _{1}-\xi _{2}-\xi _{2}^{\prime
}\right) ^{2}+\xi _{2}^{2}-\left( \xi _{2}^{\prime }\right) ^{2}\right) }{%
\langle \xi _{1}-\xi _{2}-\xi _{2}^{\prime }\rangle ^{\varepsilon }\langle
\xi _{2}\rangle ^{\varepsilon }\langle \xi _{2}^{\prime }\rangle
^{\varepsilon }}\hat{\phi}(\xi _{1}-\xi _{2}-\xi _{2}^{\prime },\xi _{2},\xi
_{2}^{\prime })d\xi _{2}d\xi _{2}^{\prime }
\end{equation*}%
Write out the left hand side of (\ref{Estimate:3*1d global in time (dual)}).%
\begin{eqnarray*}
&&\int_{\mathbb{R}^{1+1}}J\left( \eta ,\xi _{1}\right) g(\eta ,\xi
_{1})d\eta d\xi _{1} \\
&=&\int d\xi _{1}d\xi _{2}d\xi _{2}^{\prime }\hat{\phi}(\xi _{1},\xi
_{2},\xi _{2}^{\prime }) \\
&&\times \left( \int d\eta \frac{\langle \xi _{1}+\xi _{2}+\xi _{2}^{\prime
}\rangle ^{\varepsilon }\delta \left( \eta +\xi _{1}^{2}+\xi _{2}^{2}-\left(
\xi _{2}^{\prime }\right) ^{2}\right) }{\langle \xi _{1}\rangle
^{\varepsilon }\langle \xi _{2}\rangle ^{\varepsilon }\langle \xi
_{2}^{\prime }\rangle ^{\varepsilon }}g(\eta ,\xi _{1}+\xi _{2}+\xi
_{2}^{\prime })\right)
\end{eqnarray*}%
Thus estimate (\ref{Estimate:3*1d global in time (dual)}) is equivalent to
the estimate%
\begin{equation*}
\int \frac{\langle \xi _{1}+\xi _{2}+\xi _{2}^{\prime }\rangle ^{2\epsilon }%
}{\langle \xi _{1}\rangle ^{2\varepsilon }\langle \xi _{2}\rangle
^{2\varepsilon }\langle \xi _{2}^{\prime }\rangle ^{2\varepsilon }}%
\left\vert g(-\xi _{1}^{2}-\xi _{2}^{2}+\left( \xi _{2}^{\prime }\right)
^{2},\xi _{1}+\xi _{2}+\xi _{2}^{\prime })\right\vert ^{2}d\xi _{1}d\xi
_{2}d\xi _{2}^{\prime }\leqslant C\left\Vert g\right\Vert _{L_{\tau ,\xi
_{1}}^{2}}
\end{equation*}%
Performing the change of variables in (\ref{change of variable:u-v in 3*1D})
to the left hand side we get%
\begin{eqnarray*}
&&\int \frac{\langle \xi _{1}+u\rangle ^{2\epsilon }}{\langle \xi
_{1}\rangle ^{2\varepsilon }\langle u+v\rangle ^{2\varepsilon }\langle
u-v\rangle ^{2\varepsilon }}\left\vert g(-\xi _{1}^{2}-uv,\xi
_{1}+u)\right\vert ^{2}d\xi _{1}dudv \\
&=&\int \left( \int \frac{\langle \xi _{1}\rangle ^{2\epsilon }}{\langle \xi
_{1}-u\rangle ^{2\varepsilon }\langle u-\frac{\eta +\left( \xi _{1}-u\right)
^{2}}{u}\rangle ^{2\varepsilon }\langle u+\frac{\eta +\left( \xi
_{1}-u\right) ^{2}}{u}\rangle ^{2\varepsilon }}\frac{1}{\left\vert
u\right\vert }du\right) \left\vert g(\eta ,\xi _{1})\right\vert ^{2}d\xi
_{1}d\eta .
\end{eqnarray*}%
Over the region $\left\vert \eta \right\vert \lesssim 1,$ $\left\vert \xi
_{1}\right\vert \lesssim 1,$ the $du$ integral effectively becomes%
\begin{equation*}
\int \frac{1}{\langle u\rangle ^{4\varepsilon }}\frac{1}{\left\vert
u\right\vert }du
\end{equation*}%
which diverges to $\infty $. Whence, we have disproved estimate (\ref%
{Estimate:3*1d global in time}).

\subsubsection{The $T<\infty $ and $\protect\varepsilon =0$ Case}

Here, we disprove the following estimate:%
\begin{equation}
\Vert \theta (\tau )B_{1,2}U^{(2)}(\tau )\phi ^{(2)}\Vert _{L_{\tau }^{2}L_{%
\mathbf{y,y}^{\prime }}^{2}}\leqslant C\Vert \phi ^{(2)}\Vert _{L_{\mathbf{%
y,y}^{\prime }}^{2}}.  \label{Estimate:3*1d no derivative}
\end{equation}%
We proceed as in the $T=\infty $ and $\varepsilon >$ $0$ case. This time 
\begin{equation*}
J\left( \eta ,\xi _{1}\right) =\int \hat{\theta}\left( \eta +\left( \xi
_{1}-\xi _{2}-\xi _{2}^{\prime }\right) ^{2}+\xi _{2}^{2}-\left( \xi
_{2}^{\prime }\right) ^{2}\right) \hat{\phi}(\xi _{1}-\xi _{2}-\xi
_{2}^{\prime },\xi _{2},\xi _{2}^{\prime })d\xi _{2}d\xi _{2}^{\prime }
\end{equation*}%
and hence (\ref{Estimate:3*1d no derivative}) is equivalent to the estimate
that 
\begin{equation*}
\int \left\vert \left( \hat{\theta}\ast g\right) (-\xi _{1}^{2}-\xi
_{2}^{2}+\left( \xi _{2}^{\prime }\right) ^{2},\xi _{1}+\xi _{2}+\xi
_{2}^{\prime })\right\vert ^{2}d\xi _{1}d\xi _{2}d\xi _{2}^{\prime
}\leqslant C\left\Vert g\right\Vert _{L_{\tau ,\xi _{1}}^{2}}
\end{equation*}%
for all $g\in L_{\eta }^{2}L_{\xi _{1}}^{2}$. By the change of variables (%
\ref{change of variable:u-v in 3*1D}), the left side of the above estimate is%
\begin{eqnarray*}
&&\int \left\vert \left( \hat{\theta}\ast g\right) (-\xi _{1}^{2}-uv,\xi
_{1}+u)\right\vert ^{2}d\xi _{1}dudv \\
&=&\int \left( \int \frac{1}{\left\vert u\right\vert }du\right) \left\vert
\left( \hat{\theta}\ast g\right) (\eta ,\xi _{1})\right\vert ^{2}d\xi
_{1}d\eta .
\end{eqnarray*}%
This disproves estimate (\ref{Estimate:3*1d no derivative}). Together with
the $T=\infty $ and $\varepsilon \geqslant 0$ case, we have attained the
optimality of Theorem \ref{THM:Optimal1DCollpasing}.

\appendix

\section{Basic Operator Facts}

\begin{lemma}[{\protect\cite[Lemma A.1]{Kirpatrick}}]
\label{Lemma:SoblevLemma}Let $x_{i},x_{j}\in \mathbb{R}$, 
\begin{equation*}
\left\Vert L_{i}^{-1}L_{j}^{-1}V(x_{i}-x_{j})L_{i}^{-1}L_{j}^{-1}\right\Vert
_{op}\leqslant \left\Vert V\right\Vert _{L^{1}}
\end{equation*}
\end{lemma}

\begin{lemma}
\label{L:commuting-positive-ops} If $A_{2}\geqslant A_{1}\geqslant 0,$ $%
B_{2}\geqslant B_{1}\geqslant 0,$ and $[A_{i},B_{j}]=0$, $\forall i,j=1,2$,
i.e. all $A$-$B$ pairs commute. Then $A_{2}B_{2}\geqslant A_{1}B_{1}$.
\end{lemma}

\begin{proof}
We compute directly that%
\begin{equation*}
\left\langle u,A_{1}B_{1}u\right\rangle =\left\langle B_{1}^{\frac{1}{2}%
}u,A_{1}B_{1}^{\frac{1}{2}}u\right\rangle \leqslant \left\langle B_{1}^{%
\frac{1}{2}}u,A_{2}B_{1}^{\frac{1}{2}}u\right\rangle =\left\langle A_{2}^{%
\frac{1}{2}}u,B_{1}A_{2}^{\frac{1}{2}}u\right\rangle \leqslant \left\langle
u,A_{2}B_{2}u\right\rangle .
\end{equation*}
\end{proof}

\begin{lemma}
\label{Lemma:ComparingDeltaFunctions} Let $\rho \in L^{1}\left( \mathbb{R}%
\right) $ such that $\int_{\mathbb{R}}\rho \left( x\right) dx=1$ and $\int_{%
\mathbb{R}}\left\langle x\right\rangle \left\vert \rho \left( x\right)
\right\vert dx<\infty ,$ and let $\rho _{\alpha }\left( x\right) =\frac{1}{%
\alpha }\rho \left( \frac{x}{\alpha }\right) .$ Then, for every $\kappa \in
\left( 0,1\right) $, there exists $C_{\kappa }>0$ s.t.%
\begin{eqnarray*}
&&\left\vert \limfunc{Tr}J^{(k)}\left( \rho _{\alpha }\left(
x_{j}-x_{k+1}\right) -\delta \left( x_{j}-x_{k+1}\right) \right) \gamma
^{(k+1)}\right\vert \\
&\leqslant &C\left( \int \left\vert \rho \left( x\right) \right\vert
\left\vert x\right\vert ^{\kappa }dx\right) \alpha ^{\kappa }\left(
\left\Vert L_{j}J^{(k)}L_{j}^{-1}\right\Vert _{\func{op}}+\left\Vert
L_{j}^{-1}J^{(k)}L_{j}\right\Vert _{\func{op}}\right) \\
&&\times \limfunc{Tr}L_{j}L_{k+1}\gamma ^{(k+1)}L_{j}L_{k+1}
\end{eqnarray*}%
for all nonnegative $\gamma ^{(k+1)}\in \mathcal{L}_{k+1}^{1}.$
\end{lemma}

\begin{proof}
Kirkpatrick, Schlein, and Staffilani stated a similar lemma (\cite[Lemma A.2]%
{Kirpatrick}) with $\rho \geqslant 0$. Their proof, slightly modified, gives
Lemma \ref{Lemma:ComparingDeltaFunctions}. For completeness, we include the
details. It suffices to prove the estimate for $k=1$. We represent $\gamma
^{(2)}$ by $\gamma ^{(2)}=\sum_{j}\lambda _{j}\left\vert \varphi
_{j}\right\rangle \left\langle \varphi _{j}\right\vert $, where $\varphi
_{j}\in L^{2}\left( \mathbb{R}\right) $ and $\lambda _{j}\geqslant 0.$ We
write 
\begin{align*}
\hspace{0.3in}& \hspace{-0.3in}\limfunc{Tr}J^{(1)}\left( \rho _{\alpha
}\left( x_{1}-x_{2}\right) -\delta \left( x_{1}-x_{2}\right) \right) \gamma
^{(2)} \\
& =\sum_{j}\lambda _{j}\left\langle \varphi _{j},J^{(1)}\left( \rho _{\alpha
}\left( x_{1}-x_{2}\right) -\delta \left( x_{1}-x_{2}\right) \right) \varphi
_{j}\right\rangle \\
& =\sum_{j}\lambda _{j}\left\langle \psi _{j},\left( \rho _{\alpha }\left(
x_{1}-x_{2}\right) -\delta \left( x_{1}-x_{2}\right) \right) \varphi
_{j}\right\rangle
\end{align*}%
where $\psi _{j}=\left( J^{(1)}\otimes 1\right) \varphi _{j}$. By Parseval,
we find%
\begin{eqnarray*}
&&|\langle \psi _{j},(\rho _{\alpha }(x_{1}-x_{2})-\delta
(x_{1}-x_{2}))\varphi _{j}\rangle | \\
&=&|\int \overline{\hat{\psi}}_{j}(\xi _{1},\xi _{2})\hat{\varphi}_{j}(\xi
_{1}^{\prime },\xi _{2}^{\prime })\left( \hat{\rho}(\alpha (\xi _{1}-\xi
_{1}^{\prime }))-1\right) \\
&&\times \delta (\xi _{1}+\xi _{2}-\xi _{1}^{\prime }-\xi _{2}^{\prime
})d\xi _{1}d\xi _{2}d\xi _{1}^{\prime }d\xi _{2}^{\prime }|.
\end{eqnarray*}%
With $\int \rho =1$, we rewrite%
\begin{eqnarray*}
&=&|\int \overline{\hat{\psi}}_{j}(\xi _{1},\xi _{2})\hat{\varphi}_{j}(\xi
_{1}^{\prime },\xi _{2}^{\prime })\rho (x)(e^{i\alpha x\cdot (\xi
_{1}^{\prime }-\xi _{1})}-1) \\
&&\times \delta (\xi _{1}+\xi _{2}-\xi _{1}^{\prime }-\xi _{2}^{\prime
})dxd\xi _{1}d\xi _{2}d\xi _{1}^{\prime }d\xi _{2}^{\prime }| \\
&\leqslant &\int |\hat{\psi}_{j}(\xi _{1},\xi _{2})||\hat{\varphi}_{j}(\xi
_{1}^{\prime },\xi _{2}^{\prime })|\delta (\xi _{1}+\xi _{2}-\xi
_{1}^{\prime }-\xi _{2}^{\prime }) \\
&&\times |\int \rho (x)(e^{i\alpha x\cdot (\xi _{1}^{\prime }-\xi
_{1})}-1)dx|d\xi _{1}d\xi _{2}d\xi _{1}^{\prime }d\xi _{2}^{\prime }.
\end{eqnarray*}%
Using the inequality that $\forall \kappa \in \left( 0,1\right) $ 
\begin{eqnarray*}
\left\vert e^{i\alpha x\cdot (\xi _{1}^{\prime }-\xi _{1})}-1\right\vert
&\leqslant &\alpha ^{\kappa }\left\vert x\right\vert ^{\kappa }\left\vert
\xi _{1}-\xi _{1}^{\prime }\right\vert ^{\kappa } \\
&\leqslant &\alpha ^{\kappa }\left\vert x\right\vert ^{\kappa }\left(
\left\vert \xi _{1}\right\vert ^{\kappa }+\left\vert \xi _{1}^{\prime
}\right\vert ^{\kappa }\right) ,
\end{eqnarray*}%
we get%
\begin{eqnarray*}
&&\hspace{-0.3in}|\langle \psi _{j},(\rho _{\alpha }(x_{1}-x_{2})-\delta
(x_{1}-x_{2}))\varphi _{j}\rangle | \\
&\leqslant &\alpha ^{\kappa }(\int \left\vert \rho (x)\right\vert
|x|^{\kappa }dr) \\
&&\times \int |\xi _{1}|^{\kappa }|\hat{\psi}_{j}(\xi _{1},\xi _{2})||\hat{%
\varphi}_{j}(\xi _{1}^{\prime },\xi _{2}^{\prime })|\delta (\xi _{1}+\xi
_{2}-\xi _{1}^{\prime }-\xi _{2}^{\prime })d\xi _{1}d\xi _{2}d\xi
_{1}^{\prime }d\xi _{2}^{\prime } \\
&&+\alpha ^{\kappa }(\int \left\vert \rho (x)\right\vert |x|^{\kappa }dr) \\
&&\times \int |\xi _{1}^{\prime }|^{\kappa }|\hat{\psi}_{j}(\xi _{1},\xi
_{2})||\hat{\varphi}_{j}(\xi _{1}^{\prime },\xi _{2}^{\prime })|\delta (\xi
_{1}+\xi _{2}-\xi _{1}^{\prime }-\xi _{2}^{\prime })d\xi _{1}d\xi _{2}d\xi
_{1}^{\prime }d\xi _{2}^{\prime } \\
&=&\alpha ^{\kappa }(\int \left\vert \rho (x)\right\vert |x|^{\kappa
}dr)(I+II).
\end{eqnarray*}%
The estimates for $I$ and $II$ are similar, so we only deal with $I$
explicitly. We rewrite $I$ as 
\begin{eqnarray*}
I &=&\int \delta \left( \xi _{1}+\xi _{2}-\xi _{1}^{\prime }-\xi
_{2}^{\prime }\right) \frac{\left\langle \xi _{1}\right\rangle \left\langle
\xi _{2}\right\rangle }{\left\langle \xi _{1}^{\prime }\right\rangle
\left\langle \xi _{2}^{\prime }\right\rangle }\left\vert \hat{\psi}%
_{j}\left( \xi _{1},\xi _{2}\right) \right\vert \\
&&\times \frac{\left\langle \xi _{1}^{\prime }\right\rangle \left\langle \xi
_{2}^{\prime }\right\rangle }{\left\langle \xi _{1}\right\rangle ^{1-\kappa
}\left\langle \xi _{2}\right\rangle }\left\vert \hat{\varphi}_{j}\left( \xi
_{1}^{\prime },\xi _{2}^{\prime }\right) \right\vert d\xi _{1}d\xi _{2}d\xi
_{1}^{\prime }d\xi _{2}^{\prime }.
\end{eqnarray*}%
Apply Cauchy-Schwarz: 
\begin{eqnarray*}
&\leqslant &\varepsilon \int \delta \left( \xi _{1}+\xi _{2}-\xi
_{1}^{\prime }-\xi _{2}^{\prime }\right) \frac{\left\langle \xi
_{1}\right\rangle ^{2}\left\langle \xi _{2}\right\rangle ^{2}}{\left\langle
\xi _{1}^{\prime }\right\rangle ^{2}\left\langle \xi _{2}^{\prime
}\right\rangle ^{2}}\left\vert \hat{\psi}_{j}\left( \xi _{1},\xi _{2}\right)
\right\vert ^{2}d\xi _{1}d\xi _{2}d\xi _{1}^{\prime }d\xi _{2}^{\prime } \\
&&+\frac{1}{\varepsilon }\int \delta \left( \xi _{1}+\xi _{2}-\xi
_{1}^{\prime }-\xi _{2}^{\prime }\right) \frac{\left\langle \xi _{1}^{\prime
}\right\rangle ^{2}\left\langle \xi _{2}^{\prime }\right\rangle ^{2}}{%
\left\langle \xi _{1}\right\rangle ^{2\left( 1-\kappa \right) }\left\langle
\xi _{2}\right\rangle ^{2}}\left\vert \hat{\varphi}_{j}\left( \xi
_{1}^{\prime },\xi _{2}^{\prime }\right) \right\vert ^{2}d\xi _{1}d\xi
_{2}d\xi _{1}^{\prime }d\xi _{2}^{\prime }.
\end{eqnarray*}%
Rearrange terms:%
\begin{eqnarray*}
&=&\varepsilon \int \left\langle \xi _{1}\right\rangle ^{2}\left\langle \xi
_{2}\right\rangle ^{2}\left\vert \hat{\psi}_{j}\left( \xi _{1},\xi
_{2}\right) \right\vert ^{2}d\xi _{1}d\xi _{2}\int \frac{1}{\left\langle \xi
_{1}+\xi _{2}-\xi _{2}^{\prime }\right\rangle ^{2}\left\langle \xi
_{2}^{\prime }\right\rangle ^{2}}d\xi _{2}^{\prime } \\
&&+\frac{1}{\varepsilon }\int \left\langle \xi _{1}^{\prime }\right\rangle
^{2}\left\langle \xi _{2}^{\prime }\right\rangle ^{2}\left\vert \hat{\varphi}%
_{j}\left( \xi _{1}^{\prime },\xi _{2}^{\prime }\right) \right\vert ^{2}d\xi
_{1}^{\prime }d\xi _{2}^{\prime }\int \frac{1}{\left\langle \xi _{1}^{\prime
}+\xi _{2}^{\prime }-\xi _{2}\right\rangle ^{2\left( 1-\kappa \right)
}\left\langle \xi _{2}\right\rangle ^{2}}d\xi _{2} \\
&\leqslant &\varepsilon \left\langle \psi _{j},L_{1}^{2}L_{2}^{2}\psi
_{j}\right\rangle \sup_{\xi }\int_{\mathbb{R}}\frac{1}{\left\langle \xi
-\eta \right\rangle ^{2}\left\langle \eta \right\rangle ^{2}}d\eta \\
&&+\frac{1}{\varepsilon }\left\langle \varphi _{j},L_{1}^{2}L_{2}^{2}\varphi
_{j}\right\rangle \sup_{\xi }\int_{\mathbb{R}}\frac{1}{\left\langle \xi
-\eta \right\rangle ^{2\left( 1-\kappa \right) }\left\langle \eta
\right\rangle ^{2}}d\eta .
\end{eqnarray*}%
When $\kappa \in \left[ 0,1\right] ,$ 
\begin{eqnarray*}
\sup_{\xi }\int_{\mathbb{R}}\frac{1}{\left\langle \xi -\eta \right\rangle
^{2\left( 1-\kappa \right) }\left\langle \eta \right\rangle ^{2}}d\eta
&<&\infty , \\
\sup_{\xi }\int_{\mathbb{R}}\frac{1}{\left\langle \xi -\eta \right\rangle
^{2}\left\langle \eta \right\rangle ^{2}}d\eta &<&\infty ,
\end{eqnarray*}%
hence we have 
\begin{eqnarray*}
&&\left\vert \limfunc{Tr}J^{(1)}\left( \rho _{\alpha }\left(
x_{1}-x_{2}\right) -\delta \left( x_{1}-x_{2}\right) \right) \gamma
^{(k+1)}\right\vert \\
&\leqslant &C\left( \int \left\vert \rho \left( x\right) \right\vert
\left\vert x\right\vert ^{\kappa }dx\right) \alpha ^{\kappa }\left(
\varepsilon \limfunc{Tr}J^{(1)}L_{1}^{2}L_{2}^{2}J^{(1)}\gamma ^{(2)}+\frac{1%
}{\varepsilon }\limfunc{Tr}L_{1}^{2}L_{2}^{2}\gamma ^{(2)}\right) \\
&=&C\left( \int \left\vert \rho \left( x\right) \right\vert \left\vert
x\right\vert ^{\kappa }dx\right) \alpha ^{\kappa } \\
&&\times \left( \varepsilon \limfunc{Tr}%
L_{1}^{-1}L_{2}^{-1}J^{(1)}L_{1}L_{1}J^{(1)}L_{1}^{-1}L_{1}L_{2}^{2}\gamma
^{(2)}L_{1}L_{2}+\frac{1}{\varepsilon }\limfunc{Tr}L_{1}^{2}L_{2}^{2}\gamma
^{(2)}\right) \\
&\leqslant &C\left( \int \left\vert \rho \left( x\right) \right\vert
\left\vert x\right\vert ^{\kappa }dx\right) \alpha ^{\kappa }\left(
\varepsilon \left\Vert L_{1}^{-1}J^{(1)}L_{1}\right\Vert _{\func{op}%
}\left\Vert L_{1}J^{(1)}L_{1}^{-1}\right\Vert _{\func{op}}+\frac{1}{%
\varepsilon }\right) \\
&&\times \limfunc{Tr}L_{1}^{2}L_{2}^{2}\gamma ^{(2)}.
\end{eqnarray*}%
Let $\varepsilon =\Vert L_{1}J^{(1)}L_{1}^{-1}\Vert _{\func{op}}^{-1}$, we
reach%
\begin{eqnarray*}
&\leqslant &C\left( \int \left\vert \rho \left( x\right) \right\vert
\left\vert x\right\vert ^{\kappa }dx\right) \alpha ^{\kappa }\left(
\left\Vert L_{1}^{-1}J^{(1)}L_{1}\right\Vert _{\func{op}}+\left\Vert
L_{1}J^{(1)}L_{1}^{-1}\right\Vert _{\func{op}}\right) \\
&&\times \limfunc{Tr}L_{1}^{2}L_{2}^{2}\gamma ^{(2)}
\end{eqnarray*}%
as claimed.
\end{proof}

\section{Deducing Theorem \protect\ref{THM:MainTHMFiniteKinetic} from
Theorem \protect\ref{THM:Main Theorem}\label{Appendix:equivalence}}

If $\psi _{N}\left( 0\right) $ satisfies (a), (b), and (c) in Theorem \ref%
{THM:MainTHMFiniteKinetic}, then $\psi _{N}\left( 0\right) $ checks the
requirements of Lemma \ref{Lemma:B2}. Thus we can define an approximation $%
\psi _{N}^{\kappa }(0)$ of $\psi _{N}\left( 0\right) $ as in (\ref%
{def:smooth approximation}). Via (i) and (iii) of Lemma \ref{Lemma:B2}, $%
\psi _{N}^{\kappa }(0)$ verifies the hypothesis of Theorem \ref{THM:Main
Theorem} for small enough $\kappa >0.$ Therefore, for $\gamma _{N}^{\kappa
,(k)}\left( t\right) ,$ the marginal density associated with $%
e^{itH_{N}}\psi _{N}^{\kappa }(0)$, Theorem \ref{THM:Main Theorem} gives the
convergence 
\begin{equation}
\lim_{N\rightarrow \infty }\limfunc{Tr}\left\vert \gamma _{N}^{\kappa
,(k)}\left( t\right) (t,\mathbf{x}_{k};\mathbf{x}_{k}^{\prime
})-\dprod_{j=1}^{k}\phi (t,x_{j})\overline{\phi }(t,x_{j}^{\prime
})\right\vert =0.  \label{convergence:smooth}
\end{equation}%
for all small enough $\kappa >0,$ all $k\geqslant 1$, and all $t\in \mathbb{R%
}$.

For $\gamma _{N}^{(k)}\left( t\right) $ in Theorem \ref%
{THM:MainTHMFiniteKinetic}, we notice that, $\forall J^{(k)}\in \mathcal{K}%
_{k}$, $\forall t\in \mathbb{R}$, we have 
\begin{eqnarray*}
&&\left\vert \limfunc{Tr}J^{(k)}\left( \gamma _{N}^{(k)}\left( t\right)
-\left\vert \phi \left( t\right) \right\rangle \left\langle \phi \left(
t\right) \right\vert ^{\otimes k}\right) \right\vert \\
&\leqslant &\left\vert \limfunc{Tr}J^{(k)}\left( \gamma _{N}^{(k)}\left(
t\right) -\gamma _{N}^{\kappa ,(k)}\left( t\right) \right) \right\vert \\
&&+\left\vert \limfunc{Tr}J^{(k)}\left( \gamma _{N}^{\kappa ,(k)}\left(
t\right) -\left\vert \phi \left( t\right) \right\rangle \left\langle \phi
\left( t\right) \right\vert ^{\otimes k}\right) \right\vert \\
&=&\text{I}+\text{II}.
\end{eqnarray*}%
Convergence \eqref{convergence:smooth} then takes care of $\text{II}$. To
handle $\text{I}$ , part (ii) of Lemma \ref{Lemma:B2} yields 
\begin{equation*}
\left\Vert e^{itH_{N}}\psi _{N}^{\kappa }(0)-e^{itH_{N}}\psi
_{N}(0)\right\Vert _{L^{2}}=\left\Vert \psi _{N}^{\kappa }(0)-\psi
_{N}(0)\right\Vert _{L^{2}}\leqslant C\kappa ^{\frac{1}{2}}
\end{equation*}%
which implies 
\begin{equation*}
I=\left\vert \limfunc{Tr}J^{(k)}\left( \gamma _{N}^{(k)}\left( t\right)
-\gamma _{N}^{\kappa ,(k)}\left( t\right) \right) \right\vert \leqslant
C\left\Vert J^{(k)}\right\Vert _{op}\kappa ^{\frac{1}{2}}.
\end{equation*}%
Since $\kappa >0$ is arbitrary, we deduce that 
\begin{equation*}
\lim_{N\rightarrow \infty }\left\vert \limfunc{Tr}J^{(k)}\left( \gamma
_{N}^{(k)}\left( t\right) -\left\vert \phi \left( t\right) \right\rangle
\left\langle \phi \left( t\right) \right\vert ^{\otimes k}\right)
\right\vert =0.
\end{equation*}%
i.e. as trace class operators 
\begin{equation*}
\gamma _{N}^{(k)}\left( t\right) \rightarrow \left\vert \phi \left( t\right)
\right\rangle \left\langle \phi \left( t\right) \right\vert ^{\otimes k}%
\text{ weak*.}
\end{equation*}%
Then again, the Gr\"{u}mm's convergence theorem upgrades the above weak*
convergence to strong. Thence, we have concluded Theorem \ref%
{THM:MainTHMFiniteKinetic} via Theorem \ref{THM:Main Theorem}.

\begin{lemma}
\label{Lemma:B2}Assume $\psi _{N}\left( 0\right) $ satisfies (a), (b), and
(c) in Theorem \ref{THM:MainTHMFiniteKinetic}. Let $\chi \in C_{0}^{\infty
}\left( \mathbb{R}\right) $ be a cut-off such that $0\leqslant \chi
\leqslant 1$, $\chi \left( s\right) =1$ for $0\leqslant s\leqslant 1$ and $%
\chi \left( s\right) =0$ for $s\geqslant 2.$ For $\kappa >0,$ we define an
approximation $\psi _{N}^{\kappa }(0)$ of $\psi _{N}\left( 0\right) $ by 
\begin{equation}
\psi _{N}^{\kappa }(0)=\frac{\chi \left( \kappa H_{N}/N\right) \psi
_{N}\left( 0\right) }{\left\Vert \chi \left( \kappa H_{N}/N\right) \psi
_{N}\left( 0\right) \right\Vert }.  \label{def:smooth approximation}
\end{equation}%
This approximation has the following properties:

(i) $\psi _{N}^{\kappa }(0)$ verifies the energy condition 
\begin{equation*}
\langle \psi _{N}^{\kappa }(0),H_{N}^{k}\psi _{N}^{\kappa }(0)\rangle
\leqslant \frac{2^{k}N^{k}}{\kappa ^{k}}.
\end{equation*}

(ii) 
\begin{equation*}
\sup_{N}\left\Vert \psi _{N}^{\kappa }(0)-\psi _{N}(0)\right\Vert
_{L^{2}}\leqslant C\kappa ^{\frac{1}{2}}.
\end{equation*}

(iii) For small enough $\kappa >0$, $\psi _{N}^{\kappa }(0)$ is
asymptotically factorized as well 
\begin{equation*}
\lim_{N\rightarrow \infty }\limfunc{Tr}\left\vert \gamma _{N}^{\kappa
,(1)}(0,x_{1};x_{1}^{\prime })-\phi _{0}(x_{1})\overline{\phi _{0}}%
(x_{1}^{\prime })\right\vert =0,
\end{equation*}%
where $\gamma _{N}^{\kappa ,(1)}\left( 0\right) $ is the marginal density
associated with $\psi _{N}^{\kappa }(0)$ and $\phi _{0}$ is the same as in
assumption (b) in Theorem \ref{THM:MainTHMFiniteKinetic}.
\end{lemma}

\begin{proof}
Let us write $\chi \left( \kappa H_{N}/N\right) $ as $\chi $ and $\psi
_{N}\left( 0\right) $ as $\psi _{N}.$ This proof of Lemma \ref{Lemma:B2}
closely follows the proof of \cite[Proposition 9.1 (i) and (ii)]{E-S-Y5} and 
\cite[Proposition 5.1 (iii)]{E-S-Y2}.

(i) is from definition. In fact, denote the characteristic function of $%
\left[ 0,\lambda \right] $ with $\mathbf{1}(s\leqslant \lambda ).$ We see
that $\mathbf{1}(H_{N}\leqslant 2N/\kappa )\chi \left( \kappa H_{N}/N\right)
=\chi \left( \kappa H_{N}/N\right) .$ Thus%
\begin{eqnarray*}
\left\langle \psi _{N}^{\kappa }(0),H_{N}^{k}\psi _{N}^{\kappa
}(0)\right\rangle &=&\left\langle \frac{\chi \psi _{N}}{\left\Vert \chi \psi
_{N}\right\Vert },\mathbf{1}(H_{N}\leqslant 2N/\kappa )H_{N}^{k}\frac{\chi
\psi _{N}}{\left\Vert \chi \psi _{N}\right\Vert }\right\rangle \\
&\leqslant &\left\Vert \mathbf{1}(H_{N}\leqslant 2N/\kappa
)H_{N}^{k}\right\Vert _{op} \\
&\leqslant &\frac{2^{k}N^{k}}{\kappa ^{k}}.
\end{eqnarray*}

We prove (ii) with a slightly modified proof of \cite[Proposition 9.1 (ii)]%
{E-S-Y5}. We still have%
\begin{eqnarray*}
\left\Vert \psi _{N}^{\kappa }-\psi _{N}\right\Vert _{L^{2}} &\leqslant
&\left\Vert \chi \psi _{N}-\psi _{N}\right\Vert _{L^{2}}+\left\Vert \frac{%
\chi \psi _{N}}{\left\Vert \chi \psi _{N}\right\Vert }-\chi \psi
_{N}\right\Vert _{L^{2}} \\
&\leqslant &\left\Vert \chi \psi _{N}-\psi _{N}\right\Vert
_{L^{2}}+\left\vert 1-\left\Vert \chi \psi _{N}\right\Vert \right\vert \\
&\leqslant &2\left\Vert \chi \psi _{N}-\psi _{N}\right\Vert _{L^{2}},
\end{eqnarray*}%
where%
\begin{eqnarray*}
\left\Vert \chi \psi _{N}-\psi _{N}\right\Vert _{L^{2}}^{2} &=&\left\langle
\psi _{N},\left( 1-\chi \left( \kappa H_{N}/N\right) \right) ^{2}\psi
_{N}\right\rangle \\
&\leqslant &\left\langle \psi _{N},\mathbf{1}(\frac{\kappa H_{N}}{N}%
\geqslant 1)\psi _{N}\right\rangle .
\end{eqnarray*}%
To continue estimating, we notice that if $C\geqslant 0$, then $\mathbf{1}%
(s\geqslant 1)\leqslant \mathbf{1}(s+C\geqslant 1)$ for all $s$. So%
\begin{eqnarray*}
\left\Vert \chi \psi _{N}-\psi _{N}\right\Vert _{L^{2}}^{2} &\leqslant
&\left\langle \psi _{N},\mathbf{1}(\frac{\kappa H_{N}}{N}\geqslant 1)\psi
_{N}\right\rangle \\
&\leqslant &\left\langle \psi _{N},\mathbf{1}(\frac{\kappa \left(
H_{N}+N\alpha +N\right) }{N}\geqslant 1)\psi _{N}\right\rangle .
\end{eqnarray*}
With the inequality that $\mathbf{1}(s\geqslant 1)\leqslant s$ for all $%
s\geqslant 0$ and the fact that 
\begin{equation*}
H_{N}+N\alpha +N\geqslant 0
\end{equation*}%
proved in Lemma \ref{L:1}, we arrive at%
\begin{eqnarray*}
\left\Vert \chi \psi _{N}-\psi _{N}\right\Vert _{L^{2}}^{2} &\leqslant &%
\frac{\kappa }{N}\left\langle \psi _{N},\left( H_{N}+N\alpha +N\right) \psi
_{N}\right\rangle \\
&\leqslant &\frac{\kappa }{N}\left\langle \psi _{N},H_{N}\psi
_{N}\right\rangle +\left( 1+\alpha \right) \kappa \left\langle \psi
_{N},\psi _{N}\right\rangle ,
\end{eqnarray*}%
where%
\begin{eqnarray*}
\frac{1}{N}\left\langle \psi _{N},H_{N}\psi _{N}\right\rangle
&=&\left\langle \psi _{N},\left( -\partial _{x_{1}}^{2}+\omega
^{2}x_{1}^{2}\right) \psi _{N}\right\rangle \\
&&+\frac{1}{N^{2}}\sum_{1\leqslant i<j\leqslant N}\int N^{\beta }V(N^{\beta
}\left( x_{i}-x_{j}\right) )\left\vert \psi _{N}\right\vert ^{2}d\mathbf{x}%
_{N} \\
&\leqslant &\left\langle \psi _{N},\left( -\partial _{x_{1}}^{2}+\omega
^{2}x_{1}^{2}\right) \psi _{N}\right\rangle \\
&&+C\left\Vert V\right\Vert _{L^{1}}\int \left( \left\vert \psi
_{N}\right\vert ^{2}+\left\vert \partial _{x_{1}}\psi _{N}\right\vert
^{2}\right) d\mathbf{x}_{N} \\
&\leqslant &C\left\langle \psi _{N},\left( -\partial _{x_{1}}^{2}+\omega
^{2}x_{1}^{2}\right) \psi _{N}\right\rangle +C.
\end{eqnarray*}%
Using (a) and (c) in the assumptions of Theorem \ref%
{THM:MainTHMFiniteKinetic}, we deduce that%
\begin{equation*}
\left\Vert \chi \psi _{N}-\psi _{N}\right\Vert _{L^{2}}^{2}\leqslant C\kappa
\end{equation*}%
which implies%
\begin{equation*}
\left\Vert \psi _{N}^{\kappa }-\psi _{N}\right\Vert _{L^{2}}\leqslant
C\kappa ^{\frac{1}{2}}.
\end{equation*}

(iii) does not follow from the proof of \cite[Proposition 9.1 (iii)]{E-S-Y5}
in which the positivity of $V$ is used. (iii) follows from the proof of \cite%
[Proposition 5.1 (iii)]{E-S-Y2} which does not require $V$ to hold a
definite sign.\footnote{%
See \cite[(5.19)]{E-S-Y2}.} Notice that we are working in one dimension, we
get a $N^{\frac{\beta }{2}}$ instead of a $N^{\frac{3\beta }{2}}$ in \cite[%
(5.20)]{E-S-Y2} and hence we get a $N^{\frac{\beta }{2}-1}$ in the estimate
of \cite[(5.18)]{E-S-Y2} which goes to zero for $\beta \in \left( 0,1\right)
.$
\end{proof}


\begin{thebibliography}{99}
\bibitem{AGT} R. Adami, F. Golse, and A. Teta, \emph{Rigorous derivation of
the cubic NLS in dimension one}, J. Stat. Phys. 127 (2007), 1194--1220.

\bibitem{Ammari2} Z. Ammari and F. Nier, \emph{Mean Field Propagation of
Wigner Measures and BBGKY Hierarchies for General Bosonic States}, J. Math.
Pures. Appl. \textbf{95} (2011), 585-626.

\bibitem{Ammari1} Z. Ammari and F. Nier, \emph{Mean Field Propagation of
Infinite Dimensional Wigner Measures with a Singular Two-body Interaction
Potential, }arXiv:1111.5918.

\bibitem{Anderson} M. H. Anderson, J. R. Ensher, M. R. Matthews, C. E.
Wieman, and E. A. Cornell, \emph{Observation of Bose-Einstein Condensation
in a Dilute Atomic Vapor}, Science \textbf{269 }(1995), 198--201.

\bibitem{Beckner} W. Beckner, \emph{Multilinear Embedding -- Convolution
Estimates on Smooth Submanifolds}, to appear in Proc. Amer. Math. Soc.

\bibitem{SchleinNew} N. Benedikter, G. Oliveira, and B. Schlein, \emph{\
Quantitative Derivation of the Gross-Pitaevskii Equation, }62pp,
arXiv:1208.0373.

\bibitem{Carles} R. Carles, \emph{Nonlinear Schr\"{o}dinger Equation with
Time Dependent Potential}, Commun. Math. Sci. \textbf{9} (2011), 937-964.

\bibitem{LChen} L. Chen, J. O. Lee and B. Schlein, \emph{Rate of Convergence
Towards Hartree Dynamics}, J. Stat. Phys. \textbf{144 (}2011\textbf{)},
872--903.

\bibitem{TChenAndNpGP1} T. Chen and N. Pavlovi\'{c}, \emph{On the Cauchy
Problem for Focusing and Defocusing Gross-Pitaevskii Hierarchies}, Discrete
Contin. Dyn. Syst. \textbf{27} (2010), 715--739.

\bibitem{TChenAndNP} T. Chen and N. Pavlovi\'{c}, \emph{The Quintic NLS as
the Mean Field Limit of a Boson Gas with Three-Body Interactions}, J. Funct.
Anal. \textbf{260} (2011), 959--997.

\bibitem{TChenAndNPSpace-Time} T. Chen and N. Pavlovi\'{c}, \emph{Derivation
of the cubic NLS and Gross-Pitaevskii hierarchy from manybody dynamics in }$%
d=3$\emph{\ based on spacetime norms}, Annales Henri Poincar\'{e} 2013,
46pp. DOI: 10.1007/s00023-013-0248-6. (arXiv:1111.6222)

\bibitem{TCNPNT} T. Chen, N. Pavlovi\'{c}, and N. Tzirakis, \emph{Energy
Conservation and Blowup of Solutions for Focusing Gross--Pitaevskii
Hierarchies}, Ann. I. H. Poincar\'{e} \textbf{27} (2010), 1271-1290.

\bibitem{TCNPNT1} T. Chen, N. Pavlovi\'{c}, and N. Tzirakis, \emph{\
Multilinear Morawetz identities for the Gross-Pitaevskii hierarchy},
Contemporary Mathematics \textbf{581} (2012), 39-62.

\bibitem{ChenDie} X. Chen, \emph{Classical Proofs Of Kato Type Smoothing
Estimates for The Schr\"{o}dinger Equation with Quadratic Potential in }$%
\mathbb{R}^{n+1}$\emph{\ with Application}, Differential and Integral
Equations \textbf{24} (2011), 209-230.

\bibitem{Chen2ndOrder} X. Chen, \emph{Second Order Corrections to Mean Field
Evolution for Weakly Interacting Bosons in the Case of Three-body
Interactions\textbf{,} }Arch. Rational Mech. Anal. \textbf{203} (2012),
455-497.

\bibitem{ChenAnisotropic} X. Chen, \emph{Collapsing Estimates and the
Rigorous Derivation of the 2d Cubic Nonlinear Schr\"{o}dinger Equation with
Anisotropic Switchable Quadratic Traps}, J. Math. Pures Appl. \textbf{98}
(2012), 450--478.

\bibitem{Chen3DDerivation} X. Chen, \emph{On the Rigorous Derivation of the
3D Cubic Nonlinear Schr\"{o}dinger Equation with A Quadratic Trap}, Arch.
Rational Mech. Anal. 2013, 44pp. DOI: 10.1007/s00205-013-0645-5.
(arXiv:1204.0125)

\bibitem{C-H3Dto2D} X. Chen and J. Holmer, \emph{On the Rigorous Derivation
of the 2D Cubic Nonlinear Schr\"{o}dinger Equation from 3D Quantum Many-Body
Dynamics}, Arch. Rational Mech. Anal., 2013, 46pp. DOI:
10.1007/s00205-013-0667-z. (arXiv:1212.0787)

\bibitem{C-H2/3} X. Chen and J. Holmer, \emph{On the Klainerman-Machedon
Conjecture of the Quantum BBGKY Hierarchy with Self-interaction}, 34pp,
arXiv:1303.5385, submitted.

\bibitem{Cornish} S. L. Cornish, N. R. Claussen, J. L. Roberts, E. A.
Cornell, and C. E. Wieman, \emph{Stable }$^{85}$\emph{Rb Bose-Einstein
Condensates with Widely Turnable Interactions, }Phys. Rev. Lett. \textbf{85 }
(2000), 1795-1798.

\bibitem{Davis} K. B. Davis, M. -O. Mewes, M. R. Andrews, N. J. van Druten,
D. S. Durfee, D. M. Kurn, and W. Ketterle, \emph{Bose-Einstein condensation
in a gas of sodium atoms}, Phys. Rev. Lett. \textbf{75 }(1995), 3969--3973.

\bibitem{E-E-S-Y1} A. Elgart, L. Erd\"{o}s, B. Schlein, and H. T. Yau, \emph{%
\ Gross-Pitaevskii Equation as the Mean Field Limit of Weakly Coupled Bosons,%
} Arch. Rational Mech. Anal. \textbf{179\ }(2006), 265--283.

\bibitem{E-Y1} L. Erd\"{o}s and H. T. Yau, \emph{Derivation of the
Non-linear Schr\"{o}dinger Equation from a Many-body Coulomb System,} Adv.
Theor. Math. Phys. \textbf{5\ }(2001), 1169--1205.

\bibitem{E-S-Y1} L. Erd\"{o}s, B. Schlein, and H. T. Yau, \emph{Derivation
of the Gross-Pitaevskii Hierarchy for the Dynamics of Bose-Einstein
Condensate}, Comm. Pure Appl. Math. \textbf{59} (2006), 1659--1741.

\bibitem{E-S-Y2} L. Erd\"{o}s, B. Schlein, and H. T. Yau, \emph{Derivation
of the Cubic non-linear Schr\"{o}dinger Equation from Quantum Dynamics of
Many-body Systems}, Invent. Math. \textbf{167} (2007), 515--614.

\bibitem{E-S-Y5} L. Erd\"{o}s, B. Schlein, and H. T. Yau, \emph{Rigorous
Derivation of the Gross-Pitaevskii Equation with a Large Interaction
Potential}, J. Amer. Math. Soc. \textbf{22} (2009), 1099-1156.

\bibitem{E-S-Y3} L. Erd\"{o}s, B. Schlein, and H. T. Yau, \emph{Derivation
of the Gross-Pitaevskii Equation for the Dynamics of Bose-Einstein Condensate%
}, Annals Math. \textbf{172} (2010), 291-370.{}

\bibitem{Frolich} J. Fr\"{o}hlich, A. Knowles, and S. Schwarz, \emph{On the
Mean-Field Limit of Bosons with Coulomb Two-Body Interaction, }Commun. Math.
Phys. \textbf{288} (2009), 1023--1059.

\bibitem{Sohinger} P. Gressman, V. Sohinger, and G. Staffilani, \emph{On the
Uniqueness of Solutions to the Periodic 3D Gross-Pitaevskii Hierarchy, }%
39pp, arXiv:1212.2987.

\bibitem{GM} M. G. Grillakis and D. Margetis, \emph{A Priori Estimates for
Many-Body Hamiltonian Evolution of Interacting Boson System}, J. Hyperb.
Diff. Eqs. \textbf{5} (2008), 857-883.

\bibitem{GM1} M. G. Grillakis and M. Machedon, \emph{Pair excitations and
the mean field approximation of interacting Bosons, I, }to appear in Commun.
Math. Phys. (arXiv:1208.3763).\textbf{\ }

\bibitem{GMM1} M. G. Grillakis, M. Machedon, and D. Margetis, \emph{Second
Order Corrections to Mean Field Evolution for Weakly Interacting Bosons. I,}
Commun. Math. Phys. \textbf{294} (2010), 273-301.

\bibitem{GMM2} M. G. Grillakis, M. Machedon, and D. Margetis, \emph{Second
Order Corrections to Mean Field Evolution for Weakly Interacting Bosons. II,}
Adv. Math. \textbf{228} (2011) 1788-1815.

\bibitem{Hayashi-Naumkin} N. Hayashi and P. Naumkin, \emph{Asymptotics for
large time of solutions to the nonlinear Schr\"{o}dinger and Hartree
equations}, Amer. J. Math. \textbf{120} (1998), 369-389.

\bibitem{Ketterle} W. Ketterle and N. J. van Druten, \emph{Evaporative
Cooling of Trapped Atoms, }Advances In Atomic, Molecular, and Optical
Physics \textbf{37} (1996), 181-236.

\bibitem{Khaykovich} L. Khaykovich, F. Schreck, G. Ferrari, T. Bourdel, J.
Cubizolles, L. D. Carr, Y. Castin, C. Salomon, \emph{Formation of a
Matter-Wave Bright Soliton}, Science \textbf{296 }(2002), 1290-1293.

\bibitem{Kirpatrick} K. Kirkpatrick, B. Schlein and G. Staffilani, \emph{\
Derivation of the Two Dimensional Nonlinear Schr\"{o}dinger Equation from
Many Body Quantum Dynamics}, Amer. J. Math. \textbf{133} (2011), 91-130.

\bibitem{KlainermanAndMachedon} S. Klainerman and M. Machedon, \emph{On the
Uniqueness of Solutions to the Gross-Pitaevskii Hierarchy}, Commun. Math.
Phys. \textbf{279} (2008), 169-185.

\bibitem{KnowlesAndPickl} A. Knowles and P. Pickl, \emph{Mean-Field
Dynamics: Singular Potentials and Rate of Convergence}, Commum. Math. Phys. 
\textbf{298} (2010), 101-138.

\bibitem{Lieb2} E. H. Lieb, R. Seiringer, J. P. Solovej and J. Yngvason, 
\emph{The Mathematics of the Bose Gas and Its Condensation}, Basel,
Switzerland: Birkha\"{u}ser Verlag, 2005.

\bibitem{Medley} P. Medley, D.M. Weld, H. Miyake, D.E. Pritchard, and W.
Ketterle, \emph{Spin Gradient Demagnetization Cooling of Ultracold Atoms,}
Phys. Rev. Lett. \textbf{106} (2011), 195301.

\bibitem{MichelangeliSchlein} A. Michelangeli and B. Schlein,\emph{\
Dynamical Collapse of Boson Stars}, Commum. Math. Phys. \textbf{311} (2012),
645-687.

\bibitem{Pickl} P. Pickl, \emph{A Simple Derivation of Mean Field Limits for
Quantum Systems}, Lett. Math. Phys. \textbf{97} (2011), 151-164.

\bibitem{RodnianskiAndSchlein} I. Rodnianski and B. Schlein, \emph{Quantum
Fluctuations and Rate of Convergence Towards Mean Field Dynamics}, Commun.
Math. Phys. \textbf{291} (2009), 31-61.

\bibitem{Simon} B. Simon, \emph{Trace Ideals and Their Applications: Second
Edition}, Mathematical Surveys Monogr. \textbf{120}, Amer. Math. Soc.,
Providence, RI, 2005.

\bibitem{Spohn} H. Spohn, \emph{Kinetic Equations from Hamiltonian Dynamics}%
, Rev. Mod. Phys. \textbf{52} (1980), 569-615.

\bibitem{Stamper} D. M. Stamper-Kurn, M. R. Andrews, A. P. Chikkatur, S.
Inouye, H. -J. Miesner, J. Stenger, and W. Ketterle, \emph{Optical
Confinement of a Bose-Einstein Condensate}, Phys. Rev. Lett. \textbf{80 }
(1998), 2027-2030.

\bibitem{Streker} K.E. Strecker; G.B. Partridge; A.G. Truscott; R.G. Hulet, 
\emph{Formation and Propagation of Matter-wave Soliton Trains, }Nature 
\textbf{417} (2002), 150-153.
\end{thebibliography}
\end{document}